\documentclass[10pt,twoside,english]{extarticle}

\usepackage[T1]{fontenc}
\usepackage[utf8]{inputenc}
\usepackage{float}
\usepackage{amsmath}
\usepackage{amsthm}
\usepackage{amssymb}
\usepackage{graphicx}
\usepackage[numbers]{natbib}

\makeatletter

\providecommand{\tabularnewline}{\\}

\theoremstyle{plain}
\newtheorem{thm}{\protect\theoremname}
\theoremstyle{remark}
\newtheorem{rem}{\protect\remarkname}
\theoremstyle{plain}
\newtheorem{prop}{\protect\propositionname}

\usepackage{babel}

\makeatother

\usepackage{babel}
\providecommand{\propositionname}{Proposition}
\providecommand{\remarkname}{Remark}
\providecommand{\theoremname}{Theorem}

\begin{document}
\title{Nonlinear Model Reduction to \\
Fractional and Mixed-Mode Spectral Submanifolds}
\author{George Haller\thanks{Corresponding author email: georgehaller@ethz.ch},
Bálint Kaszás, Aihui Liu and Joar Axås\\
Institute for Mechanical Systems\\
ETH Zürich\\
Leonhardstrasse 21, 8092 Zürich, Switzerland\\
\\
(To appear in \emph{Chaos)}}
\maketitle
\begin{abstract}
A primary spectral submanifold (SSM) is the unique smoothest nonlinear
continuation of a nonresonant spectral subspace $E$ of a dynamical
system linearized at a fixed point. Passing from the full nonlinear
dynamics to the flow on an attracting primary SSM provides a mathematically
precise reduction of the full system dynamics to a very low-dimensional,
smooth model in polynomial form. A limitation of this model reduction
approach has been, however, that the spectral subspace yielding the
SSM must be spanned by eigenvectors of the same stability type. A
further limitation has been that in some problems, the nonlinear behavior
of interest may be far away from the smoothest nonlinear continuation
of the invariant subspace $E$. Here we remove both of these limitations
by constructing a significantly extended class of SSMs that also contains
invariant manifolds with mixed internal stability types and of lower
smoothness class arising from fractional powers in their parametrization.
We show on examples how fractional and mixed-mode SSMs extend the
power of data-driven SSM reduction to transitions in shear flows,
dynamic buckling of beams and periodically forced nonlinear oscillatory
systems. More generally, our results reveal the general function library
that should be used beyond integer-powered polynomials in fitting
nonlinear reduced-order models to data. 
\end{abstract}
\textbf{Extraction of reduced-order models for dynamical systems known
only from data has been receiving significant interest in multiple
physical disciplines. Most current approaches fit either neural networks
or linear systems to the data. The former approach tends to produce
overly complex and non-interpretable models with little predictive
power outside the range of the training data. The latter approach
is destined to fail over domains of non-linearizable behavior, such
as those with coexisting stationary states and transitions between
them. Reduction to spectral submanifolds (SSMs), i.e., special, attracting
invariant surfaces constructed over the strongest modal interactions,
is known to overcome these shortcomings but only under certain assumptions
on the modes involved. Here we remove these assumptions by uncovering
a much larger family of SSMs than previously noted. This family now
includes mixed-mode and fractional SSMs that provide data-driven reduced
models for new classes of non-linearizable physical phenomena. These
phenomena include transitions in fluid flows, buckling of beams and
global dynamics in externally forced nonlinear mechanical oscillators.}

\section{Introduction}

Reduced-order modeling seeks to identify a low-dimensional dynamical
system that captures important features of a much higher (often infinite)
dimensional physical system. This problem is too ambitious to have
a general solution and hence a number of approaches with different
assumptions have been under development. Recent reviews cover a number
of these approaches (see, e.g., \citet{amsallem15,lelievre19}) while
our focus here is specifically data-driven model reduction, which
targets physical systems defined by data sets rather than equations
(see, e.g., the focus issue by \citet{ghadami22}).

Machine learning (ML) methods are broadly used for reduced-order modeling
(\citet{brunton2016,hartman17,mohamed18,daniel20,calka21}) but tend
to produce overly complex models that lack physical interpretability
and perform poorly outside their training range (\citet{loiseau20}).
As an improvement, physics-inspired machine learning (PIML) seeks
to encode physical features (such as governing equations, symmetries
and conservation laws) into the learning process (\citet{karniadakis21}).
Remaining challenges for PIML are the handling of high-frequency features,
application to large-scale examples and the incorporation of multi-physics. 

A broadly used alternative to ML is the dynamic mode decomposition
(DMD) and its variants, reviewed recently by \citet{schmid22}. The
original DMD algorithm of \citet{schmid10} seeks to fit a linear
dynamical system to a set of observables, while the extended DMD (EDMD)
of \citet{williams15}) performs such a fit to a user-specified set
of functions defined on those observables. If the observables coincide
with the phase space variables of the system, then DMD can be justified
near fixed points as an approximation to the first-order term in the
spatial expansion of the flow map at the fixed point. In the same
setting, EDMD defined with respect to a set of polynomial functions
serves as an approximation of a truncated linearizing transformation
for a dynamical system in a neighborhood of the fixed points. For
more general observables, Koopman operator theory is often invoked
to view EDMD as an attempt to construct special functions of the observables
that span an eigenspace of the Koopman operator (\citet{rowley09,budisic12,mauroy20}).
All this rationalizes the operational use of DMD or EDMD for continuous
or discrete dynamical systems on domains where those systems exhibit
linearizable dynamics. For such problems, DMD and EDMD are easily
implementable, although the efficacy of the latter depends on the
functions of observables used in the procedure.

A number of physical phenomena, however, are fundamentally non-linearizable
over their region of interest, given that they display coexisting
isolated stationary states of various stability types, transitions
among them and bifurcations affecting them. No linearizing coordinate
change may ever cover such a range of behaviors simultaneously, because
no linear system can have more than one isolated stationary state.
Most of the interest in data-driven reduced modeling stems originally
from the need to describe, predict and control coexisting states and
the transitions connecting them, in problems ranging from transition
to turbulence (\citet{avila23}) to tipping points in climate (\citet{lenton08}).
DMD is often proposed for model reduction in such a non-linearizable
setting as well, because a high-enough dimensional linear dynamical
system can always be fitted with high accuracy to a small number of
trajectories of a nonlinear system. Outside the linear range of the
system, however, such a fit simply represents a black-box-type pattern
matching with little predictive power (see, e.g., \citet{cenedese22a}
for a demonstration of this on experimental data). 

A logical next step towards the reduced modeling of non-linearizable
phenomena is to fit a simple, low-dimensional nonlinear differential
equation to a small set of observables via regression. The resulting
algorithm, the sparse identification of nonlinear dynamics (or SINDy;
\citet{brunton16}), has been an influential tool for identifying
low-dimensional systems whose nonlinearities are expected to fall
in specific (typically polynomial) function classes. For unsteady
dynamics with a priori unknown dimensionality and nonlinearities,
however, the inherent sensitivities of this approach limit its applicability.
Recently, model reduction approaches based on manifold learning have
also been proposed (\citet{floryan22}, \citet{farzamnik23}). These
methods can be very effective in locating candidate surfaces for model
reduction operationally, without identifying the underlying dynamics
of the phase space that creates those surfaces manifolds. It remains,
therefore, unclear if these manifolds are robust under parameter changes
or under the addition of external, time-dependent forcing.

A recent alternative to these model reduction methods is reduction
to \emph{spectral submanifolds} (SSMs), which were defined by \citet{haller16}
as the unique smoothest invariant manifolds of a continuous or discrete
nonlinear system that are tangent and locally diffeomorphic to the
spectral subspaces of the linearization of the system at a fixed point.
Therefore, the internal dynamics of attracting SSMs constructed over
a slowest set of eigenspaces provides a perfect, mathematically justifiable
reduced-order model for a large set of trajectories. SSMs can cut
through the boundaries of the domain of attraction of a fixed point
and hence can carry non-linearizable dynamics. These invariant manifolds,
termed nonlinear normal modes (NNMs) at the time, were first envisioned
and formally computed in seminal work by \citet{shaw93} on nonlinear
mechanical vibrations. The existence of such manifolds was later proved
rigorously under nonresonance conditions even for infinite-dimensional
systems by \citet{cabre03}, who also showed the SSMs are unique in
a high enough smoothness class. Automated computations of SSMs and
their reduced dynamics for large finite-element problems are now available
via the open source package \emph{SSMTool} of \citet{jain2022}, whereas
an automated data-driven extraction of SSMs from arbitrary sets of
observables can be carried out via the open source \emph{SSMLearn}
and \emph{fastSSM} developed by \citet{cenedese22a} and \citet{axas22}.

Applications of data-driven SSM-reduction in mechanical system identification
and control of soft robots show strong performance of this approach
in the domain of non-linearizable behavior (see \citet{cenedese20b,kaszas22,alora23}).
Based on their mathematical construction, SSMs are also known to be
robust with respect to parameter changes and the addition of external
periodic or quasiperiodic forcing. The latter property enables them
to predict experimentally observed forced response based solely on
unforced training data (see \citet{cenedese20a,cenedese20b}) and
closed-loop response based on open-loop training data in applications
to control (\citet{alora23}). 

There have been, however, two main limitations in the underlying SSM
theory that have been hindering the application of SSM-reduction to
some other important problems. First, the spectral subspace underlying
the SSM had to be of \emph{like-mode} type, i.e., spanned by eigenvectors
of the same stability type, which does not hold for important transitions
among exact coherent states in turbulence (see, e.g., \citet{hof04}
and \citet{graham21}). Similar, mixed-mode connections are also known
to exist in the dynamic buckling of structures (\citet{abramovich18})
and dissipation induced instabilities of gyroscopic systems (\citet{krechetnikov07}).
A further limitation has been that in problems with long transitions
among different steady state, the actual transition manifold may divert
substantially from the SSM defined as the smoothest nonlinear continuation
of the invariant subspace $E$. In that case, an approximate polynomial
reduced model will fail to signal the target of the transition. Beyond
transitions in Couette flows (\citet{kaszas22}), this limitation
also creates an issue with certain constructs in renormalization group
theory (\citet{li95}), as noted by \citet{delallave97}.

Here we remove both of these limitations for generic finite-dimensional
dynamical systems. Specifically, we construct explicit parametrizations
for a significantly extended class of SSMs that also contains invariant
manifolds with mixed internal stability types and of lower smoothness
class that arises from fractional powers in their parametrization.
For both continuous and discrete dynamical systems, we derive these
representations explicitly near hyperbolic, nonresonant and semisimple
fixed points of continuous dynamical systems and discrete maps. These
three conditions may certainly be violated in small analytic examples
with symmetries, but all of them will hold generically in data sets
of physical systems. Importantly, resonances arising from repeated
eigenvalues are covered by our results. Via the construct of Poincaré
maps, the secondary SSMs obtained for maps imply the existence of
such manifolds around periodic orbits of continuous dynamical systems. 

Beyond the realm of SSM-reduction, the specific fractional-powered
representations we derive for the full SSM family should be useful
for other data-driven model reduction methods, such as EDMD. Indeed,
all methods fitting dynamical systems to data have exclusively been
using integer-powered polynomial functions. Here we show what fractional
powers will also appear generically in reduced order models up to
a given order of truncation. Using these additional fractional terms
in the fit enables one to keep the order of the polynomial model relatively
low and hence mitigate large derivatives in the model and spurious
oscillations in the SSM in larger distances from the origin. We also
discuss normal forms for the reduced dynamics in fractional SSMs,
which have not been treated in the literature.

In Section 2 we outline two simple, low-dimensional examples that
illustrate the need for extending the prior theory of primary SSMs.
In Section 3 we state our main results for continuous dynamical systems
along with several remarks on their interpretation and application.
The proofs are lengthy and involve detailed calculations which prompts
us to relegate them to appendices of a separate \emph{Supplementary
Material}. Section 4 reformulates the same results for discrete dynamical
systems. Section 5 develops a normal form theory for fractional reduced-order
models on slow two-dimensional SSMs, which is the most frequent setting
for dissipative mechanical systems exhibiting underdamped oscillatory
behavior. Section 6 discusses three examples of data-driven model
reduction to fractional and mixed-mode SSMs arising in three physical
problems: transition in Couette flows, dynamic buckling of a von Kármán
beam, and transition in a forced nonlinear oscillator system. The
first two examples involve autonomous systems, whereas the last example
involves non-autonomous dynamics.

\section{Two motivating examples\label{sec:Two-motivating-examples}}

\begin{figure}[H]
\begin{centering}
\includegraphics[width=1\textwidth]{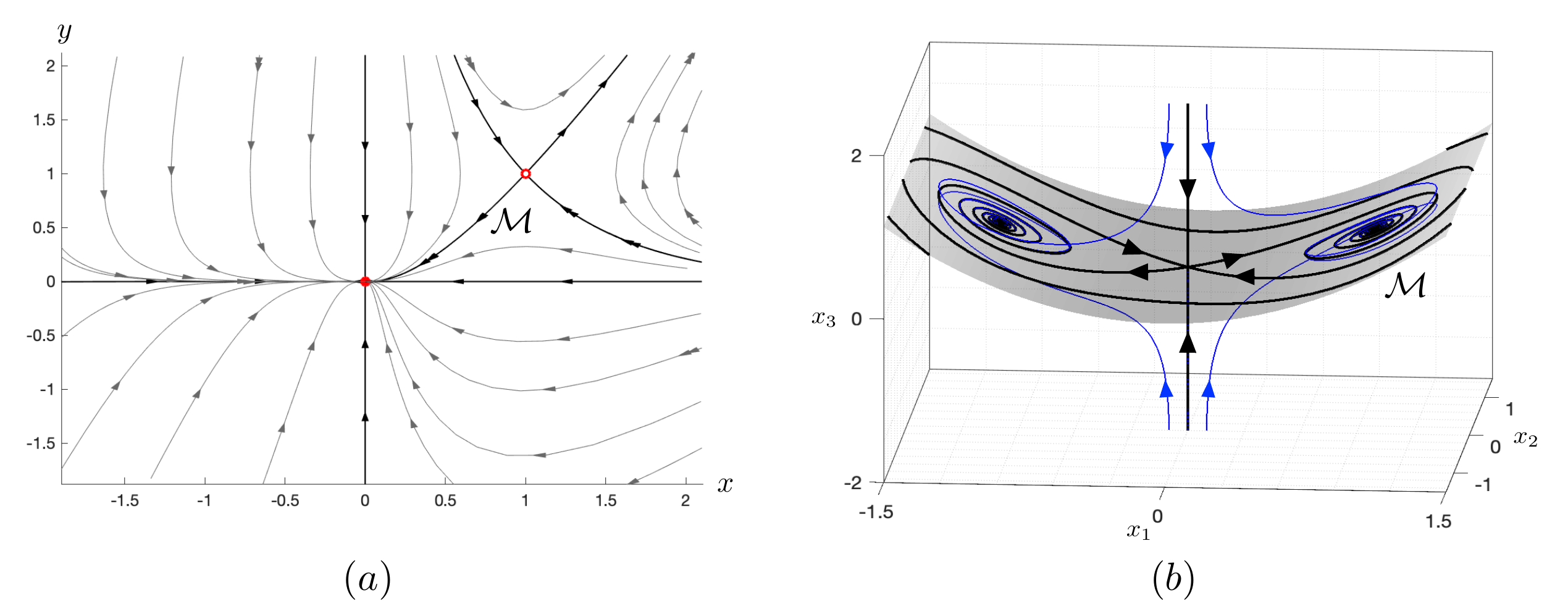}
\par\end{centering}
\caption{Simple examples illustrating the need to extend the primary SSMs used
in nonlinear model reduction to a more general family of invariant
manifolds (a) The phase portrait of system (\ref{eq:example}) for
$a=b=1$ and $c=2.5$. (b) The phase space geometry of system (\ref{eq:3D motivating example})
for $k=\sqrt{2}/2$, $a=0.5$ and $c=0.2$. \label{fig:Simple_examples}}
\end{figure}

Consider the planar system of ODEs 
\begin{align}
\dot{x} & =x(y-b),\nonumber \\
\dot{y} & =cy(x-a),\label{eq:example}
\end{align}
with the parameters $a,b,c>0$ This system has two fixed points: a
stable node $\mathbf{x}_{0}=(0,0)$ and a saddle-type fixed point
at $\mathbf{x}_{1}=\left(a,b\right)$. Also note that both the $x$
and the $y$ axes are invariant lines. A representative phase portrait
of (\ref{eq:example}) for $a=b=1$ and $c=2.5$ is shown in Fig.
\ref{fig:Simple_examples}a. This phase portrait mimics the geometry
of a number of important transition problems in fluids, such as transitions
between stationary states in Couette flows (see, e.g., \citet{page19})
and transitions to turbulence in pipe flow (see, e.g., \citet{skufca06}).
In the latter context, the saddle point mimics an exact coherent state
(ECS) and its stable manifold models the\emph{ edge of chaos}, the
border between laminar and turbulent behavior. An important objective
of reduced-order modeling in these flows is to establish whether specific
initial conditions decay to the origin or move towards the interior
of the turbulent region.

A Galerkin projection onto POD modes does not produce a feasible reduced-order
model for this simple problem, as we show in Appendix A of the \emph{Supplementary
Material}. Linear data-driven techniques, such as DMD and its variants
(\citet{schmid22}), are clearly unable to capture simultaneously
the two fixed points and the orbit connecting them, as was explicitly
illustrated in the context of Couette flows by \citet{page19}. 

The slow SSM of the origin, i.e., the unique smoothest nonlinear continuation
of the slower eigenspace, is the $x$ axis, on which the reduced dynamics
is $\dot{x}=-bx$. This gives a correct reduced-order model for the
asymptotic behavior of all trajectories in the domain of attraction
of the origin, yet fails to provide any information about the truly
nonlinear dynamics. This is undesirable as both the coexisting saddle
point and the separatrix $\mathcal{M}$ connecting it to the origin
can be brought arbitrarily close to the origin by selecting small
enough values for the parameters $a$ and $b$. 

A simple calculation reveals that any trajectory of the linear part
of \eqref{eq:example} can be written as a graph $y=K\left|x\right|^{\frac{ac}{b}}$
for some choice of the parameter $K\in\mathbb{R}$. This suggest that
the manifold $\mathcal{M}$ containing the two fixed points may also
just be $\left[\frac{ac}{b}\right]$ times differentiable at the origin,
with the brackets referring to the integer part of a real number.
Our first main objective in this paper is to establish the general
existence and smoothness of secondary SSMs like $\mathcal{M}$ in
nonlinear dynamical systems and obtain data-driven reduced models
on such manifolds. These SSMs, however, turn out to lie outside the
scope of usual equation- and data-driven approximation methods that
use polynomials with integer powers. As we will see for this example
in Appendix A of the Supplementary material, use of fractional SSMs
enables the data-driven construction for an accurate 1D reduced-order
model for trajectories near the transition orbit. Additionally, we
will show for the planar Couette flow in Section \ref{sec:Couette flow example}
that use of fractional SSMs increases the accuracy of the SSM-reduced
model without the inclusion of higher-order terms in that model.

As a second example, consider the system 
\begin{align}
\dot{x}_{1} & =x_{2},\nonumber \\
\dot{x}_{2} & =x_{1}-cx_{2}-x_{1}^{3},\nonumber \\
\dot{x}_{3} & =-kx_{3}+kax_{1}^{2}+2ax_{1}x_{2}.\label{eq:3D motivating example}
\end{align}
for which $\mathcal{M}=\left\{ x\in\mathbb{R}^{3}\colon x_{3}=ax_{1}^{2}\right\} $
is a unique, two-dimensional (2D), attracting $C^{\infty}$ invariant
manifold that is tangent to the $x_{3}=0$ spectral subspace. Within
this manifold, the restricted dynamics is just that of a damped Duffing
oscillator with two stable spirals and one saddle point at the origin,
as shown in Fig. \ref{fig:Simple_examples}b. 

The $x_{3}=0$ spectral subspace contains both a stable and an unstable
eigenvector of the origin, and hence the results of \citet{cabre03}
are not applicable to conclude the existence of $\mathcal{M}$ purely
based on the spectrum of the linearized ODE at the origin. The more
classic result of Irwin on pseudo-unstable manifolds (\citet{delallave95})
do not apply for the parameter values used in Fig. \ref{fig:Simple_examples}
(see Appendix B of the \emph{Supplementary Material}). Our second
main objective in this paper is to derive reduced-order models on
secondary SSMs such as $\mathcal{M}$ in Fig. \ref{fig:Simple_examples}b.
As we will see in Section \ref{sec:Buckling beam example}, an important
physical application of such mixed-mode SSMs is the reduced-order
modeling of dynamic buckling of beams.

\section{Generalized SSMs near fixed points of autonomous dynamical systems\label{sec:Generalized-SSMs-for-continuous-DS}}

We consider a smooth system of $n$-dimensional ODEs of the form
\begin{equation}
\dot{x}=\mathcal{A}x+f(x),\qquad x\in\mathbb{R}^{n},\quad\mathcal{A}\in\mathbb{R}^{n\times n},\quad f=\mathcal{O}\left(\left|x\right|^{2}\right)\in C^{\infty},\qquad1\leq n<\infty,\label{eq:original system}
\end{equation}
and assume that its fixed point at $x=0$ is hyperbolic, i.e., the
spectrum $\mathrm{spect}\left(\mathcal{A}\right)=\left\{ \lambda_{1},\ldots,\lambda_{n}\right\} $
of $\mathcal{A}$ satisfies 
\begin{equation}
0\notin\mathrm{Re}\left[\mathrm{spect}\left(\mathcal{A}\right)\right],\label{eq:hyperbolicity}
\end{equation}
with $\mathrm{Re}\left[\mathrm{spect}\left(\mathcal{A}\right)\right]$
denoting the real part of the spectrum of the linear operator $\mathcal{A}$.
We note that the $C^{\infty}$ assumption on $f$ can be relaxed to
finite smoothness if $\mathrm{Re}\left[\mathrm{spect}\left(\mathcal{A}\right)\right]\subset(-\infty,0)$
or $\mathrm{Re}\left[\mathrm{spect}\left(\mathcal{A}\right)\right]\subset(0,\infty)$
holds (see statement (v) of Theorem \ref{thm:main1}). We also note
that our results in Section \ref{sec:Generalized-SSMs-for-discrete-DS}
will also allow for the addition of small time-periodic forcing term
to (\ref{eq:original system}) (see Remark \ref{rem: extension of Thm. 1 to time-periodic perturbations}).

We further assume that $\mathcal{A}$ is semisimple, i.e., the geometric
and algebraic multiplicities of any potential repeated eigenvalue
of $\mathcal{A}$ are equal. This is generally true for physical oscillatory
systems in which possible repeated eigenvalues arise from symmetry
and hence each eigenvalue has its own associated oscillatory mode
generated by an independent eigenvector. We finally assume that with
the exception of possible $1\colon1$ resonances created by repeated
eigenvalues of $\mathcal{A}$, there are no further resonances within
$\mathrm{spect}\left(\mathcal{A}\right)$, i.e., we have
\begin{equation}
\lambda_{j}\neq\sum_{k=1}^{n}m_{k}\lambda_{k},\quad m_{k}\in\mathbb{N},\quad\sum_{k=1}^{n}m_{k}\geq2,\quad j=1,\ldots,n.\label{eq:nonresonance conditions}
\end{equation}

While these nonresonance conditions will typically hold for generic
parameter configurations of typical dissipative systems, they can
easily fail on simple toy examples with non-generic parameter choices.
For instance, a simple 2D saddle-type fixed point with eigenvalues
$\lambda_{1}=-1$ and $\lambda_{2}=1$ will violate infinitely many
of the conditions (\ref{eq:nonresonance conditions}) for $j=1$,
$m_{1}=m_{2}+1$ and $j=2$, $m_{2}=m_{1}+1$. The main motivation
for the present paper is the data-driven modeling of dissipative physical
systems in which such exceptional parameter configurations are unlikely.
We stress, however, that repeated eigenvalues arising from a physical
symmetry are \emph{not} excluded by conditions (\ref{eq:nonresonance conditions}).

We now select and fix an arbitrary spectral subspace $E\subset\mathbb{R}^{n}$
of $\mathcal{A}$ and seek to find its nonlinear continuation in the
full system (\ref{eq:original system}). In practical terms, we are
interested in nonlinear interactions along which all modes can be
enslaved to a set of linear master modes spanning $E$. We assume
that $\mathrm{spect}\left(\mathcal{A}\vert_{E}\right)$ contains $p$
purely real eigenvalues and $q$ pairs of complex conjugate eigenvalues.
Similarly, we assume that $\mathrm{spect}\left(\mathcal{A}\right)-\mathrm{spect}\left(\mathcal{A}\vert_{E}\right)$
contains $r$ purely real eigenvalues and $s$ pairs of complex conjugate
eigenvalues. In that case, after a linear change of coordinates, we
can assume a partition of the $x$ coordinate in the form
\begin{equation}
x=\left(u,\left(a_{1},b_{1}\right),\ldots,\left(a_{q},b_{q}\right),v,\left(c_{1},d_{1}\right),\ldots,\left(c_{s},d_{s}\right)\right)\in\mathbb{R}^{p}\times\underbrace{\mathbb{R}^{2}\times\ldots\times\mathbb{R}^{2}}_{q}\times\mathbb{R}^{r}\times\underbrace{\mathbb{R}^{2}\times\ldots\times\mathbb{R}^{2}}_{s},\label{eq:linear change of coordinates}
\end{equation}
with $\quad p+2q+r+2s=n,$ in which the coefficient matrix $\mathcal{A}$
of the linear part of system (\ref{eq:original system}) takes its
real Jordan canonical form 
\begin{equation}
\mathcal{A}=\mathrm{diag}\left[A,B_{1},\ldots,B_{q},C,D_{1},\ldots,D_{s}\right],\label{eq:linsystem1-1}
\end{equation}
with
\begin{align}
A & =\mathrm{diag}\left[\lambda_{1},\ldots,\lambda_{p}\right]\in\mathbb{R}^{p\times p},\quad B_{k}=\left(\begin{array}{cc}
\alpha_{k} & -\omega_{k}\\
\omega_{k} & \alpha_{k}
\end{array}\right)\in\mathbb{R}^{2\times2},\quad k=1,\ldots,q,\nonumber \\
C & =\mathrm{diag}\left[\kappa_{1},\ldots,\kappa_{r}\right]\in\mathbb{R}^{r\times r},\quad D_{m}=\left(\begin{array}{cc}
\beta_{m} & -\nu_{m}\\
\nu_{m} & \beta_{m}
\end{array}\right)\in\mathbb{R}^{2\times2},\quad m=1,\ldots,s.\label{eq:linsystem1-2}
\end{align}
Some of the eigenvalues $\lambda_{j},$ $\alpha_{k}+i\omega_{k}$,
$\kappa_{\ell}$ and $\beta_{m}+i\nu_{m}$ in the spectrum of $\mathcal{A}$
may be repeated, as $1:1$ resonances are allowed by the nonresonance
conditions (\ref{eq:nonresonance conditions}).

Our main result is the simplest to state in complexified coordinates.
To this end, we introduce the variables
\begin{equation}
z=a+ib\in\mathbb{C}^{q},\quad w=c+id\in\mathbb{C}^{s},\label{eq:complexification}
\end{equation}
and the piecewise constant functions
\begin{equation}
K_{\ell j}(\zeta)=\left\{ \begin{array}{c}
K_{\ell j}^{+}\qquad\zeta>0,\\
\\
K_{\ell j}^{-}\qquad\zeta\leq0,
\end{array}\right.\qquad O_{mj}(\zeta)=\left\{ \begin{array}{c}
O_{mj}^{+}\qquad\zeta>0,\\
\\
O_{mj}^{-}\qquad\zeta\leq0,
\end{array}\right.\label{eq:piecewise constant functions}
\end{equation}
with constant families $K_{\ell j}^{\pm}\in\mathbb{R}$ and $O_{mj}^{\pm}\in\mathbb{C}$.
We also define the continuous function families
\begin{align}
\mathcal{V}_{\ell}(u,z) & =\sum_{j=1}^{p}K_{\ell j}(u_{j})\left|u_{j}\right|^{\frac{\kappa_{\ell}}{\lambda_{j}}}+\sum_{k=1}^{q}L_{\ell k}\left|z_{k}\right|^{\frac{\kappa_{\ell}}{\alpha_{k}}},\quad\ell=1,\ldots,r,\nonumber \\
\mathcal{\mathcal{E}}_{m}(u,z) & =\left[\sum_{j=1}^{p}O_{mj}(u_{j})\left|u_{j}\right|^{\frac{\beta_{m}}{\lambda_{j}}}+\sum_{k=1}^{q}Q_{mk}\left|z_{k}\right|^{\frac{\beta_{m}}{\alpha_{k}}}\right]e^{i\left(\sum_{j=1}^{p}\frac{\nu_{m}}{\lambda_{j}}\log\left|u_{j}\right|+\sum_{k=1}^{q}\frac{\nu_{m}}{\alpha_{k}}\log\left|z_{k}\right|\right)},\quad m=1,\ldots,s,\nonumber \\
\label{eq:complex SSM 1}
\end{align}
where $Q_{mk}\in\mathbb{C}$ and 
\begin{align}
K_{\ell j} & (u_{j})\equiv0,\quad\mathrm{if}\quad\frac{\kappa_{\ell}}{\lambda_{j}}<0,\qquad O_{mj}(u_{j})\equiv0,\quad\mathrm{if}\quad\frac{\beta_{m}}{\lambda_{j}}<0,\nonumber \\
L_{\ell k} & =0,\quad\mathrm{if}\quad\frac{\kappa_{\ell}}{\alpha_{k}}<0,\qquad Q_{mk}=0,\quad\mathrm{if}\quad\frac{\beta_{m}}{\alpha_{k}}<0.\label{eq:general_fractional_SSM2}
\end{align}
As seen in the proof of our upcoming main theorem, the function families
(\ref{eq:complex SSM 1}) describe all continuous invariant graphs
over $E$ in the linear part of system (\ref{eq:original system}). 

For any positive integer $N$, any vector $e\in\mathbb{C}^{N}$ and
nonnegative integer vector $\mathbf{k}\in\mathbb{N}^{N}$, we will
use the notation
\[
e^{\mathbf{k}}:=e_{1}^{k_{1}}e_{2}^{k_{2}}\cdot\ldots\cdot e_{N}^{k_{N}}\in\mathbb{C}.
\]
 Finally, we will use the notation 
\begin{equation}
\sum_{2\leq\mathrm{monomial}\mathrm{\,order}\leq\mathcal{K}}M_{\mathbf{k}}\left(u,z\right)\label{eq:little o notation}
\end{equation}
for expansions involving monomials $M_{\mathbf{k}}\left(u,z\right)$
of general positive (but potentially fractional) powers of $u$ and
$z$. These monomials will be indexed by a multi-index $\mathbf{k}\in\mathbb{N}^{p+2q+r+2s}$
but their order (referred to as \emph{monomial order}) will be a homogeneous
linear function of $\mathbf{k}$, rather than $\mathbf{k}$ itself.
We will explicitly write out the monomial order as functions in simpler
settings (see, e.g., Propositions \ref{prop:1D real  spectral subspace}
and \ref{prop:2D complex spectral subspace}) in which the notation
does not become too involved.

With these definitions, our main result can be stated as follows.
\begin{thm}
\label{thm:main1}
\end{thm}
\begin{description}
\item [{(i)}] For any integer order $\mathcal{K}\geq2$ of approximation,
there exists a unique set of coefficients $h_{\mathbf{k}}\in\mathbb{C}^{r+s}$
such that the full family $\mathcal{W}\left(E\right)$ of SSMs tangent
to the spectral subspace $E$ of the nonlinear system (\ref{eq:original system})
can locally be written near the origin as 
\begin{align}
\left(\begin{array}{c}
v\\
w
\end{array}\right) & =\mathcal{G}(u,z)=\left(\begin{array}{c}
\mathcal{V}(u,z)\\
\mathcal{\mathcal{E}}(u,z)
\end{array}\right)+\sum_{2\leq\mathrm{monomial}\mathrm{\,order}\leq\mathcal{K}}h_{\mathbf{k}}u^{\mathbf{k}_{1}}z^{\mathbf{k}_{2}}\bar{z}^{\mathbf{k}_{3}}\mathcal{V}^{\mathbf{k}_{4}}\mathcal{\mathcal{E}}^{\mathbf{k}_{5}}\bar{\mathcal{\mathcal{E}}}^{\mathbf{k}_{6}}+o\left(\left|\left(u,z\right)\right|^{\mathcal{K}}\right),\label{eq:SSM in original coordinates real-1-5}\\
 & \,\,\,\,\,\,\mathbf{k}_{1}\in\mathbb{N}^{p},\quad\mathbf{k}_{2},\mathbf{k}_{3}\in\mathbb{N}^{q},\quad\mathbf{k}_{4}\in\mathbb{N}^{r},\quad\mathbf{k}_{5},\mathbf{k}_{6}\in\mathbb{N}^{s},\quad\mathbf{k}=\left(\mathbf{k}_{1},\mathbf{k}_{2},\mathbf{k}_{3},\mathbf{k}_{4},\mathbf{k}_{5},\mathbf{k}_{6}\right),\nonumber 
\end{align}
where $\mathcal{V}(u,z)\in\mathbb{R}^{r}$ and $\mathcal{\mathcal{E}}(u,z)\in\mathbb{C}^{s}$
are defined via formulas (\ref{eq:complex SSM 1}), with $K_{\ell j}$,
$O_{mj}$, $L_{\ell k}$ and $Q_{mk}$ satisfying the constraints
(\ref{eq:general_fractional_SSM2}) but otherwise chosen arbitrarily. 
\item [{(ii)}] The reduced dynamics on the individual members of the family
$\mathcal{W}\left(E\right)$ are given by
\begin{equation}
\left(\begin{array}{c}
\dot{u}\\
\dot{z}
\end{array}\right)=\left(\begin{array}{cccc}
A & 0 & 0 & 0\\
0 & B_{1} & 0 & 0\\
0 & 0 & \ddots & 0\\
0 & 0 & 0 & B_{q}
\end{array}\right)\left(\begin{array}{c}
u\\
z
\end{array}\right)+f^{\left(u,z\right)}\left(u,z,\mathcal{G}(u,z)\right),\label{eq:reduced dynamics continuous}
\end{equation}
where $f^{\left(u,z\right)}$ denotes the $(u,z)$ coordinate component
of $f=\left(f^{\left(u,z\right)},f^{v},f^{w}\right)$. 
\item [{\emph{(iii)}}] \emph{If and 
\begin{equation}
\frac{\kappa_{\ell}}{\lambda_{j}}<0,\qquad\frac{\beta_{m}}{\lambda_{j}}<0,\qquad\frac{\kappa_{\ell}}{\alpha_{k}}<0,\qquad\frac{\beta_{m}}{\alpha_{k}}<0\label{eq:condition for SSM uniqueness in any smoothness class}
\end{equation}
hold for all values of the indices $j,k,\ell$ and $m$, then $\mathcal{W}\left(E\right)$
is unique and of class $C^{\infty}$. Otherwise, $\mathcal{W}\left(E\right)$
is a family of non-unique SSMs whose typical member is of smoothness
class $C^{\eta}$ with
\begin{equation}
\eta=\mathrm{Int}\left[\min_{j,k,\ell,m}^{+}\left\{ \frac{\kappa_{\ell}}{\lambda_{j}},\quad\frac{\beta_{m}}{\lambda_{j}},\quad\frac{\kappa_{\ell}}{\alpha_{k}},\quad\frac{\beta_{m}}{\alpha_{k}}\right\} \right],\label{eq:typical smoothness measure}
\end{equation}
where the $\mathrm{Int}$ refers to the integer part of a positive
number and $\min^{+}$ refers to a minimum of the positive elements
of a list of numbers that follow. Even in that case, however, there
exists a unique member $\mathcal{W}^{\infty}(E)$ of the SSM family
that is of class $C^{\infty}$ and satisfies $K_{\ell j}^{\pm}\equiv O_{mj}^{\pm}=L_{\ell k}=Q_{mk}=0$
for all $j,k,\ell$ and $m$ in eq. (}\ref{eq:complex SSM 1}). 
\item [{\emph{(iv)}}] If system (\ref{eq:SSM in original coordinates real-1-5})
is $C^{R}$ in a parameter vector $\mu\in\mathbb{R}^{P}$, then $\mathcal{W}\left(E\right)$
is also $C^{R}$ in $\mu$.
\item [{(v)}] If the function $f$ only has finite smoothness, i.e., $f=\mathcal{O}\left(\left|x\right|^{2}\right)\in C^{\rho}$
and $\mathrm{Re\,\lambda_{j}}$ has the same nonzero sign for all
$j=1,\ldots,N$, then statements (i)-(iv) still hold for all $\mathcal{K}\leq\rho$. 
\item [{\emph{(vi)}}] If the function $f$ is analytic in a neighborhood
of the origin and $\mathrm{Re\,\lambda_{j}}$ has the same nonzero
sign for all $j=1,\ldots,N$, then the series expansion (\ref{eq:SSM in original coordinates real-1-5})
converges near the origin as $\mathcal{K}\to\infty$. 
\end{description}
\begin{proof}
See Appendix C of the \emph{Supplementary Material}.
\end{proof}
\begin{rem}
\label{rem: explanation on obtaining primary SSMs} {[}\textbf{Relation
to primary SSMs}{]} We note that to obtain classical stable and unstable
manifolds as well as the primary (smoothest) SSMs defined by \citet{haller16},
all constants $K_{\ell j}$, $O_{mj}$, $L_{\ell k}$ and $Q_{mk}$
in the definitions (\ref{eq:complex SSM 1}) have to be set to zero
and the signs of all $\lambda_{j}$ and $\alpha_{k}$ involved must
be the same. In contrast, the full family of invariant manifolds covered
by Theorem \ref{thm:main1} is significantly larger and includes secondary
(fractional and mixed-mode) SSMs. The existing open-source \emph{SSMLearn}
package of \citet{cenedese21} can be used to compute primary mixed-mode
SSMs as well, even though the SSM theory originally supporting that
package was only available for fixed points that are asymptotically
stable in forward or backward time. This is because the integer-powered
polynomial expansions for mixed-mode SSMs satisfy the same type of
invariance equations as their like-mode counterparts and hence can
be approximated from data in the same fashion. 
\end{rem}
\begin{rem}
\label{rem: explanation on relevance for nonlinearizable dynamics}{[}\textbf{Relation
to non-linearizable dynamics}{]} Increasing $\mathcal{K}$ in the
general formula (\ref{eq:SSM in original coordinates real-1-5}) generally
decreases the size of the open neighborhood of the origin in which
this representation of the SSM family is valid. Specifically, in the
limit of $\mathcal{K}\to\infty$, only the part of the SSM family
$\mathcal{W}\left(E\right)$ is captured by formula (\ref{eq:SSM in original coordinates real-1-5})
that carries linearizable dynamics, given that the main tool used
in obtaining this asymptotic limit of this formula is the $C^{\infty}$
linearization result of \citet{sternberg58}. In contrast, using a
finite $\mathcal{K}$ in a data-driven identification of $\mathcal{W}\left(E\right)$
enables the accurate approximation of SSMs on larger domains containing
non-linearizable dynamics, such as coexisting isolated compact invariant
sets and transitions among them. We will illustrate this on specific
examples in Sections \ref{sec:Couette flow example}-\ref{sec:Shaw-Pierre example}.
\end{rem}
\begin{rem}
\label{rem: handling or resonances}{[}\textbf{Resonances}{]} As noted
after formula (\ref{eq:nonresonance conditions}), $1\colon1$ resonances
are allowed among the eigenvalues of $\mathcal{A}$ by the assumptions
of Theorem 1. For other resonances, the nonresonance condition (\ref{eq:nonresonance conditions})
is substantially less restrictive than what one might be used to in
the classical vibrations literature. For instance, $2\lambda_{1}=3\lambda_{2}$
satisfies condition (\ref{eq:nonresonance conditions}) and hence
is \emph{not }a resonance in our discussion unless further matches
arise with the remaining eigenvalues, which has low probability in
a physical system. Similarly, the case $\mathrm{Im}\lambda_{1}=3$$\mathrm{Im}\lambda_{2}$
is considered a $1:3$ frequency-resonance in the nonlinear vibrations
literature, but it is only a resonance in the sense of formula (\ref{eq:nonresonance conditions})
if we also have the same resonance between the decay rates of the
two corresponding modes, i.e., $\mathrm{Re}\lambda_{1}=3\mathrm{Re}\lambda_{2}$
holds simultaneously. This is again highly unlikely for physical system,
especially in the data driven setting. If eigenvalue relations violating
(\ref{eq:nonresonance conditions}) do need to be considered for some
reason, then primary, like-mode SSMs defined over spectral subspaces
incorporating the resonant eigenvectors still exist by the theory
of \citet{cabre03} (see \citet{li22a,li22b} for specific examples).
Finally, under resonances of order $\sum_{k=1}^{n}m_{k}>\mathcal{K}$,
formula (\ref{eq:SSM in original coordinates real-1-5}) still provides
an an expression for an approximately invariant manifold with invariance
error of $o\left(\left|\left(u,z\right)\right|^{\mathcal{K}}\right)$.
This follows from the fact that a finite polynomial (and hence analytic)
transformation linearizing the system up to order $\mathcal{K}$ still
exists. In a data-driven setting, such a small error in locating invariant
manifolds is generally unavoidable for higher values of $\mathcal{K}$
even in the absence of resonances. 
\end{rem}
\begin{rem}
\label{rem: explanation for constraints on parameters}{[}\textbf{Constraints
on parameters}{]} The constraints on the parameters in (\ref{eq:general_fractional_SSM2})
are relevant when the stability types of the eigenvectors aligned
with $v_{\ell}$ is the opposite of those aligned with $u_{j}$ or
$z_{k}$, in which case any non-flat invariant graph over $E$ blows
up at the origin. If contrast, if $\frac{\kappa_{\ell}}{\lambda_{j}}>0$
and $\frac{\kappa_{\ell}}{\alpha_{k}}>0$ hold, then any nonzero $K_{\ell j}^{\pm}$
and $L_{\ell k}$ are allowed in the expression (\ref{eq:complex SSM 1}). 
\end{rem}
\begin{rem}
\label{rem:explanation on typical smoothness}{[}\textbf{Generic smoothness
in the $\mathcal{W}\left(E\right)$ family}{]} Statement of (iii)
of Theorem \ref{thm:main1} specifies the order $\eta$ of differentiability
of a generic fractional SSM in the family of $\mathcal{W}\left(E\right)$.
For such generic members of the family, all coefficients appearing
in the definition (\ref{eq:complex SSM 1}) of the functions $\mathcal{V}$
and $\mathcal{E}$ are nonzero. Nongeneric subfamilies of the $\mathcal{W}\left(E\right)$
family will have higher smoothness classes equal to the integer part
of one of the positive members in the list of quotients in formula
(\ref{eq:typical smoothness measure}). However, as we illustrate
in Fig. \ref{fig: smoothness_illustration} on a specific case, these
nongeneric subfamilies form measure zero, nowhere dense subsets of
the full $\mathcal{W}\left(E\right)$ family. 
\begin{figure}[h]
\centering{}\includegraphics[width=0.75\textwidth]{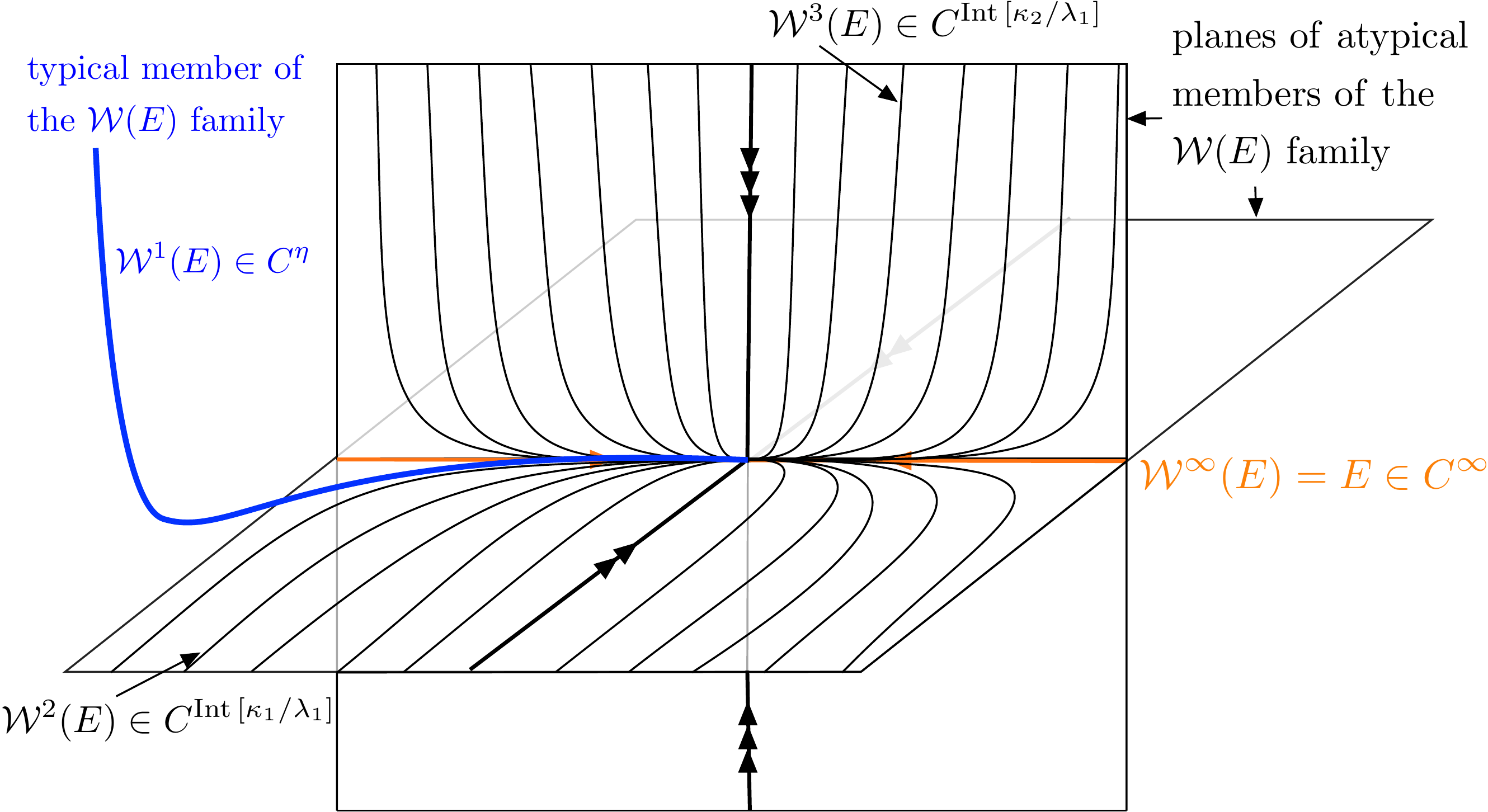}\caption{Various members of the SSM family $\mathcal{W}(E)$ in the case of
$j=1$ , $k=0$, $\ell=1$ and $m=0$, shown after a $C^{\infty}$
linearizing transformation near the origin has been performed (see
the proof of Theorem \ref{thm:main1}). $\mathcal{W}^{\infty}(E)=E$:
unique member of smoothness class $C^{\infty}$ ; $\mathcal{W}^{1}(E)$:
generic member of smoothness class $C^{\eta}$; $\mathcal{W}^{2}(E)$:
nongeneric member of smoothness class $C^{\mathrm{Int}\,\left[\kappa_{1}/\lambda_{1}\right]}$
; $\mathcal{W}^{3}(E)$: nongeneric member of class $C^{\mathrm{Int}\,\left[\kappa_{2}/\lambda_{1}\right]}$
.\label{fig: smoothness_illustration}}
\end{figure}
\end{rem}
Finally, by their construction in the proof of Theorem \ref{thm:main1},
generic members of the fractional SSM family $\mathcal{W}\left(E\right)$
guaranteed by statement (i) of the theorem are actually always of
class $C^{\infty}$ \emph{away} from the origin. At the origin, their
Hölder smoothness class is $C^{\eta,\alpha}$ for some $0\leq\alpha<1$. 

\begin{rem}
\label{rem: explanation for piecewise linear functions}{[}\textbf{Use
of piecewise constant coefficients}{]} The use of the piecewise constant
functions (\ref{eq:piecewise constant functions}) enables us to pair
up arbitrary positive and negative branches of fractional SSMs to
obtain a single, continuous invariant manifold. This is illustrated
in Fig. \ref{fig:extended manifold} for the linear part of system
(\ref{eq:original system}) which already exhibits the same geometry.
Replacing the function $K_{\ell j}(u_{j})$ simply with a single constant
$K_{\ell j}$ is possible but limits the general family of graphs
to symmetric ones, as seen in Fig. \ref{fig:extended manifold}. The
same observation is valid for $O_{mj}(u_{j})$.
\begin{figure}[H]
\centering{}\includegraphics[width=0.6\textwidth]{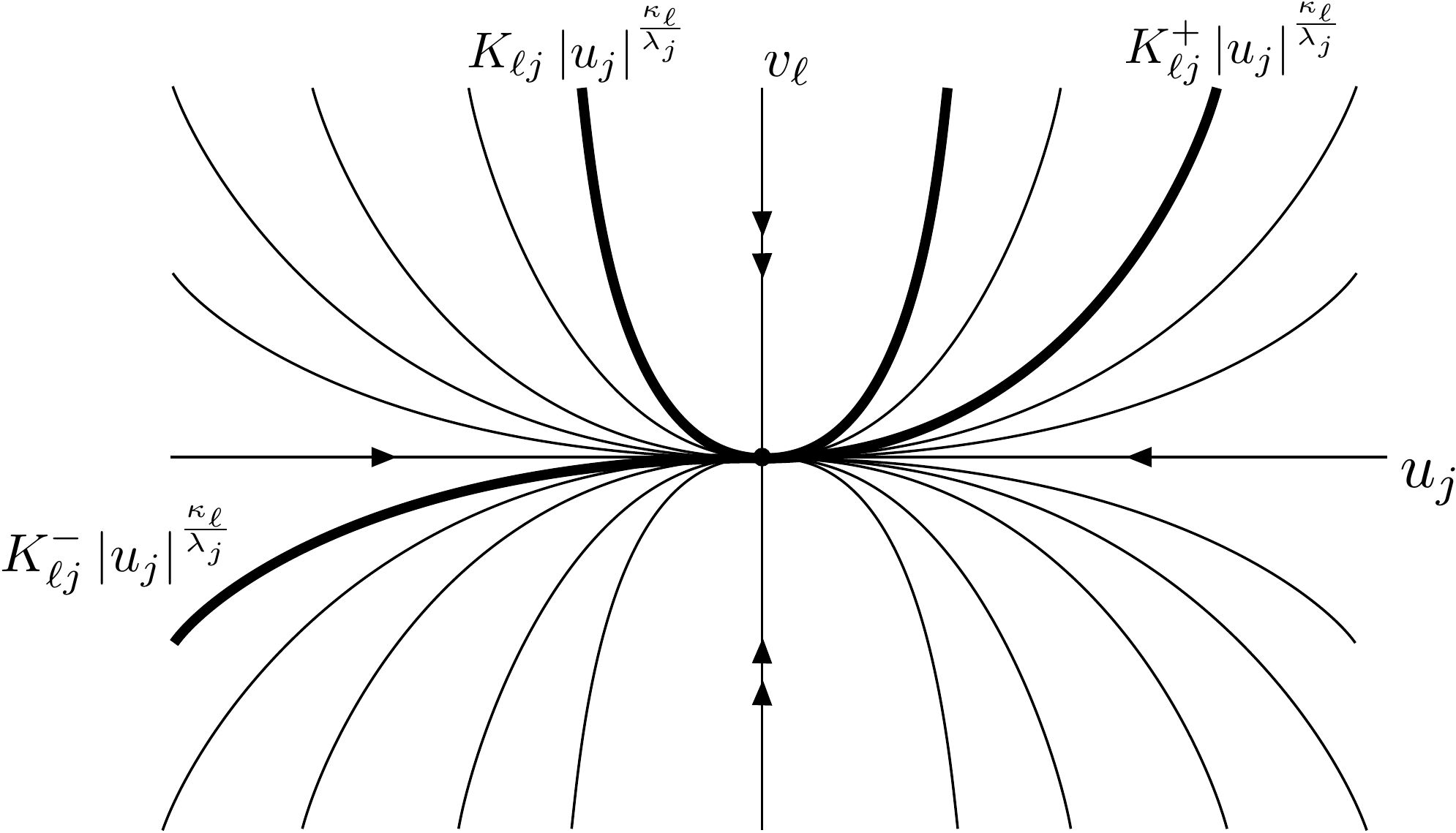}\caption{Infinitely many invariant graphs (fractional SSMs) of smoothness class
$C^{\mathrm{int}\left[\frac{\kappa_{\ell}}{\lambda_{j}}\right]}$
pass through any point of the $(u_{j},v_{\ell})$ plane in the linear
part of system \eqref{eq:original system}. The symmetric members
of this family are of the form $v_{\ell}=K_{\ell j}\left|u_{j}\right|^{\frac{\kappa_{\ell}}{\lambda_{j}}}$
while the asymmetric members are of the general form $v_{\ell}=K_{\ell j}(u_{j})\left|u_{j}\right|^{\frac{\kappa_{\ell}}{\lambda_{j}}}$.
The single $C^{\infty}$ invariant graph (primary SSM) over the $u_{j}$
axis is the $u_{j}$ axis itself. \label{fig:extended manifold}}
\end{figure}
\end{rem}
\begin{rem}
\label{rem:explanation on overfit from nearly equal powers}{[}\textbf{Data-driven
construction of fractional and mixed-mode SSMs}{]} Primary mixed-mode
SSMs require no detailed information on the spectrum of $\mathcal{A}$
outside their underlying spectral eigenspaces. For this reason, mixed-mode
SSMs can be directly constructed from data featuring decaying nonlinear
oscillations near $E$ via the open source Matlab codes \emph{SSMLearn}
and \emph{fastSSM} developed by \citet{cenedese22a} and \citet{axas22}
(see https://github.com/haller-group). In contrast, the construction
of fractional SSMs requires decay-rate information both inside and
outside $E$, the latter of which is not immediately recognizable
from decaying oscillation data near $E$. In the Conclusions section,
we mention a few data-driven strategies to obtain this spectral information
external to $E$. Once this information is available, the construction
of mixed-mode SSMs follows the same algorithm used in \emph{SSMLearn}
and \emph{fastSSM} but with the fractional-powered terms added to
the integer-powered representation of the invariant manifolds. If,
however, some of the fractional powers appearing in the formulas (\ref{eq:complex SSM 1})
are close to integers, then those fractional terms are advisable to
omit in a data-driven construction of the SSM to avoid an overfit.
All these steps in a data-driven implementation of Theorem \ref{thm:main1}
are currently under development and will appear elsewhere. In all
our upcoming examples in Sections \ref{sec:Couette flow example}-\ref{sec:Shaw-Pierre example},
we obtain linear spectral information outside $E$ from the governing
equations but construct the fractional and mixed-mode manifolds in
a data-driven fashion, using only numerically generated solutions
of those governing equations.
\end{rem}
\begin{rem}
\label{rem:uniqueness of manifolds}{[}\textbf{Uniqueness of invariant
manifolds}{]} The $C^{\infty}$ linearizing conjugacies of \citet{sternberg58}
used in proving Theorem 1 have only been known to be unique under
the assumptions of statement (v) of the theorem (\citet{kvalheim21}).
This, however, does not affect the uniqueness of \emph{$\mathcal{W}^{\infty}(E)$
}or the uniqueness of the members of the full SSM family $\mathcal{W}(E)$.
Indeed, as members of $\mathcal{W}(E)$ provide a well-defined foliation
by invariant surfaces of a full neighborhood of the origin under a
given $C^{\infty}$ linearizing conjugacy, another $C^{\infty}$ linearization
conjugacy cannot change the elements and their smoothness in this
foliation without violating their invariance. Only the specific parametrization
of the fractional and primary SSMs in $\mathcal{W}(E)$ can change
when one picks another linearizing conjugacy. This is also the reason
behind the uniqueness of $\mathcal{W}\left(E\right)$ in statement
(iii), given that conditions (\ref{eq:condition for SSM uniqueness in any smoothness class})
cause the linearized system to have a unique continuous invariant
manifold ($E$ itself) that is tangent to $E$ at the origin and has
the same dimension as $E$. This manifold then has the same unique
$C^{\infty}$ preimage under any $C^{\infty}$ linearizing transformation,
and hence the corresponding $\mathcal{W}\left(E\right)$ is a unique
manifold in the nonlinear system (\ref{rem: explanation on relevance for nonlinearizable dynamics}).
\end{rem}
\begin{rem}
\label{rem:choice of dimension for SSM}{[}\textbf{Choice of the dimension
of the SSM}{]} Theorem 1 provides a hierarchical chain of SSM families,
$W\left(E_{1}\right)\subset W\left(E_{2}\right)\subset\ldots\subset W\left(E_{n-1}\right)$,
with each family acting as a nonlinear, invariant continuation of
a corresponding spectral subspace in the hierarchical chain $E_{1}\subset E_{2}\subset\ldots\subset E_{n-1}$
of spectral subspaces. In applications, the asymptotic dynamics are
generally captured by the slowest family $W\left(E_{1}\right)$ of
these SSM families. If predictions are required for shorter time scales
as well, one can gradually increase the size of the family by considering
$W\left(E_{k}\right)$ under increasing $k.$ This procedure includes
more and more of the transient dynamics while keeping the core asymptotic
dynamics unchanged. One reason for considering $W\left(E_{k}\right)$
with $k>1$ in a data-driven setting is if the data shows significant
interaction of the first mode with higher modes up to order $k$.
In most practical applications of SSMs with larger spectral gaps,
however, including the lowest-order primary SSM suffices (see \citet{cenedese22a}
and \citet{cenedese22b}). A similar experience is emerging from applications
of SSMs in model-predictive control of soft robots, where a two-dimensional
SSM constructed along each spatial degree of freedom already provides
a highly accurate reduced-order model (see \citet{alora23}).
\end{rem}
The simplest case in applications is wherein $E$ is 1D ($p=1$, $q=0$),
and there are also $r$ further non-oscillatory modes of the same
stability type ($s=0$). For that case, Theorem \ref{thm:main1} takes
the following more specific form. 
\begin{prop}
\label{prop:1D real  spectral subspace}Assume that $p=1$, $q=0$
and $s=0$ hold. Then, for any integer order $\mathcal{K}\geq2$ of
approximation, the full family $\mathcal{W}\left(E\right)$ of SSMs
tangent to the spectral subspace $E$ of the nonlinear system (\ref{eq:original system})
can locally be written as 
\begin{align}
v & =\sum_{1\leq k_{1}+\sum_{\ell=1}^{r}k_{4\ell}\frac{\kappa_{\ell}}{\lambda_{1}}\leq\mathcal{K}}C_{\mathbf{k}}(u_{1})u_{1}^{k_{1}}\left|u_{1}\right|^{\sum_{\ell=1}^{r}k_{4\ell}\frac{\kappa_{\ell}}{\lambda_{1}}}+o\left(\left|u_{1}\right|^{\mathcal{K}}\right),\label{eq:1D v_SSM}
\end{align}
with the arbitrary, piecewise constant functions $C_{\mathbf{k}}(u_{1})\in\mathbb{R}^{r}$
satisfying
\[
C_{\mathbf{k}\ell}(u_{1})=\left\{ \begin{array}{c}
C_{\mathbf{k}\ell}^{+}\qquad u_{1}>0,\\
\\
C_{\mathbf{k}\ell}^{-}\qquad u_{1}\leq0,
\end{array}\right.\qquad C_{\mathbf{k}\ell}^{\pm}=0:\quad\frac{\kappa_{\ell}}{\lambda_{1}}<0\,\,\,\,\,\,\,\mathrm{or}\,\,\,\,\,\left|\mathbf{k}\right|=1,\,\,k_{4\ell}\neq1,\,\qquad\ell=1,\ldots,r.
\]
\end{prop}
\begin{proof}
See Appendix D of the \emph{Supplementary Material}.
\end{proof}
Note that for the primary (smoothest) SSM in the $\mathcal{W}\left(E\right)$
family covered by Proposition \ref{prop:1D real  spectral subspace},
there exists a unique parameter vector sequence $C_{\mathbf{k}}$
with $k_{1}\geq2$ and $\mathbf{k}_{4}=\mathbf{0}$. 

Another frequent case in practice is wherein $E$ is 2D, carries oscillations
($p=0$, $q=1$), and there are also $s$ further oscillatory modes
($r=0$). For that case, Theorem \ref{thm:main1} takes the following
more specific form.
\begin{prop}
\label{prop:2D complex spectral subspace}Assume that $p=0$, $q=1$
and $r=0$ hold. Then, for any integer order $\mathcal{K}\geq2$ of
approximation, the full family $\mathcal{W}\left(E\right)$ of SSMs
tangent to the spectral subspace $E$ of the nonlinear system (\ref{eq:original system})
can locally be written as 
\begin{align}
w & =\sum_{1\leq k_{2}+k_{3}+\sum_{m=1}^{s}\left(k_{5m}+k_{6m}\right)\frac{\beta_{m}}{\alpha_{1}}\leq\mathcal{K}}D_{\mathbf{k}}z_{1}^{k_{2}}\bar{z}_{1}^{k_{3}}\left|z_{1}\right|^{\sum_{m=1}^{s}\left(k_{5m}+k_{6m}\right)\frac{\beta_{m}}{\alpha_{1}}}e^{i\sum_{m=1}^{s}\left(k_{5m}-k_{6m}\right)\frac{\nu_{m}}{\alpha_{1}}\log\left|z_{1}\right|}\label{eq:SSM in original coordinates real-1-5-1}\\
 & +o\left(\left|z_{1}\right|^{\mathcal{K}}\right),\\
 & \,\,\,\,\,\,k_{2},k_{3}\in\mathbb{N},\quad\mathbf{k}_{5},\mathbf{k}_{6}\in\mathbb{N}^{s},\quad\mathbf{k}=\left(k_{2},k_{3},\mathbf{k}_{5},\mathbf{k}_{6}\right),\nonumber 
\end{align}
with the arbitrary constants $D_{\mathbf{k}}\in\mathbb{C}^{s}$ satisfying
\[
D_{\mathbf{k}m}=0:\quad\frac{\beta_{m}}{\alpha_{1}}<0\,\,\,\,\,\,\mathrm{or}\,\,\,\,\,\,\,\left|\mathbf{k}\right|=1,k_{5m}\neq1,\qquad m=1,\ldots,s.
\]
\end{prop}
\begin{proof}
See Appendix D of the \emph{Supplementary Material}.
\end{proof}
Note that for the primary SSM $\mathcal{W}^{\infty}(E)$ in the $\mathcal{W}\left(E\right)$
family covered by Proposition \ref{prop:2D complex spectral subspace},
there exists a unique parameter vector sequence $D_{\mathbf{k}}$
with $\left|\mathbf{k}\right|\geq2$, $\mathbf{k}_{5}=\mathbf{k}_{6}=\mathbf{0}$. 

\section{Generalized SSMs near fixed points of discrete dynamical systems\label{sec:Generalized-SSMs-for-discrete-DS}}

The techniques and prior results used in the proof of Theorem \ref{thm:main1}
are also applicable with appropriate modifications to discrete dynamical
systems defined by iterated mappings. This enables us to derive generalized
SSM results for sampled data sets of autonomous dynamical systems
and for the Poincaré maps of time-periodically forced dynamical systems.
The latter extension, in turn, implies the existence of time-periodic
generalized SSMs along periodic orbits of such systems. We only list
the related results for maps in this section with their proofs relegated
to appendices in the \emph{Supplementary Material.}

We consider an iterated nonlinear mapping of the form

\begin{equation}
x(\imath+1)=\mathcal{A}x(\imath)+f\left(x(\imath)\right),\qquad x\in\mathbb{R}^{n},\quad\mathcal{A}\in\mathbb{R}^{n\times n},\quad f=\mathcal{O}\left(\left|x\right|^{2}\right)\in C^{\infty},\qquad1\leq n<\infty,\label{eq:original system-1}
\end{equation}
and assume that its fixed point at $x=0$ is hyperbolic, i.e., the
spectrum $\mathrm{spect}\left(\mathcal{A}\right)=\left\{ \lambda_{1},\ldots,\lambda_{n}\right\} $
of $\mathcal{A}$ satisfies 
\begin{equation}
\mathrm{spect}\left(\mathcal{A}\right)\cap\left\{ \lambda\in\mathbb{C}:\,\,\left|\lambda\right|=1\right\} =\emptyset.\label{eq:hyperbolicity-2}
\end{equation}
As in the case of continuous dynamical systems, we can also allow
$f$ to be of finite smoothness if $\mathrm{spect}\left(\mathcal{A}\right)$
is fully inside or fully outside the complex unit circle (see Theorem
\ref{thm:main2}). 

We again assume that $\mathcal{A}$ is semisimple and the following
nonresonance conditions hold among its eigenvalues:
\begin{equation}
\lambda_{j}\neq\prod_{k=1}^{n}\lambda_{k}^{m_{k}},\quad m_{k}\in\mathbb{N},\quad\sum_{k=1}^{n}m_{k}\geq2,\quad j=1,\ldots,n.\label{eq:nonresonance conditions-1}
\end{equation}
As in the case of continuous dynamical systems that we have already
treated, repeated eigenvalues arising from a physical symmetry are
not excluded by conditions (\ref{eq:nonresonance conditions-1}).

We select and fix an arbitrary spectral subspace $E\subset\mathbb{R}^{n}$
and assume that the spectrum of $\mathcal{A}\vert_{E}$ contains $p$
purely real eigenvalues and $q$ pairs of complex conjugate eigenvalues.
Similarly, we assume that the remainder of the spectrum of $\mathcal{A}$
contains $r$ purely real eigenvalues and $s$ pairs of complex conjugate
eigenvalues. In that case, possibly after a linear change of coordinates,
we can assume a partition of the $x$ coordinate in the form
\begin{equation}
x=\left(u,\left(a_{1},b_{1}\right),\ldots,\left(a_{q},b_{q}\right),v,\left(c_{1},d_{1}\right),\ldots,\left(c_{s},d_{s}\right)\right)\in\mathbb{R}^{p}\times\underbrace{\mathbb{R}^{2}\times\ldots\times\mathbb{R}^{2}}_{q}\times\mathbb{R}^{r}\times\underbrace{\mathbb{R}^{2}\times\ldots\times\mathbb{R}^{2}}_{s},\label{eq:linear change of coordinates-1}
\end{equation}
with $\quad p+2q+r+2s=n,$ in which the coefficient matrix $\mathcal{A}$
of the linear part of system (\ref{eq:original system-1}) takes its
real Jordan canonical form
\begin{equation}
\mathcal{A}=\mathrm{diag}\left[A,B_{1},\ldots,B_{q},C,D_{1},\ldots,D_{s}\right],\label{eq:linsystem1-1-1}
\end{equation}
where
\begin{align}
A & =\mathrm{diag}\left[\lambda_{1},\ldots,\lambda_{p}\right]\in\mathbb{R}^{p\times p},\quad B_{k}=\left(\begin{array}{cc}
\alpha_{k} & -\omega_{k}\\
\omega_{k} & \alpha_{k}
\end{array}\right)\in\mathbb{R}^{2\times2},\quad k=1,\ldots,q,\nonumber \\
C & =\mathrm{diag}\left[\kappa_{1},\ldots,\kappa_{r}\right]\in\mathbb{R}^{r\times r},\quad D_{m}=\left(\begin{array}{cc}
\beta_{m} & -\nu_{m}\\
\nu_{m} & \beta_{m}
\end{array}\right)\in\mathbb{R}^{2\times2},\quad m=1,\ldots,s.\label{eq:linsystem1-2-1}
\end{align}
To simplify our final formulas, we also assume that 
\begin{equation}
\kappa_{\ell}>0,\quad\ell=1,\ldots,r,\label{eq:positivre real part}
\end{equation}
which means that the linearized mapping is orientation preserving
along all of the $v_{\ell}$ directions. 

Using again the complexified variables $z=a+ib\in\mathbb{C}^{q}$
and $w=c+id\in\mathbb{C}^{s}$, as well as the piecewise constant
functions (\ref{eq:piecewise constant functions}), we define the
discrete analogues of formulas (\ref{eq:complex SSM 1}) as
\begin{align}
\mathcal{V}_{\ell}(u,z) & =\sum_{j=1}^{p}K_{\ell j}\left(u_{j}\right)\left|u_{j}\right|^{\frac{\log\kappa_{\ell}}{\log\left|\lambda_{j}\right|}}+\sum_{k=1}^{q}L_{\ell k}\left|z_{k}\right|^{\frac{\log\kappa_{\ell}}{\log\sqrt{\alpha_{k}^{2}+\omega_{k}^{2}}}},\quad\ell=1,\ldots,r,\nonumber \\
\mathcal{\mathcal{E}}_{m}(u,z) & =\left[\sum_{j=1}^{p}O_{mj}\left(u_{j}\right)\left|u_{j}\right|^{\frac{\log\sqrt{\beta_{m}^{2}+\nu_{m}^{2}}}{\log\left|\lambda_{j}\right|}}+\sum_{k=1}^{q}Q_{mk}\left|z_{k}\right|^{\frac{\log\sqrt{\beta_{m}^{2}+\nu_{m}^{2}}}{\log\sqrt{\alpha_{k}^{2}+\omega_{k}^{2}}}}\right]\nonumber \\
 & \,\,\,\,\,\,\times e^{i\left(\frac{1}{p+q}\arctan\left(\frac{\nu_{m}}{\beta_{m}}\right)\left[\sum_{j=1}^{p}\frac{\log\left|u_{j}\right|}{\log\left|\lambda_{j}\right|}+\sum_{k=1}^{q}\frac{\log\left|z_{k}\right|}{\log\sqrt{\alpha_{k}^{2}+\omega_{k}^{2}}}\right]\right)},\quad m=1,\ldots,s,\label{eqcomplex SSM 1-2-1}
\end{align}
 where 
\begin{align}
K_{\ell j} & (u_{j})\equiv0,\quad\mathrm{if}\quad\frac{\log\kappa_{\ell}}{\log\left|\lambda_{j}\right|}<0,\qquad O_{mj}(u_{j})\equiv0,\quad\mathrm{if}\quad\frac{\log\sqrt{\beta_{m}^{2}+\nu_{m}^{2}}}{\log\left|\lambda_{j}\right|}<0,\nonumber \\
L_{\ell k} & =0,\quad\mathrm{if}\quad\frac{\log\kappa_{\ell}}{\log\sqrt{\alpha_{k}^{2}+\omega_{k}^{2}}}<0,\qquad Q_{mk}=0,\quad\mathrm{if}\quad\frac{\log\sqrt{\beta_{m}^{2}+\nu_{m}^{2}}}{\log\sqrt{\alpha_{k}^{2}+\omega_{k}^{2}}}<0.\label{eq:general_fractional_SSM2-2}
\end{align}
As in the case of Theorem \ref{thm:main1}, the function families
(\ref{eqcomplex SSM 1-2-1}) turn out to describe all continuous invariant
graphs over $E$ in the linear part of the mapping (\ref{eq:original system-1}).
Our main result on generalized SSMs for the full nonlinear mapping
can then be stated as follows.
\begin{thm}
\label{thm:main2}
\end{thm}
\begin{description}
\item [{\emph{(i)}}] For any integer order $\mathcal{K}\geq2$ of approximation\emph{
there exists a unique set of coefficients $h_{\mathbf{k}}\in\mathbb{C}^{r+s}$
such that the full family $\mathcal{W}\left(E\right)$ of SSMs tangent
to the spectral subspace $E$ of the nonlinear mapping (\ref{eq:original system-1})
can locally be written as 
\begin{align}
\left(\begin{array}{c}
v\\
w
\end{array}\right) & =\mathcal{G}(u,z)=\left(\begin{array}{c}
\mathcal{V}(u,z)\\
\mathcal{\mathcal{E}}(u,z)
\end{array}\right)+\sum_{2\leq\mathrm{monomial}\mathrm{\,order}\leq\mathcal{K}}h_{\mathbf{k}}u^{\mathbf{k}_{1}}z^{\mathbf{k}_{2}}\bar{z}^{\mathbf{k}_{3}}\mathcal{V}^{\mathbf{k}_{4}}\mathcal{\mathcal{E}}^{\mathbf{k}_{5}}\bar{\mathcal{\mathcal{E}}}^{\mathbf{k}_{6}}+o\left(\left|(u,z)\right|^{\mathcal{K}}\right),\label{eq:SSM in original coordinates real-1-5-3}\\
 & \,\,\,\,\,\,\mathbf{k}_{1}\in\mathbb{N}^{p},\quad\mathbf{k}_{2},\mathbf{k}_{3}\in\mathbb{N}^{q},\quad\mathbf{k}_{4}\in\mathbb{N}^{r},\quad\mathbf{k}_{5},\mathbf{k}_{6}\in\mathbb{N}^{s},\quad\mathbf{k}=\left(\mathbf{k}_{1},\mathbf{k}_{2},\mathbf{k}_{3},\mathbf{k}_{4},\mathbf{k}_{5},\mathbf{k}_{6}\right),\nonumber 
\end{align}
where $\mathcal{V}(u,z)\in\mathbb{R}^{r}$ and $\mathcal{\mathcal{E}}(u,z)\in\mathbb{C}^{s}$
are defined via formulas (}\ref{eqcomplex SSM 1-2-1}\emph{), with
$K_{\ell j}$, $O_{mj}$, $L_{\ell k}$ and $Q_{mk}$ satisfying the
constraints (}\ref{eq:general_fractional_SSM2-2}\emph{) but otherwise
chosen arbitrarily. }
\item [{(ii)}] The reduced dynamics on the individual members of the family
$\mathcal{W}\left(E\right)$ are given by
\begin{equation}
\left(\begin{array}{c}
u(\imath+1)\\
z(\imath+1)
\end{array}\right)=\left(\begin{array}{cccc}
A & 0 & 0 & 0\\
0 & B_{1} & 0 & 0\\
0 & 0 & \ddots & 0\\
0 & 0 & 0 & B_{q}
\end{array}\right)\left(\begin{array}{c}
u(\imath)\\
z(\imath)
\end{array}\right)+f^{\left(u,z\right)}\left(u(\imath),z(\imath),\mathcal{G}(u(\imath),z(\imath))\right),\label{eq:reduced dynamics discrete}
\end{equation}
where $f^{\left(u,z\right)}$ denotes the $(u,z)$ coordinate component
of $f=\left(f^{\left(u,z\right)},f^{v},f^{w}\right)$. 
\item [{\emph{(iii)}}] \emph{If 
\[
\frac{\log\kappa_{\ell}}{\log\left|\lambda_{j}\right|}<0,\qquad\frac{\log\sqrt{\beta_{m}^{2}+\nu_{m}^{2}}}{\log\left|\lambda_{j}\right|}<0,\qquad\frac{\log\kappa_{\ell}}{\log\sqrt{\alpha_{k}^{2}+\omega_{k}^{2}}}<0,\qquad\frac{\log\sqrt{\beta_{m}^{2}+\nu_{m}^{2}}}{\log\sqrt{\alpha_{k}^{2}+\omega_{k}^{2}}}<0
\]
hold for all values of the indices $j,k,\ell$ and $m$, then $\mathcal{W}\left(E\right)$
is unique and of class $C^{\infty}$. Otherwise, $\mathcal{W}\left(E\right)$
is a family of non-unique SSMs whose typical member is of smoothness
class $C^{\eta}$ with
\begin{equation}
\eta=\mathrm{Int}\left[\min_{j,k,\ell,m}^{+}\left\{ \frac{\log\kappa_{\ell}}{\log\left|\lambda_{j}\right|},\quad\frac{\log\sqrt{\beta_{m}^{2}+\nu_{m}^{2}}}{\log\left|\lambda_{j}\right|},\quad\frac{\log\kappa_{\ell}}{\log\sqrt{\alpha_{k}^{2}+\omega_{k}^{2}}},\quad\frac{\log\sqrt{\beta_{m}^{2}+\nu_{m}^{2}}}{\log\sqrt{\alpha_{k}^{2}+\omega_{k}^{2}}}\right\} \right],\label{eq:eta formula for smoothness - discrete}
\end{equation}
where the $\mathrm{Int}$ refers to the integer part of a positive
number and $\min^{+}$ refers to a minimum of the positive elements
of a list of numbers that follow. Even in that case, however, there
exists a unique member $\mathcal{W}^{\infty}(E)$ of the SSM family
that is of class $C^{\infty}$ and satisfies $K_{\ell j}^{\pm}\equiv O_{mj}^{\pm}=L_{\ell k}=Q_{mk}=0$
for all $j,k,\ell$ and $m$ in eq. (}\ref{eqcomplex SSM 1-2-1}). 
\item [{\emph{(iv)}}] \emph{If system (\ref{eq:original system-1}) is
$C^{R}$ in a parameter vector $\mu\in\mathbb{R}^{P}$, then $\mathcal{W}\left(E\right)$
is also $C^{R}$ in $\mu$. }
\item [{(v)}] If the function $f$ only has finite smoothness, i.e., $f=\mathcal{O}\left(\left|x\right|^{2}\right)\in C^{\rho}$
and $\mathrm{\log\left|\lambda_{j}\right|}$ has the same nonzero
sign for all $j=1,\ldots,N$, then statements (i)-(iv) still hold
for all $\mathcal{K}\leq\rho$. 
\item [{\emph{(vi)}}] If the function $f$ is analytic in a neighborhood
of the origin and $\mathrm{\log\left|\lambda_{j}\right|}$has the
same nonzero sign for all $j=1,\ldots,N$, then the series expansion
(\ref{eq:SSM in original coordinates real-1-5}) converges near the
origin as $\mathcal{K}\to\infty$. 
\end{description}
\begin{proof}
See Appendix E of the \emph{Supplementary Material}.
\end{proof}
Appropriately rephrased versions of Remarks \ref{rem: explanation on obtaining primary SSMs}-\ref{rem:uniqueness of manifolds}
continue to apply in the present, discrete setting. We also add the
following remark that shows how Theorem \ref{thm:main2} can be used
to extend the results of Theorem \ref{thm:main1} to small, time-periodic
perturbations of the autonomous system (\ref{eq:original system}).
\begin{rem}
\label{rem: extension of Thm. 1 to time-periodic perturbations} {[}\textbf{Extension
to time-periodically forced continuous dynamical systems}{]} Consider
a small, temporally $T$-periodic perturbation of the system of ODEs
(\ref{eq:original system}) in the form 
\begin{equation}
\dot{x}=\mathcal{A}x+f(x)+\epsilon g(x,t),\quad f=\mathcal{O}\left(\left|x\right|^{2}\right)\in C^{\infty},g\in C^{\infty},\quad g(x,t+T)=g(x,t),\quad0\leq\epsilon\ll1.\label{eq:original system-2}
\end{equation}
For $\epsilon=0$, if the conditions of Theorem \ref{thm:main1} are
satisfied, then the resulting SSM family $\mathcal{W}\left(E\right)$
also acts as an SSM family for the period-$T$ Poincaré map defined
for the $\epsilon=0$ limit of system (\ref{eq:original system-2}).
By statement (iv) of Theorem \ref{thm:main2}, the SSMs of this Poincaré
map are smooth functions of $\epsilon$ and hence will smoothly persist
for the Poincaré map of system (\ref{eq:original system-2}) for small
$\epsilon>0$ parameter values. This in turn allows us to conclude
the existence of $T$-periodic SSMs for the full dynamical system
(\ref{eq:original system-2}) for small enough $\epsilon>0$. These
SSMs and their reduced dynamics can the sought in the form of a Taylor-expansion
in $\epsilon$ combined with a Fourier expansion in time, as explained
and demonstrated for primary, non-mixed-mode SSMs by \citet{haller16,breunung2018,ponsioen2020}. 
\end{rem}
Next, we also restate the appropriately modified versions of Propositions
\ref{prop:1D real  spectral subspace}-\ref{prop:2D complex spectral subspace}
for maps, whose proofs we write out for completeness in Appendix F
of the Supplementary Material.
\begin{prop}
\label{prop:1D real  spectral subspace-1}Assume that $p=1$, $q=0$
and $s=0$ hold. Then, for any integer order $\mathcal{K}\geq2$ of
approximation, the full family $\mathcal{W}\left(E\right)$ of SSMs
tangent to the spectral subspace $E$ of the nonlinear mapping (\ref{eq:original system-1})
can locally be written as 
\begin{align}
v & =\sum_{1\leq k_{1}+\sum_{\ell=1}^{r}k_{4\ell}\frac{\log\kappa_{\ell}}{\log\left|\lambda_{1}\right|}\leq\mathcal{K}}C_{\mathbf{k}}(u_{1})u_{1}^{k_{1}}\left|u_{1}\right|^{\sum_{\ell=1}^{r}k_{4\ell}\frac{\log\kappa_{\ell}}{\log\left|\lambda_{1}\right|}}+o\left(\left|u_{1}\right|^{\mathcal{K}}\right),\label{eq:1D v_SSM-1}
\end{align}
with the arbitrary, piecewise constant functions $C_{\mathbf{k}}(u_{1})\in\mathbb{R}^{r}$
satisfying
\[
C_{\mathbf{k}\ell}(u_{1})=\left\{ \begin{array}{c}
C_{\mathbf{k}\ell}^{+}\qquad u_{1}>0,\\
\\
C_{\mathbf{k}\ell}^{-}\qquad u_{1}\leq0,
\end{array}\right.\qquad C_{\mathbf{k}\ell}^{\pm}=0:\quad\frac{\log\kappa_{\ell}}{\log\left|\lambda_{1}\right|}<0\,\,\,\,\,\,\,\mathrm{or}\,\,\,\,\,\left|\mathbf{k}\right|=1,k_{4\ell}\neq1,\,\qquad\ell=1,\ldots,r.
\]
\end{prop}
~~~~~
\begin{prop}
\label{prop:2D complex spectral subspace-1}Assume that $p=0$, $q=1$
and $r=0$ hold. Then, for any integer order $\mathcal{K}\geq2$ of
approximation, the full family $\mathcal{W}\left(E\right)$ of SSMs
tangent to the spectral subspace $E$ of the nonlinear system (\ref{eq:original system})
can locally be written as 
\begin{align}
w & =\sum_{1\leq k_{2}+k_{3}+\Xi\left(\mathbf{k}_{5},\mathbf{k}_{6}\right)\leq\mathcal{K}}D_{\mathbf{k}}z_{1}^{k_{2}}\bar{z}_{1}^{k_{3}}\left|z_{1}\right|^{\Xi\left(\mathbf{k}_{5},\mathbf{k}_{6}\right)}e^{i\Gamma\left(\mathbf{k}_{5},\mathbf{k}_{6}\right)\log\left|z_{1}\right|}+o\left(\left|z_{1}\right|^{\mathcal{K}}\right)\label{eq:SSM in original coordinates real-1-5-1-1}\\
 & \,\,\,\,\,\,k_{2},k_{3}\in\mathbb{N},\quad\mathbf{k}_{5},\mathbf{k}_{6}\in\mathbb{N}^{s},\quad\mathbf{k}=\left(k_{2},k_{3},\mathbf{k}_{5},\mathbf{k}_{6}\right),\nonumber \\
 & \Xi\left(\mathbf{k}_{5},\mathbf{k}_{6}\right)=\sum_{m=1}^{s}\left(k_{5m}+k_{6m}\right)\frac{\log\sqrt{\beta_{m}^{2}+\nu_{m}^{2}}}{\log\sqrt{\alpha_{1}^{2}+\omega_{1}^{2}}},\quad\Gamma\left(\mathbf{k}_{5},\mathbf{k}_{6}\right)=\frac{\sum_{m=1}^{s}\left(k_{5m}-k_{6m}\right)\arctan\left(\frac{\nu_{m}}{\beta_{m}}\right)}{\log\sqrt{\alpha_{1}^{2}+\omega_{1}^{2}}},
\end{align}
 with the arbitrary constants $D_{\mathbf{k}}\in\mathbb{C}^{s}$ satisfying
\[
D_{\mathbf{k}m}=0:\quad\frac{\log\sqrt{\beta_{m}^{2}+\nu_{m}^{2}}}{\log\sqrt{\alpha_{1}^{2}+\omega_{1}^{2}}}<0\,\,\,\,\,\,\mathrm{or}\,\,\,\,\,\,\,\left|\mathbf{k}\right|=1,k_{5m}\neq1,\qquad m=1,\ldots,s.
\]
\end{prop}

\section{Normal forms on 2D fractional SSMs}

It is often advantageous to bring the dynamics on a generalized SSM
to a normal form. Normal form transformations simplify the reduced
SSM dynamics (\ref{eq:reduced dynamics continuous}) and (\ref{eq:reduced dynamics discrete})
by removing as many terms from these expansions as possible via near-identity
changes of coordinates. These transformations have been well-understood
for primary SSMs on which the reduced dynamics are expressible via
classical Taylor expansions (see \citet{ponsioen2018,jain2022,cenedese22a,cenedese22b}). 

This classical normal form procedure (see, e.g., \citet{arnold83,guckenheimer83})
is also applicable to mixed-mode primary SSMs, as those have no fractional-powered
terms in their expansion either. For that reason, available equation-
and data-driven model packages, such as SSMTool and SSMLearn, can
directly be used to compute mixed-mode SSMs and normal forms of their
reduced dynamics, even though the mixed-mode SSM theory described
here in Sections \ref{sec:Generalized-SSMs-for-continuous-DS}-\ref{sec:Generalized-SSMs-for-discrete-DS}
was not available at the time of the development of these packages.
Normal forms on fractional SSMs, however, have not yet been treated
in the literature and require a separate discussion here.

Each normal form transformation is a near-identity change of coordinates
that simplifies a dynamical system near a fixed point up to a given
polynomial order. The simplified dynamical system, however, is only
topologically conjugate to the original system on a domain on which
the transformation is at least continuously invertible. As a consequence,
each normal form transformation limits the ability of the normalized
system to capture nonlinear dynamical features away from the origin.
A clear manifestation of this limitation is the smooth linearization
of \citet{sternberg58} used in proving Theorems \ref{thm:main1}-\ref{thm:main2}
of this paper. Indeed, this linearization result succeeds in constructing
a $C^{\infty}$ normal form transformation that removes\emph{ all
}nonlinear terms from the polynomial expansion of the right-hand side
of the dynamical system. At the same time, the transformation only
converges on a domain around the fixed point in which the original
system exhibits linearizable behavior, completely missing all interesting
dynamical phenomena that arise from the coexistence of several isolated
invariant sets. For this reason, we do not advocate linearization
for model reduction and only use linearization to identify the local
signatures of SSMs which we then study globally without linearization. 

More generally, we advocate a sparing use of normal form transformations
for SSM-reduced models in order to preserve the predictive power of
SSM-reduced dynamics on larger domains around the origin. To simplify
our discussion, we only cover here the setting of 2D underdamped,
oscillatory SSMs on which low-order normal forms have proven to be
very effective in predicting coexisting, nontrivial steady states
(see \citet{ponsioen2019,ponsioen2020,jain2022,cenedese22a,cenedese22b}).
These \emph{extended normal forms} only seek a partial (as opposed
to maximal possible) simplification of the dynamics. Namely, they
do not remove near-resonant terms that would be in principle removable
but would result in small denominators that seriously limit the domain
of invertibility of the normal form transformation (see \citet{cenedese22a}
for a detailed discussion for data-driven SSMs). Specifically, in
the setting of Proposition \ref{prop:2D complex spectral subspace},
we obtain the following result. 
\begin{thm}
\label{prop:2D complex normal form}Assume that $p=0$, $q=1$ and
$r=0$ hold. Then:
\end{thm}
\begin{description}
\item [{(i)}] For any integer order $\mathcal{K}\geq2$ of approximation\emph{
there exists a near-identity change of coordinates that transforms
the 2D reduced dynamics on each member of the SSM family $\mathcal{W}\left(E\right)$
to the form
\begin{align}
\dot{\xi} & =\left(\alpha_{1}+i\beta_{1}\right)\xi+\sum_{1\leq\mathrm{monomial}\mathrm{\,order}\leq\mathcal{K}}H_{\mathbf{k}}\xi_{1}^{k_{2}}\bar{\xi}_{1}^{k_{2}-1}\left|\xi_{1}\right|^{\sum_{m=1}^{s}\left(k_{5m}+k_{6m}\right)\frac{\beta_{m}}{\alpha_{1}}}e^{i\sum_{m=1}^{s}\left(k_{5m}-k_{6m}\right)\frac{\nu_{m}}{\alpha_{1}}\log\left|\xi_{1}\right|}\label{eq:SSM in original coordinates real-1-5-1-2}\\
 & \,\,\,\,\,+o\left(\left|\xi_{1}\right|^{\mathcal{K}}\right),\nonumber \\
 & \,\,\,\,\,\,k_{2}\in\mathbb{N},\quad\mathbf{k}_{5},\mathbf{k}_{6}\in\mathbb{N}^{s},\quad\mathbf{k}=\left(k_{2},\mathbf{k}_{5},\mathbf{k}_{6}\right),\nonumber 
\end{align}
with the constants $H_{\mathbf{k}}\in\mathbb{C}$ satisfying 
\[
H_{\mathbf{k}}=0:\quad\frac{\beta_{m}}{\alpha_{1}}<0\,\,\,\,\,\,\mathrm{or}\,\,\,\,\,\,\,\left|\mathbf{k}\right|=1,k_{5m}\neq1,\qquad m=1,\ldots,s,
\]
The coordinate change is of class $C^{\infty}$ if $\min_{m}\left(\beta_{m}/\alpha_{1}\right)<0$
and of class $C^{\mathrm{Int}\left[\min_{m}\left(\beta_{m}/\alpha_{1}\right)\right]}$
if $\min_{m}\left(\beta_{m}/\alpha_{1}\right)>0$. }
\item [{(iii)}] \emph{Assume further that $s=1$ and the relation
\begin{equation}
\frac{1}{2}<\frac{\beta_{1}}{\alpha_{1}}<2\label{eq:beta over alpha assumption-1}
\end{equation}
holds between the decay exponent $\beta_{1}$ outside $E$ and the
decay exponent $\alpha_{1}$ inside $E$. Then in polar coordinates
$(r,\phi)$ defined via the relation $z=re^{i\phi}$, a truncation
of the normal form (\ref{eq:SSM in original coordinates real-1-5-1-2})
that incorporates all possible terms below $\mathcal{K}=3$ (and possibly
more, depending on the exact value of $\beta_{1}/\alpha_{1}$) can
be written as 
\begin{align}
\dot{r} & =\alpha_{1}r+Ar^{3}+P_{1}r^{\frac{\beta_{1}}{\alpha_{1}}+1}\sin2\left[Q+\frac{\nu_{1}}{\alpha_{1}}\log r\right]+P_{2}r^{\frac{2\beta_{1}}{\alpha_{1}}+1}+P_{3}r^{\frac{2\beta_{1}}{\alpha_{1}}+1}\sin2\left[Q+\frac{\nu_{1}}{\alpha_{1}}\log r\right],\nonumber \\
\dot{\phi} & =\omega_{1}+Br^{2}+R_{1}r^{\frac{\beta_{1}}{\alpha_{1}}}\sin2\left[Q+\frac{\nu_{1}}{\alpha_{1}}\log r\right]+R_{2}r^{\frac{2\beta_{1}}{\alpha_{1}}}+R_{3}r^{\frac{2\beta_{1}}{\alpha_{1}}}\sin2\left[Q+\frac{\nu_{1}}{\alpha_{1}}\log r\right],\label{eq:final polar normal form-1}
\end{align}
for appropriate constants $A,$$B$, $P_{j}$ ,$Q$ and $R_{j}$ that
depend on both the original full system and the specific SSM chosen. }
\end{description}
\begin{proof}
See Appendix G of the \emph{Supplementary Material}.
\end{proof}
\begin{rem}
In nonlinear vibration studies, the \emph{backbone curve} for a given
mode is the instantaneous frequency viewed as a function of the instantaneous
amplitude. From the polar normal form (\ref{eq:final polar normal form-1}),
this backbone curve can be approximated up to cubic order by the formula
\begin{equation}
\Omega(r)=\omega_{1}+Br^{2}+R_{1}r^{\frac{\beta_{1}}{\alpha_{1}}}\sin2\left[Q+\frac{\nu_{1}}{\alpha_{1}}\log r\right]+R_{2}r^{\frac{2\beta_{1}}{\alpha_{1}}}+R_{3}r^{\frac{2\beta}{\alpha}}\sin2\left[Q+\frac{\nu_{1}}{\alpha_{1}}\log r\right].\label{eq:Omega(r)}
\end{equation}
Note that on the unique smoothest (or primary) SSM, only the first
two terms of this expression can have nonzero coefficients. 
\end{rem}
\begin{rem}
Another customary notion in nonlinear vibration studies is the \emph{instantaneous
damping curve}, defined as quotient of the rate of change of the instantaneous
amplitude and the instantaneous amplitude itself. From formula (\ref{eq:final polar normal form-1}),
we obtain this curve up to cubic order in the following form:
\begin{equation}
\kappa(r)=\alpha_{1}+Ar^{2}+P_{1}r^{\frac{\beta_{1}}{\alpha_{1}}}\sin2\left[Q+\frac{\nu_{1}}{\alpha_{1}}\log r\right]+P_{2}r^{\frac{2\beta_{1}}{\alpha_{1}}}+P_{3}r^{\frac{2\beta_{1}}{\alpha_{1}}}\sin2\left[Q+\frac{\nu_{1}}{\alpha_{1}}\log r\right].\label{eq:kappa(r)}
\end{equation}
Again, on the primary SSM, only the first two terms $\kappa(r)$ are
nonzero; all other terms appear only on secondary SSMs. 
\end{rem}
\begin{rem}
Formulas (\ref{eq:Omega(r)})-(\ref{eq:kappa(r)}) remain valid for
the slowest, 2D SSM family of a general, multi-degree-of-freedom dynamical
system ($n\geq2)$ as long as the relation $\frac{1}{2}<\frac{\beta_{1}}{\alpha_{1}}<2$
holds between the decay rates of the first two modes, and all other
modes decay at rates satisfying $\beta_{m}\geq\beta_{1}$ for $m=2,\ldots,s.$
This is because under this relation, further factional terms are too
high in order to enter the normal form truncated at cubic order.
\end{rem}

\section{Example 1: Data-driven fractional SSM for transitions in Couette
flow\label{sec:Couette flow example}}

\begin{figure}[H]
\begin{centering}
\includegraphics[width=1\textwidth]{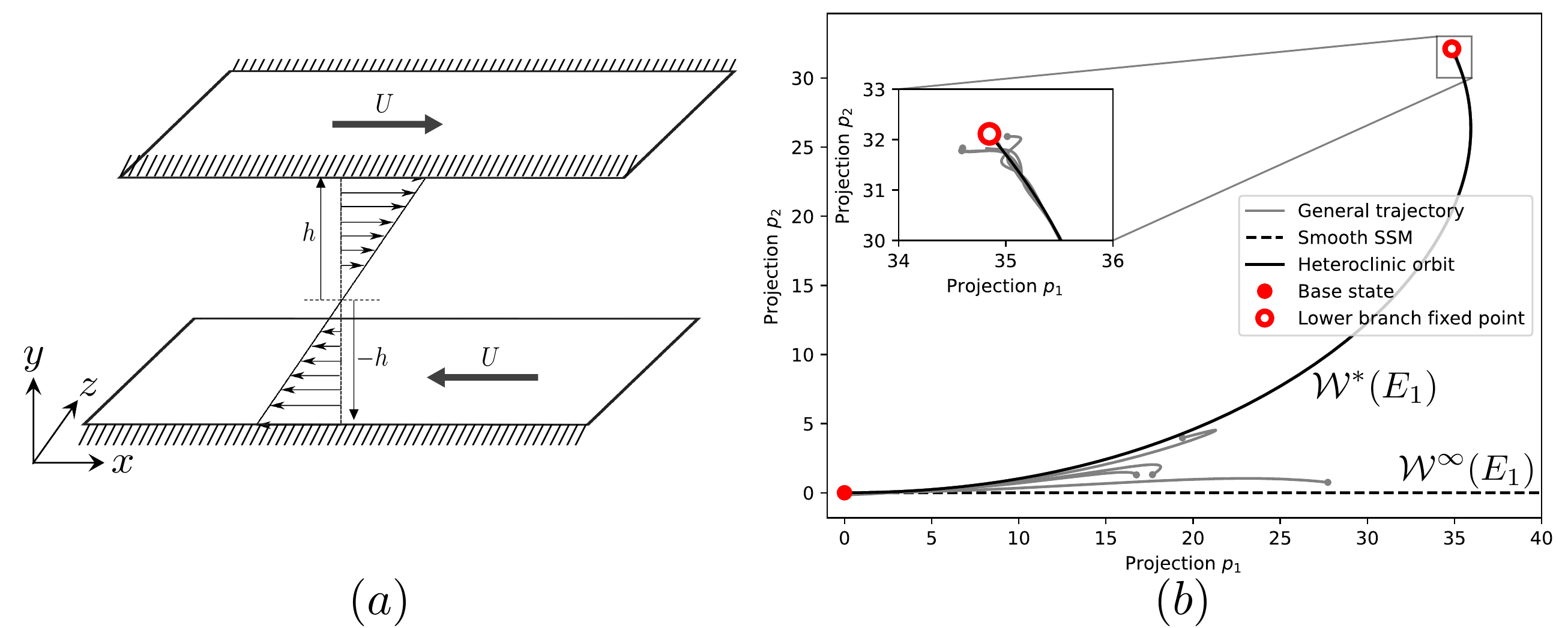}
\par\end{centering}
\caption{(a) The geometry of planar Couette flow. (b) Projection of the full
Couette flow dynamics onto the first and second POD modes of the dynamics
on the spectral subspace $E_{1}$ for $Re=134.52$. Solid red dot:
stable base state; hollow red dot: unstable equilibrium state; black
curve: 1D fractional SSM, $W^{*}\left(E_{1}\right)$ connecting the
two steady states; dashed black line: primary (analytic) slow SSM
$\mathcal{W^{\infty}}(E_{1})$ coinciding with $E_{1}$; grey curves:
nearby trajectories converging to the base state. \label{fig:primary SSM in Couette flow}}
\end{figure}

As a first example, we consider a plane Couette flow between two plates
that move in opposite directions with velocity $U$, as shown in Fig.
\ref{fig:primary SSM in Couette flow}a.
As already mentioned in Section \ref{sec:Two-motivating-examples},
\citet{page19} showed in detail the inevitable failure of linear
model reduction methods to capture transitions among coexisting steady
states in this flow. In contrast, \citet{kaszas22} used data-driven
SSM reduction to construct 1D and 2D models for such transitions.
These models were obtained from an integer-powered polynomial approximation
of the SSM. Here we will show that lower-order polynomial approximation
to these manifolds are more accurate if we account for fractional-powered
terms in the data-driven construction of the SSM.

The Navier\textendash Stokes equations governing the evolution of
the velocity field $\mathbf{u}(x,y,z,t)$ defined over the domain
$(x,y,z)\in\Omega=[0,L_{x}]\times[-h,h]\times[0,L_{z}]$ can be written
as 

\begin{align}
\frac{\partial\mathbf{u}}{\partial t}+(\mathbf{u}\nabla)\mathbf{u} & =-\nabla p+\frac{1}{\text{Re}}\Delta\mathbf{u},\label{eq:fullNS}\\
\nabla\cdot\mathbf{u} & =0,\nonumber 
\end{align}
where $\mathrm{Re}=Uh/\nu$ denotes the Reynolds number with viscosity
$\nu$ and $p$ is the pressure. To eliminate the pressure from (\ref{eq:fullNS}),
we recall the Helmholtz-decomposition to write\textbf{ $\mathbf{u}$}
as
\[
\mathbf{u}=\nabla\Phi+\mathbf{q},
\]
where $\Phi$ is a scalar field and $\mathbf{q}$ is divergence free.
The Leray projection of $\mathbf{u}$ , denoted by $\mathbb{P}[\mathbf{u}]$,
extracts the divergence-free component of $\mathbf{u}$, i.e., returns
$\mathbb{P}[\mathbf{u}]=\mathbf{q}.$ With this notation, the Leray
projection of \eqref{eq:fullNS} becomes 

\begin{equation}
\frac{\partial\mathbf{u}}{\partial t}=-\mathbb{P}\left[(\mathbf{u}\nabla)\mathbf{u}\right]+\frac{1}{\text{Re}}\Delta\mathbf{u}.\label{eq:projectedNS}
\end{equation}

Denoting perturbations from the laminar solution $\mathbf{U}(y)=(y,0,0)$
by $\mathbf{v}$, we can rewrite the Leray-projected Navier-Stokes
equation (\ref{eq:projectedNS}) as 
\begin{equation}
\frac{\partial\mathbf{v}}{\partial t}=\frac{1}{\text{Re}}\Delta\mathbf{v}-\mathbb{P}\left[\left(\mathbf{U}\nabla\right)\mathbf{v}\right]-\mathbb{P}\left[\left(\mathbf{v}\nabla\right)\mathbf{U}\right]+\mathbb{P}\left[\left(\mathbf{v}\nabla\right)\mathbf{v}\right],\label{eq:perturbationNS}
\end{equation}
 whose linear part, 
\begin{equation}
\frac{\partial\mathbf{v}}{\partial t}=\frac{1}{\text{Re}}\Delta\mathbf{v}-\mathbb{P}\left[\left(\mathbf{U}\nabla\right)\mathbf{v}\right]-\mathbb{P}\left[\left(\mathbf{v}\nabla\right)\mathbf{U}\right],\label{eq:orr-sommerfeld}
\end{equation}
 is known as the Orr\textendash Sommerfeld equation. Assuming periodicity
in $x$ and $z$ allows us to expand $\mathbf{v}$ into Fourier-series
as
\[
\mathbf{v}(x,y,z,t)=\sum_{n,m}\mathbf{v}(y)e^{\lambda t}e^{i\alpha_{x}nx+i\alpha_{z}mz},
\]
where the fundamental wave numbers are defined as $\alpha_{x}=\frac{2\pi}{L_{x}}$
and $\alpha_{z}=\frac{2\pi}{L_{z}}$. This turns the linearized equation
(\ref{eq:orr-sommerfeld}) into an eigenvalue problem, with the eigenvalue
denoted by $\lambda$.\textbf{ }Using the explicit expressions of
the differential operators involved (see, e.g., \citet{schmid01}),
one finds that the least stable mode is independent of $x$ and has
eigenvalues and eigenfunctions of the form 

\begin{equation}
\lambda_{1}=-\frac{\pi^{2}+4\alpha_{z}^{2}}{4\text{Re}}\qquad\mathbf{v}_{1}=\sin(\alpha_{z}z)\left(\begin{array}{c}
\sin(\sqrt{\lambda_{1}\text{Re}}\,y)\\
0\\
0
\end{array}\right).\label{eq:orrsommerfeld smallest}
\end{equation}
For this mode, a direct calculation shows that 
\[
\mathbb{P}\left[\left(\mathbf{v}_{1}\nabla\right)\mathbf{v}_{1}\right]=\mathbf{0},
\]
i.e., the nonlinear term in eq. (\ref{eq:perturbationNS}) vanishes
on the slowest eigenfunction. Consequently, 
\begin{equation}
\mathcal{W^{\infty}}(E_{1})=E_{1}=\left\{ \mathbf{U}(y)+\tau\mathbf{v}_{1}(y,z),\tau\in\mathbb{R}\right\} \label{eq:couette SSM}
\end{equation}
is the primary, slow SSM of the laminar state. 

Indeed, substituting an integer-powered polynomial expansion for an
invariant manifold tangent to $E_{1}$ into the nonlinear equation
(\ref{eq:perturbationNS}) yields vanishing Taylor coefficients at
all orders. In other words, such an expansion trivially converges
and yields the single flat, analytic SSM shown as a dashed line in
Fig. \ref{fig:primary SSM in Couette flow}b. The dynamics on this
primary SSM is governed by the linear system (\ref{eq:orr-sommerfeld}),
and hence no connection from the laminar base state $\mathbf{U}(y)$
to any nontrivial steady state can be contained in $\mathcal{W^{\infty}}(E_{1})$.
As a result, the transition to the base state to the other nontrivial
steady state in system \eqref{eq:perturbationNS} must be a fractional
SSM, $\mathcal{W}^{*}\left(E_{1}\right)$, anchored at the base state,
just as in our simple example \eqref{eq:example}. We confirm this
numerically in Fig. \ref{fig:primary SSM in Couette flow}b, whose
geometry bears a close similarity to that of our simple example in
Fig. \ref{fig:Simple_examples}a.

The fractional SSM $\mathcal{W}^{*}\left(E_{1}\right)$ in Fig. \ref{fig:primary SSM in Couette flow}b
has limited differentiability. Yet a data-driven global approximation
for this SSM still gives satisfactory results for high enough orders
of regression, as shown by \citet{kaszas22}. This is because every
continuous curve can be arbitrarily closely approximated by smooth
polynomials (\citet{rudin76}). The required order for such an approximation,
however, is generally higher than necessary and leads to larger errors
further away from the origin. 

Here, we will seek the globally defined manifold $\mathcal{W}^{*}\left(E_{1}\right)$
in the asymptotic form \eqref{eq:1D v_SSM-1} which has the same approximation
power as polynomial expansions but also has the correct local form
and correct degree of smoothness near the origin. We show that a data-driven
fit of this exact fractional form for $\mathcal{W}^{*}\left(E_{1}\right)$
is indeed more accurate at a given order than the one with purely
integer-powered polynomials.

We use \emph{Channelflow }of \citet{gibson18} with a spectral discretization
of $n=\mathcal{O}\left(10^{5}\right)$ for the nonlinear system \eqref{eq:perturbationNS}.
For the\textbf{ }parameter values chosen ($L_{x}=5\pi/2$, $L_{z}=4\pi/3$),
we have $p=1$, $q=0$ and $r+s=n-1$ in terms of the notation used
in Section \ref{sec:Generalized-SSMs-for-continuous-DS}. We use the
discrete formulation of our results from Section \ref{sec:Generalized-SSMs-for-discrete-DS}
with an integer time sampling, which \citet{kaszas22} found to give
more robust reduced models than the continuous time formulation. The
eigenvalues of the linearized map $\lambda_{\mathrm{map}}$ and the
linearized flow $\lambda_{\text{flow}}$ at the base state are related
to each other via $\lambda_{\mathrm{map}}=e^{\lambda_{\text{\ensuremath{\mathrm{flow}}}}}$
.

For these parameters, formulas \eqref{eq:1D v_SSM-1} give the fractional
SSM parametrization
\begin{align}
v & =\sum_{1\leq k_{1}+\sum_{\ell=1}^{r}k_{4\ell}\frac{\log\kappa_{\ell}}{\log\left|\lambda_{1}\right|}\leq\mathcal{K}}C_{\mathbf{k}}u_{1}^{k_{1}+\sum_{\ell=1}^{r}k_{4\ell}\frac{\log\kappa_{\ell}}{\log\left|\lambda_{1}\right|}}+o\left(\left|u_{1}\right|^{\mathcal{K}}\right),\qquad\mathbf{k}=\left(k_{1},\mathbf{k}_{2}\right)\in\mathbb{N}\times\mathbb{N}^{r},\quad C_{\mathbf{k}}\in\mathbb{R}^{n}.\label{eq:1D v_SSM-1-1}
\end{align}
In this expression, we have omitted the modulus of $u_{1}$ in the
fractional powered terms as we are only interested in the positive
branch of $\mathcal{W}^{*}\left(E_{1}\right)$ along which $u_{1}\geq0$
can be chosen. A linear term with $\mathbf{k}=\left(1,\mathbf{0}\right)$
also appears in the expression (\ref{eq:1D v_SSM-1-1}) because the
linear part of system \eqref{eq:perturbationNS} has not been diagonalized
in the current set of coordinates.

The spectral ratio $\log\kappa_{1}/\log\lambda_{1}=1.989703$ in Table
\ref{tab:eigenvalues} shows that $\mathcal{W}^{*}\left(E_{1}\right)$
is only once continuously differentiable by formulas \eqref{eq:eta formula for smoothness - discrete}
if the coefficient $h_{\mathbf{k}}$ with \textbf{$\mathbf{k}=\left(0,1,0,\ldots,0\right)$
}is nonzero.
\begin{table}[H]
\centering{}%
\begin{tabular}{|c|c|c|}
\hline 
 & value & $\frac{\log\kappa_{j}}{\log\left|\lambda_{1}\right|}$ (spectral ratio)\tabularnewline
\hline 
\hline 
$\log\left|\lambda_{1}\right|$ & -0.035068 & \tabularnewline
\hline 
$\log\kappa_{1}$ & -0.069776 & 1.989703\tabularnewline
\hline 
$\log\kappa_{2}$ & -0.073369 & 2.092178\tabularnewline
\hline 
$\log\kappa_{3}$ & -0.140274 & 4.000013\tabularnewline
\hline 
$\log\kappa_{4}$ & -0.168877 & 4.815674\tabularnewline
\hline 
\end{tabular}\caption{Leading eigenvalues of the base state in the plane Couette flow at
$Re=134.52$ . The second column contains the first 5 eigenvalues,
while in the third column we report the corresponding spectral ratios
$\log\kappa_{j}/\log\left|\lambda_{1}\right|$.\label{tab:eigenvalues}}
\end{table}
 Furthermore, \citet{kaszas22} show that this heteroclinic connection
between the base state is a graph over the square root of the energy
input rate, $J=\sqrt{|I|},$ with the energy input rate 
\begin{equation}
I[\mathbf{u}]=\frac{1}{2L_{x}L_{z}}\int_{0}^{L_{x}}\int_{0}^{L_{z}}\left(\left(\frac{\partial u}{\partial y}u\right)_{y=1}+\left(\frac{\partial u}{\partial y}u\right)_{y=-1}\right)\text{d}z\ \text{d}x-1.\label{eq:I}
\end{equation}
Based on this, we will construct the fractional SSM $\mathcal{W}^{*}\left(E_{1}\right)$
as a graph over $J.$ Formally, this can be achieved by a change of
coordinates. The observable $J$ is a function of the state, 
\[
J=F(x)=F(u_{1},...,u_{p},a,b,v,c,d),
\]
where we recall the decomposition of \eqref{eq:linear change of coordinates-1}.
By the inverse function theorem, on a neighborhood of $x=0$, there
exists a function $G$, such that $u_{1}=G(J,u_{2},...,u_{p},a,b,v,c,d)$,
as long as 
\[
\left.\frac{\partial F}{\partial u_{1}}\right|_{x=0}\neq0.
\]
 Let us now introduce the new state vector $\hat{x}=(J,\hat{u}_{2},...,\hat{u}_{p},\hat{a},\hat{b},\hat{v},\hat{c},\hat{d})$,
which is related to $x$ through $x=g(\hat{x}),$ where 
\[
g(\hat{x})=\left(G(J,\hat{u}_{2},...,\hat{u}_{p},\hat{a},\hat{b},\hat{v},\hat{c},\hat{d}),\hat{u}_{2},\dots,\hat{u}_{p},\hat{a},\hat{b},\hat{v},\hat{c},\hat{d}\right)^{\mathrm{T}}.
\]

Under this change of coordinates, \eqref{eq:original system} must
be transformed as 
\begin{align}
\dot{\hat{x}} & =(D_{\hat{x}}g)^{-1}\mathcal{A}(D_{\hat{x}}g)\hat{x}+\hat{f}(\hat{x}).\label{eq:transformed_eq_withObs}
\end{align}

Since the linear part of \eqref{eq:transformed_eq_withObs} is related
to the linear part of \eqref{eq:original system} through a similarity
transformation, they have the same spectrum. This means, that we may
view the manifolds as graphs over the variable $J$, and the parametrization
contains the same fractional powers as the parametrization in the
original phase space. 

To compute $\mathcal{W}^{*}\left(E_{1}\right)$ up to order $5$ ,
we need to find all terms in \eqref{eq:1D v_SSM-1-1} with integer
and fractional powers up to $\mathcal{K}=5$. These powers satisfy
\[
k_{1}+\sum_{\ell=1}^{\mathcal{K}-k_{1}}k_{4\ell}\frac{\log\kappa_{\ell}}{\log\left|\lambda_{1}\right|}\leq5.
\]
Since $\frac{\log\kappa_{5}}{\log|\lambda_{1}|}>5,$ we only have
to use fractional terms arising from the eigenvalues with $\kappa_{1},...,\text{\ensuremath{\kappa_{4}} in }$
\eqref{eq:1D v_SSM-1-1}. We also see from Table \ref{tab:eigenvalues}
that the fractional powers $\frac{\log\kappa_{j}}{\log\left|\lambda_{1}\right|}$
and $\frac{\log\kappa_{j}}{\log\left|\lambda_{1}\right|}$ are very
close to the integer powers $2$ and $4$, which are also present
in the expansion \eqref{eq:1D v_SSM-1-1}. While some of the near-integer
spectral ratios arise due to resonances in the infinite-dimensional
system, they are never exact resonances in the finite truncation of
the system that we work with here. Nevertheless, to avoid sensitivity
and overfit (see Remark \ref{rem:explanation on overfit from nearly equal powers}),
we set $k_{41}=k_{43}=0$. Based on these considerations, a quintic-order
approximation for the SSM-reduced map $J(\imath)\mapsto J(\imath+1)$
can be written as 
\begin{equation}
J(\imath+1)=\sum_{k_{1}+k_{42}\frac{\log\kappa_{2}}{\log\left|\lambda_{1}\right|}+k_{44}\frac{\log\kappa_{4}}{\log\left|\lambda_{1}\right|}\leq5}R_{k_{1},k_{42},k_{45}}J(\imath)^{k_{1}+k_{42}\frac{\log\kappa_{2}}{\log|\lambda_{1}|}+k_{44}\frac{\log\kappa_{4}}{\log|\lambda_{1}|}}.\label{eq:reddynamics_Couette}
\end{equation}

We determine the coefficients $R_{k_{1},k_{42},k_{44}}$ from linear
regression by minimizing the squared error 

\[
\sum_{\imath=1}^{N_{s}}\left\Vert J(\imath+1)-\sum_{k_{1}+k_{42}\frac{\log\kappa_{2}}{\log\left|\lambda_{1}\right|}+k_{44}\frac{\log\kappa_{4}}{\log\left|\lambda_{1}\right|}\leq5}R_{k_{1},k_{42},k_{45}}J(\imath)^{k_{1}+k_{42}\frac{\log\kappa_{2}}{\log|\lambda_{1}|}+k_{44}\frac{\log\kappa_{4}}{\log|\lambda_{1}|}}\right\Vert ^{2},
\]
along a training trajectory started on the heteroclinic orbit for
the available data snapshots indexed by $\imath=1,\ldots,N_{s}$.
This regression gives the coefficients listed in Table \ref{tab:coefficients}
for the reduced mapping (\ref{eq:reddynamics_Couette}): 

\begin{table}[H]
\centering{}%
\begin{tabular}{|c|c|c|c|}
\hline 
$k_{1}$ & $k_{42}$ & $k_{44}$ & $R_{k_{1},k_{42},k_{44}}$\tabularnewline
\hline 
\hline 
1 & 0 & 0 & \textrm{0.96182}\tabularnewline
\hline 
2 & 0 & 0 & \textrm{0.82616}\tabularnewline
\hline 
0 & 1 & 0 & \textrm{-1.03809}\tabularnewline
\hline 
3 & 0 & 0 & \textrm{0.43671}\tabularnewline
\hline 
1 & 1 & 0 & 0.35338\tabularnewline
\hline 
2 & 1 & 0 & -0.61552\tabularnewline
\hline 
4 & 0 & 0 & \textrm{-1.65568}\tabularnewline
\hline 
0 & 2 & 0 & \textrm{0.67765}\tabularnewline
\hline 
0 & 0 & 1 & \textrm{4.52449}\tabularnewline
\hline 
5 & 0 & 0 & \textrm{-3.47146}\tabularnewline
\hline 
\end{tabular}\caption{Coefficients $R_{k_{1},k_{42},k_{44}}$ in the reduced model (\ref{eq:reddynamics_Couette})
for transitions along the fractional SSM $W^{*}\left(E_{1}\right)$
in the planar Couette flow.\label{tab:coefficients}}
\end{table}

\begin{figure}[H]
\centering{}\includegraphics[width=1\textwidth]{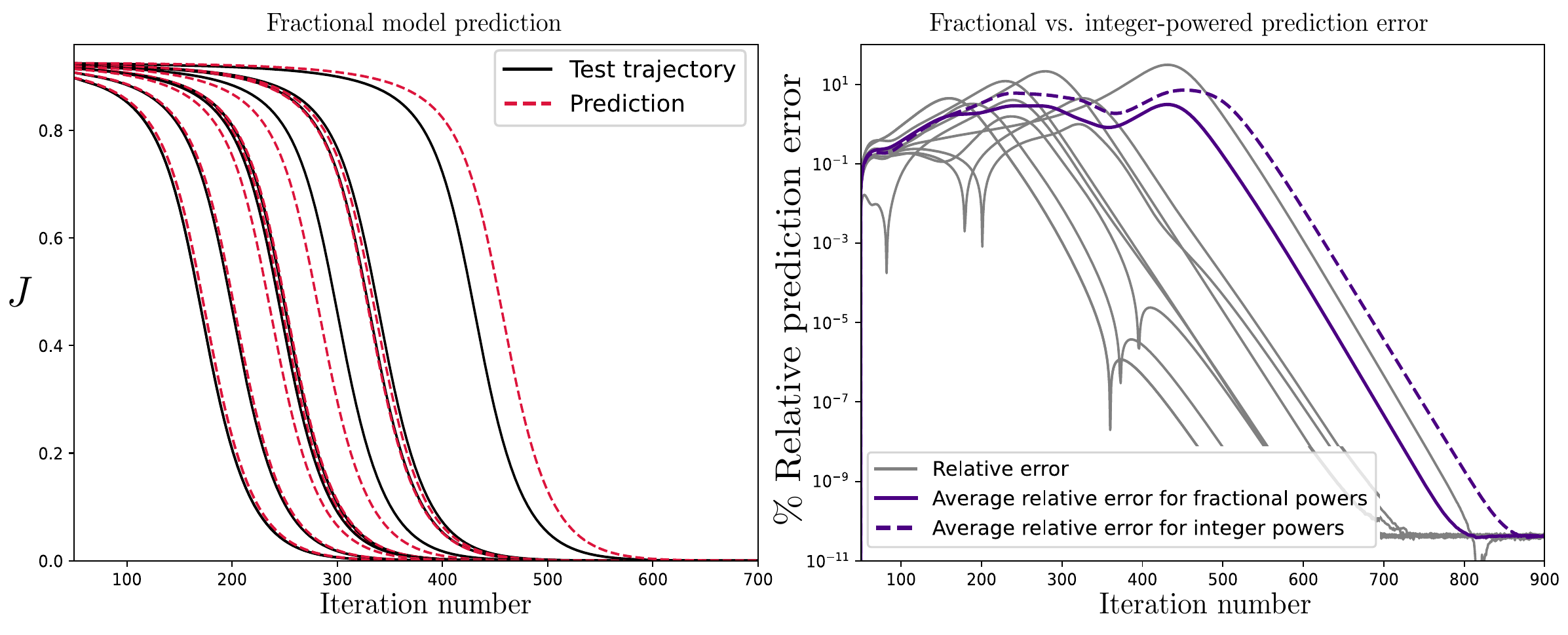}\caption{Performance of the fractional SSM-reduced model (\ref{eq:reddynamics_Couette})
for the planar Couette flow transition problem. Left: Reduced, 1D
fractional model predictions (red) for 10 trajectories of the full
system that were not used in regressing the model. Right: Instantaneous
(gray) and mean errors (solid purple) for the fractional model and
for a purely integer-powered fit (dashed purple) to the training trajectories.
The relative error is defined as $|J_{true}(\imath)-J_{predicted}(\imath)|/\max J_{true}(\imath)$.\label{fig: couettePredictions-1}}
\end{figure}

In Fig. \ref{fig: couettePredictions-1}, we compare the predictions
of the model \eqref{eq:reddynamics_Couette} with those of the model 

\begin{align}
J(\imath+1) & =0.96262J(\imath)+0.04265J^{2}(\imath)-0.10282J^{3}(\imath)+0.24485J^{4}(\imath)-0.14738J^{5}(\imath),\label{eq:CouetteIntegerModel}
\end{align}
which we obtain from a classical, integer-powered polynomial regression
up to order $\mathcal{K}=5$. Note that this model is purely a numerical
fit, given that the only SSM with an integer-powered expansion in
the problem is $\mathcal{W}\left(E_{1}\right)=E_{1}$, whose expansion
coefficients are all zero up to any order. 

As seen from Fig. \ref{fig: couettePredictions-1}, there is a noticeable
improvement in the correct, fractional expansion of the model relative
to the formal, polynomial fit to the same data set up to the same
order. The polynomial fit ignores five terms in the regression that
do appear in the actual reduced dynamics below quintic order. The
accuracy of the polynomial fit can only be increased by adding higher-order,
integer-powered polynomials, which in turn add higher derivatives
and hence more sensitivity to the model.

\section{Example 2: Data-driven, mixed-mode, SSM-reduced model for the dynamic
buckling of a beam\label{sec:Buckling beam example}}

In our second example, we discuss the data-driven construction of
a mixed-mode, SSM-reduced model for buckling using a discretized nonlinear
von Kármán beam model. This model is a second-order approximation
of geometrically exact beam theory that reproduces the bending-stretching
coupling arising from bending displacements in the order of the beam
thickness. We use the finite-element code of \citet{jain2018}, who
applied a combination of slow-fast reduction and SSM reduction to
study oscillations of this beam model. In contrast, here we apply
a sufficiently large longitudinal compression force and seek to derive
an SSM-reduced nonlinear model for the ensuing buckling dynamics.

\begin{figure}[ht]
\centering{}\includegraphics[width=1\textwidth]{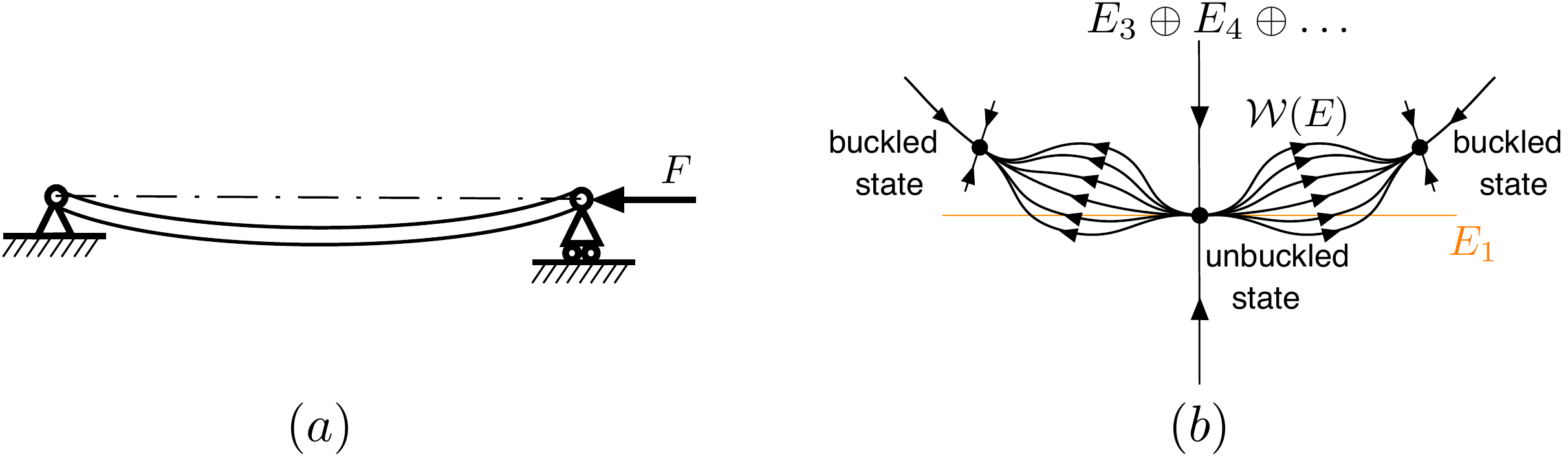}\caption{(a) Schematic of the geometry of the beam buckling under a horizontal
force $F$. (b) Schematic of the family $\mathcal{W}\left(E\right)$
of 2D mixed-mode SSMs connecting the unstable unbuckled state to the
two stable buckled states. The whole family $\mathcal{W}\left(E\right)$
is tangent to the unstable subspace $E_{1}$, as shown. However, only
$\mathcal{W}^{\infty}(E)$ is tangent to the full spectral subspace
$E$. \label{fig:vk_beam}}
\end{figure}

Schematically shown in Fig. \ref{fig:vk_beam}, the finite element
model of the von Kármán beam consists of $13$ nodes, each with three
degrees of freedom (DOF): axial, transversal and rotational. By setting
boundary conditions as pinned-pinned, we constrain 2 DOFs of the leftmost
node and 1 DOF of the rightmost node. Therefore, the full system before
model reduction has 12 elements and 72 DOFs. The beam is made of steel,
with a Young's modulus $E=1.9\cdot10^{11}\,\mathrm{Pa}$, density
$\rho=7,850\,kg/m^{3}$, viscous damping rate $\kappa=7\cdot10^{6}\,kg\,m/s$,
length $L=2\,m$, height $h=0,01\,m$ and width $w=0.05\,m$. Euler's
critical horizontal load for buckling is given by $P_{N}=N^{2}\pi^{2}EI/L^{2}$,
with $I$ denoting the area moment of inertia of the cross section
of the beam. Accordingly, to generate the first buckling mode ($N=1$),
we apply the force
\begin{equation}
F=N^{2}\pi^{2}EI/L^{2}=\frac{1}{12}\pi^{2}Ebh^{3}/L^{2}.\label{eq:applied buckling force}
\end{equation}
. 

\subsection{2D primary mixed-mode SSM\label{subsec:2D-primary-mixed-mode}}

Under the external force (\ref{eq:applied buckling force}), the beam
has three equilibria: one unstable equilibrium corresponding to purely
axial displacement and two stable equilibria corresponding to the
upper and lower buckled states. The first $10$ eigenvalues of the
linearized finite element model at the unstable equilibrium are listed
in Table \ref{tab: first ten eigenvalues for buckling beam}, where
we have already labelled the eigenvalues based on the role they will
play in the construction of a mixed-mode SSM tangent to the 2D eigenspace
\begin{equation}
E=E_{1}\oplus E_{2},\label{eq:2D E for beam}
\end{equation}
with the 1D eigenspaces, $E_{1}$ and $E_{2}$ corresponding to the
eigenvalues $\lambda_{1}=11.06$ and $\lambda_{2}=-11.10$, respectively.
The remaining eigenvalues are all of the form $\beta_{m}\pm i\nu_{m}$
with $m>4$. 
\begin{table}[H]
\begin{centering}
\begin{tabular}{|c|c|c|c|c|c|}
\hline 
$\lambda_{1}$ & $\lambda_{2}$  & $\beta_{1}\pm i\nu_{1}$  & $\beta_{2}\pm i\nu_{2}$  & $\beta_{3}\pm i\nu_{3}$  & $\beta_{4}\pm i\nu_{4}$ \tabularnewline
\hline 
\hline 
$11.06$ & $-11.10$ & $-0.36\pm i119.36$ & $-1.83\pm i295.56$ & $-5.80\pm i541.50$ & $-14.19\pm i858.19$\tabularnewline
\hline 
\end{tabular}
\par\end{centering}
\caption{The first $10$ eigenvalues of the linearized finite-element model
for a buckling beam at its unstable fixed point. Our notation is relative
to the 2D eigenspace $E$ defined in (\ref{eq:2D E for beam}). \label{tab: first ten eigenvalues for buckling beam}}
 
\end{table}

To obtain a reduced model that captures the dynamics of buckling,
we seek reduction to a member of the SSM family $\mathcal{W}\left(E\right)$.
Therefore, we have $p=2$, $q=r=0$ and $s=n-2$ with $n=144$ in
terms of the notation used in Theorem \ref{thm:main1}. By statement
(ii) Theorem \ref{thm:main1}, a typical member of this mixed-mode
family is only continuous, given that
\[
\eta=\mathrm{Int}\left[\min_{j,k,\ell,m}^{+}\left\{ \frac{\beta_{m}}{\lambda_{j}}\right\} \right]=\mathrm{Int}\left[\min_{m\geq1}\left\{ \frac{\beta_{m}}{\lambda_{2}}\right\} \right]=\mathrm{Int}\left[\frac{0.36}{11}\right]=0.
\]
All these fractional SSMs and the primary SSM, $\mathcal{W}^{\infty}(E)$,
contain all three fixed points and are all tangent to the 2D stable
subspaces of the two buckled fixed points, as illustrated in Fig.
\ref{fig:vk_beam}b. For this reason, we select $\mathcal{W}^{\infty}(E)$
for constructing a reduced-order model because we can then use a classical
polynomial expansion without fractional powers.

\begin{figure}[h]
\centering{}\includegraphics[width=1\textwidth]{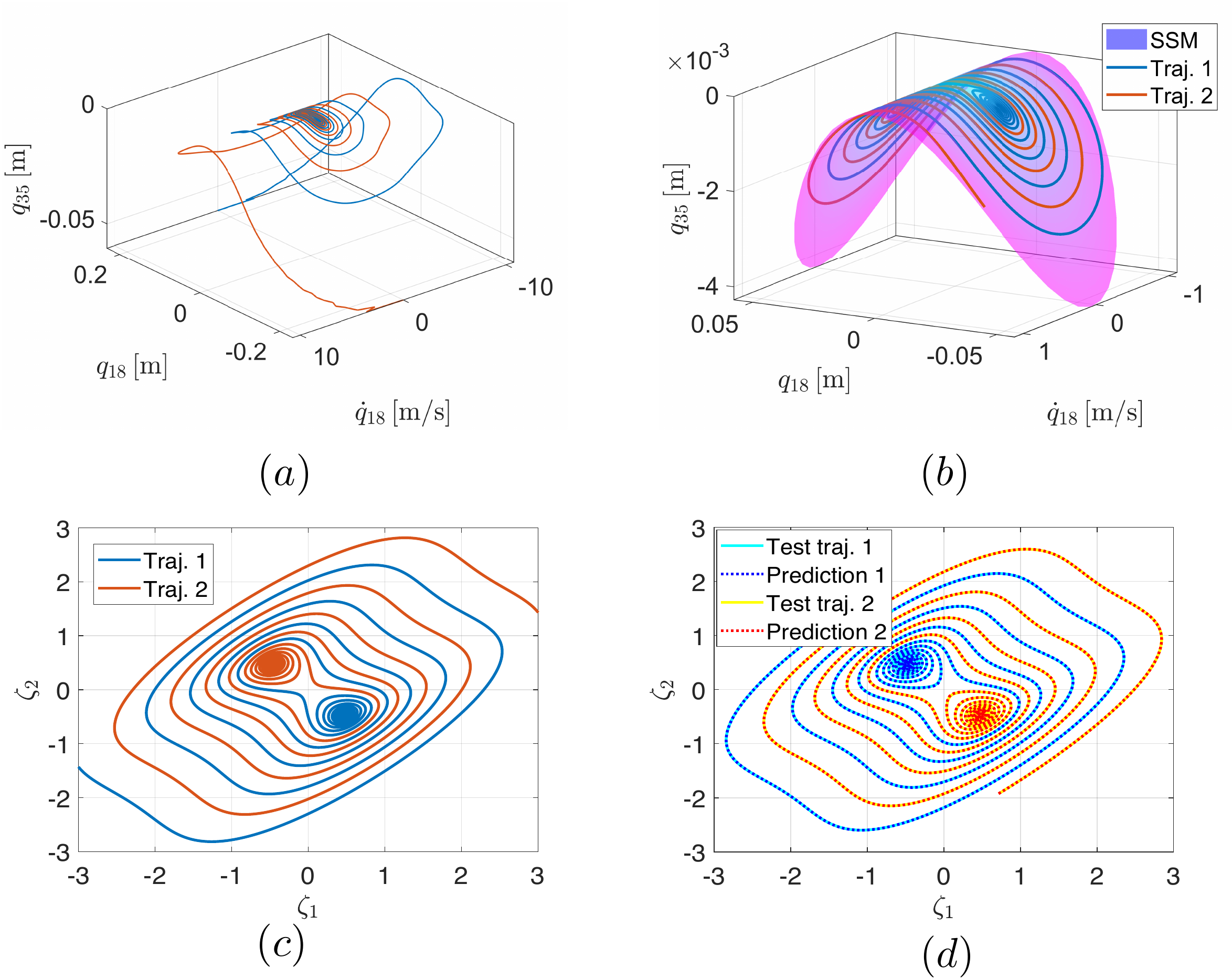}\caption{Construction and testing of a reduced model on the 2D, primary, mixed-mode
SSM of the buckling beam problem. (a) Training trajectories. (b) The
primary mixed-mode SSM, $\mathcal{W}^{\infty}(E)$, together with
the two training trajectories. (c) Training trajectories in model
coordinates. (d) Test trajectories of the full system along with model
predictions for them. \label{fig:vk_beam 2D SSM}}
\end{figure}

By Remark \ref{rem: explanation on obtaining primary SSMs}, the \emph{SSMLearn}
package of \citet{cenedese21} can be used directly to approximate
the primary mixed-mode SSM in this example. We use two training trajectories
whose initial conditions are obtained via a positive and a negative
vertical force applied to the midpoint of the beam, respectively.
We then turn off the vertical forces and record the approach of the
resulting two trajectories to the upper and lower buckled state of
the beam. The two trajectories obtained in this fashion are shown
in Fig. \ref{fig:vk_beam 2D SSM}a, projected to a 3D coordinate space
comprising the transverse displacement $q_{18}$ and velocity $\dot{q}_{18}$
of the midpoint, as well as the axial displacement $q_{35}$ of right
endpoint. Figure \ref{fig:vk_beam 2D SSM}b shows the primary mixed-mode
SSM, $\mathcal{W}^{\infty}(E),$ projected onto the same coordinate
system, along with truncated versions of the two training trajectories.
Figure \ref{fig:vk_beam 2D SSM}c shows truncated training trajectories
in model coordinates from the spectral subspace $E=E_{1}\oplus E_{2}$.
Finally, Figure \ref{fig:vk_beam 2D SSM}d shows two truncated test
trajectories of the full mechanical model system and predictions for
them by the $\mathcal{W}^{\infty}(E)$-reduced and normalized 2D model.

Following the procedure outlined by \citet{cenedese21} for data-driven
primary SSMs, we find that both the invariance error of the manifold
and the accuracy of the reduced-order model for the dynamics on the
manifold is minimal at a polynomial order of $11$. This minimum,
however, is quite weak for the invariance error and allows us to choose
an order $\mathcal{K}=7$ approximation that reaches almost the same
accuracy, resulting in a $0.1\%$ normalized error in fitting $\mathcal{W}^{\infty}(E)$
to the trajectories (see \citet{cenedese21} on computing this error).
Figure \ref{fig:vk_beam 2D SSM}b shows the result of the data-driven
identification of $\mathcal{W}^{\infty}(E)$ by SSMLearn, along with
two training trajectories used in the procedure. Each training trajectory
is truncated to a segment with no initial transients, which ensures
that the trajectory is close to the SSM to be inferred from the data. 

The training trajectories are also shown in Fig. \ref{fig:vk_beam 2D SSM}c
in modal coordinates that parametrize the spectral subspace $E.$
Finally, Fig. \ref{fig:vk_beam 2D SSM}d shows a highly accurate prediction
by the 2D $7^{th}$-order reduced model on $\mathcal{W}^{\infty}(E)$
for two test trajectories of the full system that were not used in
the construction of the SSM-reduced model. We applied no normal form
transformation to this model, because the two stable buckled states
are far from the origin and hence are fully expected to lie outside
the validity of a truncated normal form.

\subsection{4D primary mixed-mode SSM\label{subsec:4D-primary-mixed-mode}}

The reduced-order model derived in the previous section already captures
the core dynamics (ie., all stable and unstable steady states) of
the bucking beam.To illustrate how enlarging the dimension of the
SSM reveals more transient behavior for larger initial conditions
(see Remark \ref{rem:choice of dimension for SSM}), we repeat the
above procedure for a 4D primary SSM constructed over the enlarged
spectral subspace
\begin{equation}
E=E_{1}\oplus E_{2}\oplus E_{3}.\label{eq:4D E}
\end{equation}
Here the newly added, 2D real eigenspace $E_{3}$ corresponds to the
eigenvalue pair $-0.36\pm i119.36$ in Fig. \ref{tab: first ten eigenvalues for buckling beam}.
This means that we now have $p=2$, $q=1$, $r=0$ and $s=n-4$ with
$n=144$ in the application of Theorem \ref{thm:main1}. Accordingly,
our notation for the eigenvalues featured in the theorem changes from
that in Fig. \ref{tab: first ten eigenvalues for buckling beam} to
that in Table \ref{tab: first ten eigenvalues for buckling beam-1}.
\begin{table}[H]
\begin{centering}
\begin{tabular}{|c|c|c|c|c|c|}
\hline 
$\lambda_{1}$ & $\lambda_{2}$ & $\alpha_{1}\pm i\omega_{1}$ & $\beta_{1}\pm i\nu_{1}$ & $\beta_{2}\pm i\nu_{2}$ & $\beta_{3}\pm i\nu_{3}$\tabularnewline
\hline 
\hline 
$11.06$ & $-11.10$ & $-0.36\pm i119.36$ & $-1.83\pm i295.56$ & $-5.80\pm i541.50$ & $-14.19\pm i858.19$\tabularnewline
\hline 
\end{tabular}
\par\end{centering}
\caption{Same eigenvalues as in Fig. \ref{tab: first ten eigenvalues for buckling beam},
but now labelled with respect to the modified definition (\ref{eq:4D E})
of the spectral subspace $E$. \label{tab: first ten eigenvalues for buckling beam-1}}
\end{table}

By Theorem \ref{thm:main1}, a typical fractional SSM in the $\mathcal{W}\left(E\right)$
family is still only continuous, given that \emph{
\begin{equation}
\eta=\mathrm{Int}\left[\min_{j,k,\ell,m}^{+}\left\{ \frac{\beta_{m}}{\lambda_{j}},\quad\frac{\beta_{m}}{\alpha_{k}}\right\} \right]=\mathrm{Int}\left[\min_{m\geq1}\left\{ \frac{\beta_{m}}{\lambda_{2}},\frac{\beta_{m}}{\alpha_{1}}\right\} \right]=\mathrm{Int}\left[\frac{1.83}{11}\right]=0.\label{eq:typical smoothness measure-1}
\end{equation}
}The geometry depicted in Fig. \ref{fig:vk_beam}b remains similar
for the choice of $E$ in (\ref{eq:4D E}) which again prompts us
to choose the primary, mixed-mode SSM $\mathcal{W}^{\infty}\left(E\right)$
for reduced-order modeling. This SSM, however, is now 4D and hence
cannot be easily visualized. 

\begin{figure}[ht]
\centering{}\includegraphics[width=1\textwidth]{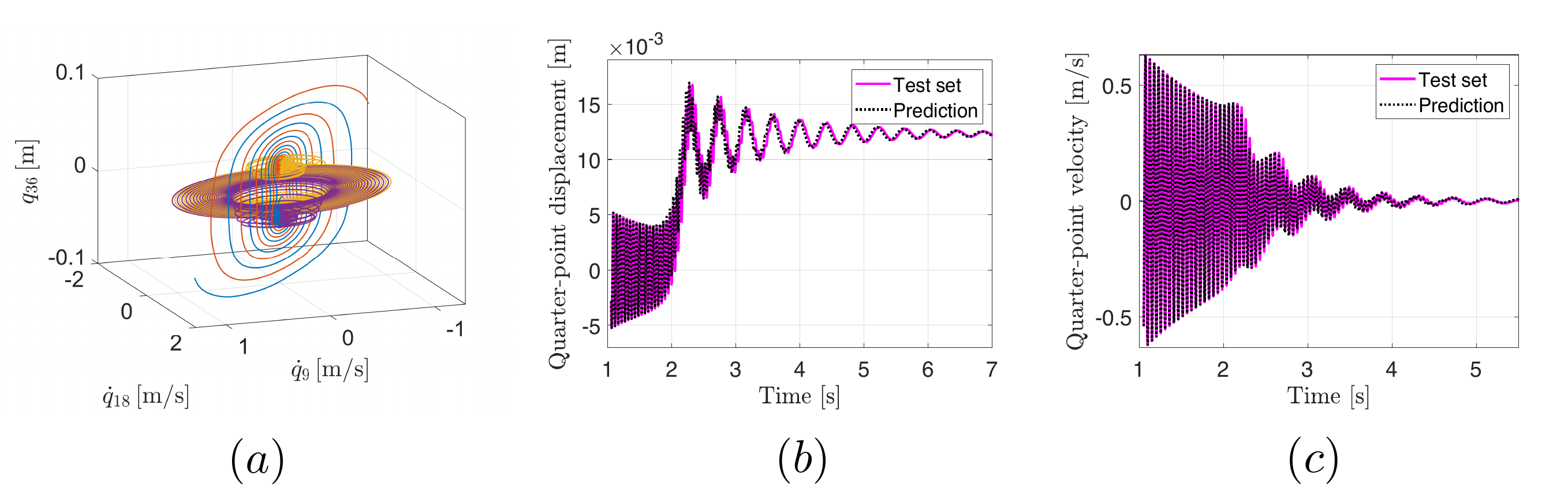}\caption{Construction and testing of a reduced model on the 4D, primary mixed-mode
SSM of the buckling beam problem. (a) Four training trajectories projected
to a 3D plane of the $\dot{q}_{18}$ transverse velocity of midpoint,
the $\dot{q}_{9}$ transverse velocity of quarter-node and the $q_{36}$
rotation of the endpoint. (b) The quarter-node displacement of a test
trajectory of the full beam model and a prediction for it by the $\mathcal{W}^{\infty}(E)$-reduced
and normalized 4D model. (c) Same as (b) but for the quarter-node
velocity of a test trajectory. \label{fig:vk_beam 4D SSM}}
\end{figure}

We initialize four new training trajectories in a way that they involve
oscillations from the second buckling mode as well. To generate each
training trajectory, we apply two vertical forces of opposing orientations
at $1/4$ and $3/4$ of the total length of the beam to create a static,
horizontal S-shaped equilibrium state for the beam. We then turn off
these forces and use the resulting trajectories shown in Fig. \ref{fig:vk_beam 4D SSM}a
for extracting the 4D mixed-mode SSM $\mathcal{W}^{\infty}(E)$ and
its reduced dynamics. Figures \ref{fig:vk_beam 4D SSM}b-c show extended
model prediction accuracy that now even capture transients of the
test trajectories that were not used in constructing the model.

We close by recalling that within mixed-mode SSMs families, just as
within like-mode SSMs, there are primary and secondary (fractional)
SSMs. For instance, the SSMs computed in Sections \ref{subsec:2D-primary-mixed-mode}
and \ref{subsec:4D-primary-mixed-mode} are primary mixed-mode SSMs.
They also have fractional mixed-mode counterparts, which we did not
have to compute in this case because of the manifold geometry shown
in Fig. \ref{fig:vk_beam}b.

\section{Example 3: Fractional and mixed-mode SSMs in a forced oscillator
system\label{sec:Shaw-Pierre example}}

\begin{figure}[ht]
\centering{}\includegraphics[width=0.6\textwidth]{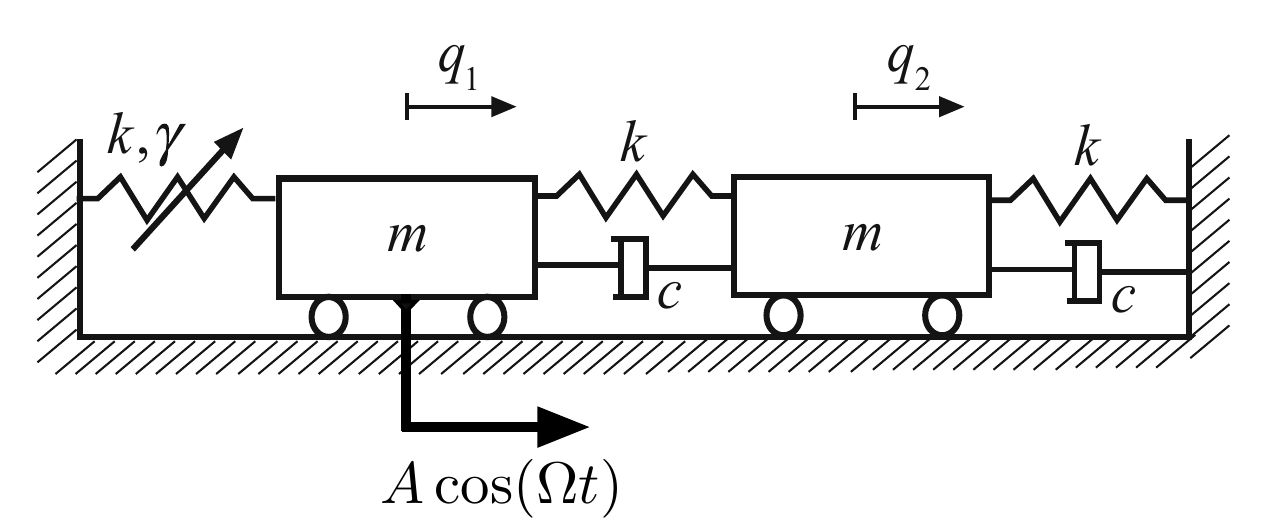}\caption{The periodically forced Shaw\textendash Pierre oscillator from \citet{ponsioen2019}.
\label{fig:Shaw-Pierre sketch}}
\end{figure}

Our third example involves a forced version of the 2 DOF nonlinear
oscillator system originally proposed by \citet{shaw93}, which has
been widely used to illustrate the concept of nonlinear normal modes
in mechanical systems. A forced and slightly modified version of this
example shown in Fig. \ref{fig:Shaw-Pierre sketch} was considered
by \citet{ponsioen2019} with small forcing amplitudes for SSM-based
model reduction. That system comprises two identical rigid bodies
of mass $m>0$ that move horizontally without friction, with their
positions marked by the coordinates $q_{1}$ and $q_{2}$. The bodies
are connected to each other and two walls via identical springs of
linear stiffness $k>0$. The first spring also has a softening cubic
nonlinearity with coefficient $\gamma>0$, while the other two springs
are subject to linear damping with coefficient $c>0$. The first body
is subjected to external forcing of amplitude $A$ and frequency $\Omega$.

With the notation $p_{1}=\dot{q}_{1}$, $p_{2}=\dot{q}_{2}$ , the
first-order form of the equations of motion for the system is
\begin{equation}
\left(\begin{array}{c}
\dot{q_{1}}\\
\dot{p}_{1}\\
\dot{q_{2}}\\
\dot{p_{2}}
\end{array}\right)=\left(\begin{array}{cccc}
0 & 1 & 0 & 0\\
-2\frac{k}{m} & -\frac{c}{m} & \frac{k}{m} & \frac{c}{m}\\
0 & 0 & 0 & 1\\
\frac{k}{m} & \frac{c}{m} & -2\frac{k}{m} & -2\frac{c}{m}
\end{array}\right)\left(\begin{array}{c}
q_{1}\\
p_{1}\\
q_{2}\\
p_{2}
\end{array}\right)+\left(\begin{array}{c}
0\\
A\cos\Omega t-\frac{\gamma}{m}q_{1}^{3}\\
0\\
0
\end{array}\right).\label{eq:shawpierreODE}
\end{equation}
We first identify fractional SSMs near the trivial equilibrium of
this system in the absence of external forcing $(A=0$). As a next
step, under the addition of external time-periodic forcing ($A>0$),
we construct a data-driven reduced model on a primary mixed-mode SSM
of an unstable periodic orbit to capture transitions between coexisting
periodic orbits. This second part of the example, therefore, serves
as an illustration of our mixed-mode SSM results for discrete dynamical
systems.

\subsection{Fractional SSMs in the unforced Shaw\textendash Pierre example ($A=0$) }

Following \citet{shaw93} and \citet{ponsioen2018}, we fix the non-dimensionalized
parameter values $m=1,$ $c=0.3,k=1,\gamma=0.5$. In that case, the
eigenvalues of the linearized system near the origin are 
\begin{equation}
\alpha\pm i\omega=-0.0741\pm1.0027i\quad\beta\pm i\nu=-0.3759\pm1.6812i,\label{eq:eigenvalues for Shaw-Pierre examole}
\end{equation}
where we have already applied the notation of Theorem \ref{thm:main1}
to the 2D spectral subspace $E$ corresponding to the slower decaying
pair of complex eigenvalues. Then, after a linear change of coordinates,
system (\ref{eq:shawpierreODE}) with $A=0$ will take the form
\begin{equation}
\left(\begin{array}{c}
\dot{a}\\
\dot{b}\\
\dot{c}\\
\dot{d}
\end{array}\right)=\left(\begin{array}{cccc}
\alpha & \omega & 0 & 0\\
-\omega & \alpha & 0 & 0\\
0 & 0 & \beta & \nu\\
0 & 0 & -\nu & \beta
\end{array}\right)\left(\begin{array}{c}
a\\
b\\
c\\
d
\end{array}\right)+f(a,b,c,d).\label{eq:shawpierreODEtransformed}
\end{equation}
This system is of the general form (\ref{eq:original system}) with
$j=0$, $k=1$, $\ell=0$ and $m=1$. Therefore, formula (\ref{eq:complex SSM 1})
simplifies to the single complex function 
\begin{align*}
\mathcal{\mathcal{E}}(z) & =Q\left|z\right|^{\frac{\beta}{\alpha}}e^{i\frac{\nu}{\alpha}\log\left|z\right|},
\end{align*}
where $Q\in\mathbb{C}$ is a complex parameter. 

In terms of the real $(a,b,c,d)$ coordinates, the full family of
2D invariant graphs of the linear part over the $E$ spectral subspace
take the form 
\begin{equation}
\left(\begin{array}{c}
c(a,b)\\
d(a,b)
\end{array}\right)=P_{1}\left(a^{2}+b^{2}\right)^{\frac{\beta}{2\alpha}}\left(\begin{array}{c}
\cos\left(P_{2}+\frac{\nu}{2\alpha}\log\left[a^{2}+b^{2}\right]\right)\\
\sin\left(P_{2}+\frac{\nu}{2\alpha}\log\left[a^{2}+b^{2}\right]\right)
\end{array}\right),\quad P_{1}=\left|Q\right|,\quad P_{2}=\arg Q,\label{eq:shaw-pierre linear mfds}
\end{equation}
which defines a two-parameter-family of invariant graphs in the linearized
system, parametrized by $P_{1},P_{2}\in\mathbb{R}$. 

As guaranteed by statement (i) of Theorem \ref{thm:main1}, the corresponding
2D, slow SSM family $\mathcal{W}\left(E\right)$ in the full nonlinear
system (\ref{eq:shawpierreODEtransformed}) can be written near the
origin as\emph{
\begin{align}
w & =\mathcal{\mathcal{E}}(z)+\sum_{2\leq k_{1}+k_{2}+\frac{\left(k_{3}+k_{4}\right)\beta}{2\alpha}\leq\mathcal{K}}h_{\mathbf{k}}z^{k_{2}}\bar{z}^{k_{3}}\mathcal{\mathcal{E}}^{k_{5}}\bar{\mathcal{\mathcal{E}}}^{k_{6}}+o\left(\left|z\right|^{\mathcal{K}}\right),\qquad\mathbf{k}=\left(k_{2},k_{3},k_{5},k_{6}\right).\label{eq:SSM in original coordinates real-1-5-3-1}
\end{align}
}By the complexification formula (\ref{eq:complexification}), when
truncated at order $\mathcal{K}$ and re-indexed, this SSM family
can be written in real coordinates as
\begin{align}
\left(\begin{array}{c}
c\\
d
\end{array}\right) & =\sum_{2\leq k_{1}+k_{2}+\frac{\left(k_{3}+k_{4}\right)\beta}{2\alpha}\leq\mathcal{K}}P_{1}^{k_{3}+k_{4}}\left(\begin{array}{c}
C_{\mathbf{k}}\\
D_{\mathbf{k}}
\end{array}\right)a^{k_{1}}b^{k_{2}}\left(a^{2}+b^{2}\right)^{\frac{\beta\left(k_{3}+k_{4}\right)}{2\alpha}}\times\nonumber \\
 & \,\,\,\,\,\,\,\,\,\,\,\,\,\,\,\,\times\cos^{k_{3}}\left(P_{2}+\frac{\nu}{2\alpha}\log\left[a^{2}+b^{2}\right]\right)\sin^{k_{4}}\left(P_{2}+\frac{\nu}{2\alpha}\log\left[a^{2}+b^{2}\right]\right),\label{eq:(c,d) expansion for Shaw-Pierre}
\end{align}
 where 
\[
C_{0010}=D_{0001}=1,
\]
and all other $C_{\mathbf{k}}$ and $D_{\mathbf{k}}$ coefficients
with $\left|\mathbf{k}\right|=1$ are zero. 

In all nonlinear normal mode and SSM calculations available in the
literature for this model, only the $k_{3}=k_{4}=0$ terms of the
above expansion have been used and hence only the unique primary SSM,
$W^{\infty}\left(E\right)$, has been exploited for reduced-order
modeling. This restriction to $W^{\infty}\left(E\right)$ is fully
justified when $\beta/\alpha$ is relatively large ( $\beta/\alpha>\mathcal{K}$),
in which case the fractional terms do not yet appear in the truncated
expansion (\ref{eq:(c,d) expansion for Shaw-Pierre}). For the parameter
configuration used by \citet{shaw93} and \citet{ponsioen2018}, the
eigenvalues in (\ref{eq:eigenvalues for Shaw-Pierre examole}) give
\[
\frac{\beta}{\alpha}=5.0729,
\]
which means that fractional powers only appear in the expansion for
any member of the $\mathcal{W}\left(E\right)$ family for orders $\mathcal{K}>5$.
The $6^{th}$ order Taylor coefficients happen to vanish in this case,
but the fractional coefficient $\sqrt{a^{2}+b^{2}}^{5.0729\left(k_{3}+k_{4}\right)}$
does appear and starts shaping the secondary SSMs and their reduced
dynamics already below order $6$. 

\begin{figure}[h]
\centering{}\includegraphics[width=0.9\textwidth]{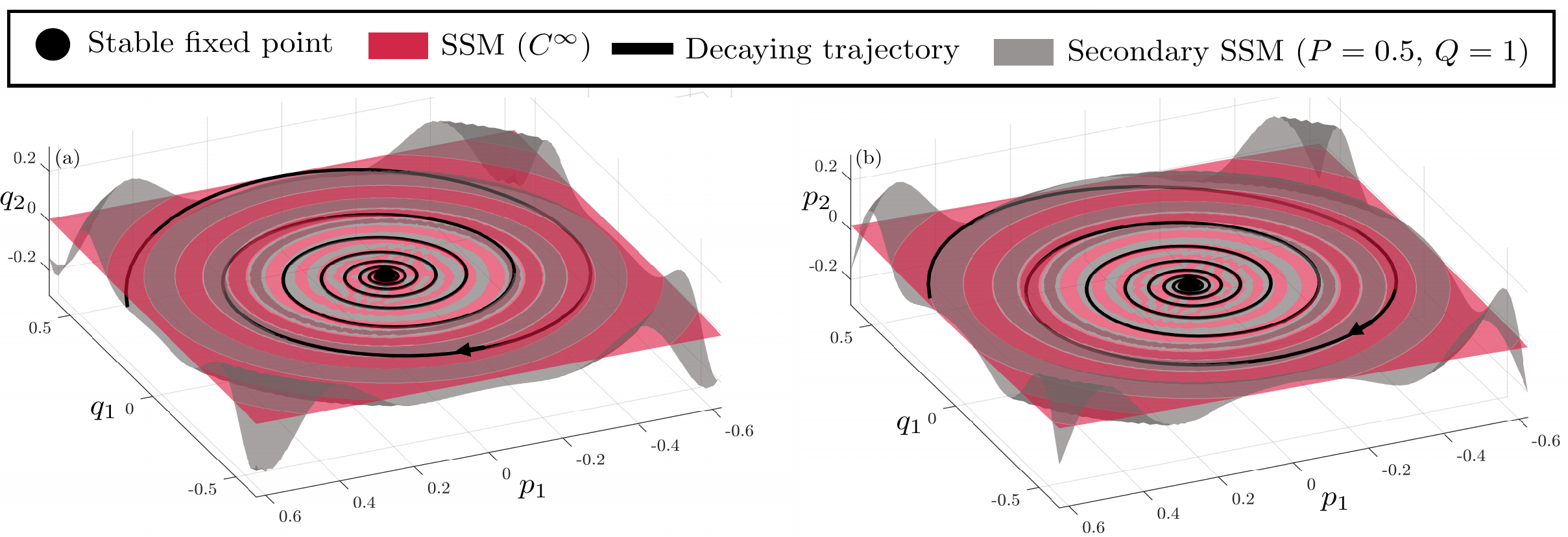}\caption{Primary and fractional slow SSMs in the unforced ($A=0$) limit of
the Shaw-Pierre model (\ref{eq:shawpierreODE}). The $q_{1}$ and
$p_{1}$ coordinates of the slow SSMs are plotted in panels (a) and
(b), respectively. Red surface: order $K=7$ approximation of the
analytic SSM $W^{\infty}\left(E\right)$. Grey surface: fractional
SSM of order $\mathcal{K}=7$ with $P=0.5$ and $Q=1$. Black curve:
a general, decaying trajectory off the primary SSM. \label{fig:shawpierre}}
\end{figure}

To show the difference between the primary SSM $W^{\infty}\left(E\right)$
and the secondary (fractional) SSMs in the full SSM family (\ref{eq:(c,d) expansion for Shaw-Pierre}),
we explicitly calculate the Taylor expansion of the $C^{\infty}$
linearizing transformation underlying the results of Theorem \ref{thm:main1}
for the present model up to and including order 7. In Fig. \ref{fig:shawpierre}
we show the transformed images of some of the SSMs of the linearized
system \eqref{eq:shaw-pierre linear mfds} under the inverse of that
linearizing transformation. Note the difference between the smooth
$W^{\infty}\left(E\right)$ ($P=0$, $Q=0$ in \eqref{eq:shaw-pierre linear mfds})
and the general fractional SSMs, which arises because fractional terms
do start appearing beyond order $\mathcal{K}=6$. The figure also
shows a general, decaying trajectory that follows the fractional SSM
and only approaches $W^{\infty}\left(E\right)$ closer to the origin.

\subsection{Mixed-mode SSM of a periodic orbit in the forced Shaw\textendash Pierre
example ($A>0$) }

\begin{figure}[h]
\centering{}\includegraphics[width=0.6\textwidth]{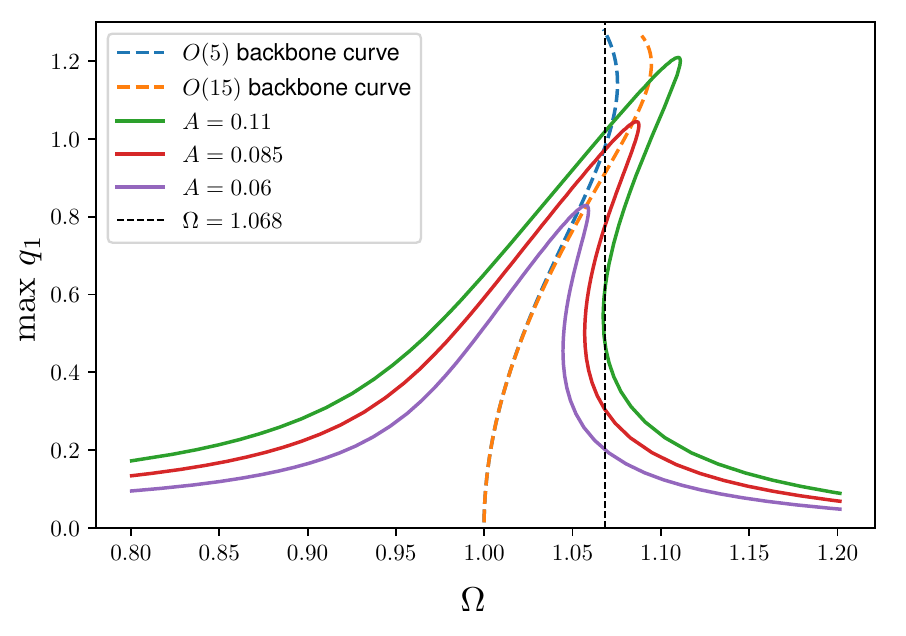}\caption{Forced response curves (i.e., maximal amplitudes of the observed steady-state
oscillations of the first mass as a function of the forcing frequency
$\Omega$) of the Shaw-Pierre system for different forcing amplitudes
$A$. Dashed lines represent the backbone curves (i.e., theoretical
centerpieces of the forced response curves) predicted by SSMTool from
a reduced-order model on the analytic slow SSM.\label{fig:shawpierre_FRC1}}
\end{figure}

Let us now consider system \eqref{eq:shawpierreODE} with a lower
value of damping, $c=0.03$, as in \citep{ponsioen2019} for which
the forced response curves of the model are given in Fig. \ref{fig:shawpierre_FRC1}
for several forcing amplitudes $A$. Also shown are the theoretically
computed backbone curves from SSMTool (see \citet{jain2022}), which
are based on a leading-order approximation to the time-periodic perturbation
of $W^{\infty}\left(E\right)$. Note that for larger forcing amplitudes
with $A>0.085$, even the $\mathcal{O}\left(15\right)$ theoretical
predictions based on $W^{\infty}\left(E\right)$ divert from the actual
periodic response.

\begin{figure}[H]
\centering{}\includegraphics[width=0.9\textwidth]{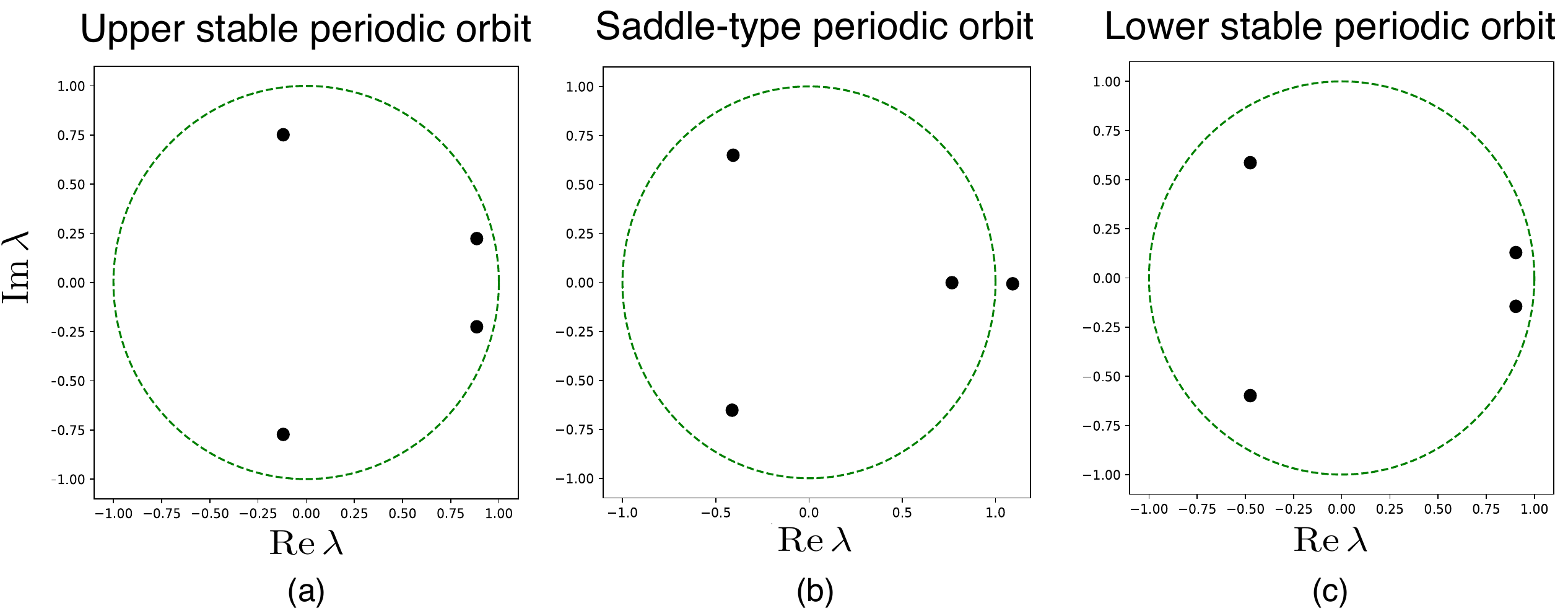}\caption{Panels (a)-(c) show, for $\Omega=1.07$ $A=0.11$, the Floquet multipliers
of the upper periodic orbit, the saddle-type periodic orbit in the
middle, and the lower periodic orbit, respectively. The green dashed
line represents the unit circle. \label{fig:shawpierre_floquet}}
\end{figure}

To improve on these predictions, we focus on the forcing frequency
$\Omega=1.07$ and forcing amplitude $A=0.11$  for which 3 periodic
orbits coexist: a high-amplitude stable periodic orbit, a saddle-type
periodic orbit and a low-amplitude stable periodic orbit. Their stability
types can be more specifically determined by computing their Floquet
multipliers numerically, which we show in Fig. \ref{fig:shawpierre_floquet}.
For our analysis, the most important is the spectrum of the saddle-type
periodic orbit. With the notation used in Theorem \ref{thm:main2},
the Floquet-multipliers of that periodic orbits are $\lambda_{1},\lambda_{2}$
and $\beta\pm i\nu$, where 
\begin{align*}
\lambda_{1} & =1.0835,\quad\lambda_{2}=0.7726,\quad\beta=-0.4132,\quad\nu=0.6474.
\end{align*}

To explore transitions among the coexisting periodic orbits, we consider
the time-$T$ Poincaré map associated with the forced system (\ref{eq:shawpierreODE}).
All three periodic orbits have period $T$ and hence all are hyperbolic
fixed points of this map, as seen from their Floquet multipliers in
Fig. \ref{fig:shawpierre_floquet}. To capture connections among the
fixed points, we seek to construct the 2D mixed-mode SSM family tangent
to the 2D spectral subspace $E$ that is spanned by the stable and
unstable real eigenvectors of the saddle points. Of these mixed-mode
SSMs, we only consider the primary one, $\mathcal{W^{\infty}}\left(E\right)$,
for the same reason as in the beam buckling problem (see Fig. \ref{fig:vk_beam}b).
The existence and uniqueness of this $\mathcal{W^{\infty}}\left(E\right)$
follows from statement (iii) of Theorem \ref{thm:main2}, given that
the Floquet multipliers of the saddle fixed point are nonresonant.

\begin{figure}[H]
\centering{}\includegraphics[width=0.9\textwidth]{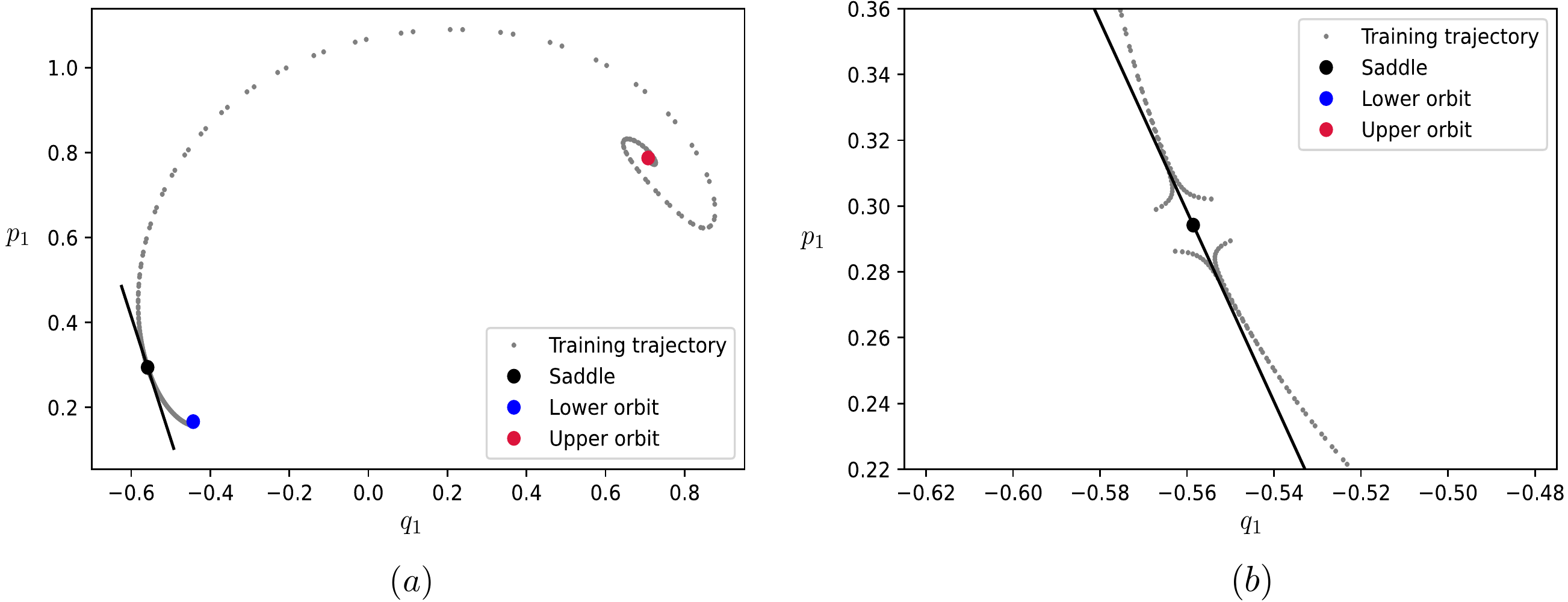}\caption{Training trajectories for constructing the 2D primary, mixed-mode
SSM, $\mathcal{W^{\infty}}\left(E\right)$, of the Poincaré map of
the Shaw\textendash Pierre model (\ref{eq:shawpierreODE}). (a) Iterates
of the four training initial conditions initialized near the saddle-type
fixed point. The unstable eigenvector $\mathbf{e}_{1}$ is shown with
the black vertical line. (b) Magnification of the region around the
saddle.\label{fig:shawpierre_training}}
\end{figure}

Since the Poincaré map is not known explicitly, we have no direct
access to the quantities featured in Section \ref{sec:Generalized-SSMs-for-continuous-DS}
for the map formulation of our main results. For a data-driven construction,
we launch four trajectories along the tangent space with initial conditions
\[
x_{0}^{(1,2,3,4)}=x_{0}\pm\delta\mathbf{e}_{1}\pm\delta\mathbf{e}_{2},
\]
where $x_{0}$ is the location of the saddle-type fixed point, $\mathbf{e}_{1}$
is its unstable eigenvector, $\mathbf{e}_{2}$ is its weakest stable
eigenvector, and $\delta=0.01.$ We use these trajectories for extracting
the primary mixed-mode SSM $\mathcal{W^{\infty}}\left(E\right)$ via
polynomial regression. In Fig. \ref{fig:shawpierre_training}, we
show the iterates of the four training trajectories under the time-$T$
map of system (\ref{eq:shawpierreODE}) via a projection onto the
$(q_{1},p_{1})$ coordinate plane.

\begin{figure}[H]
\centering{}\includegraphics[width=0.8\textwidth]{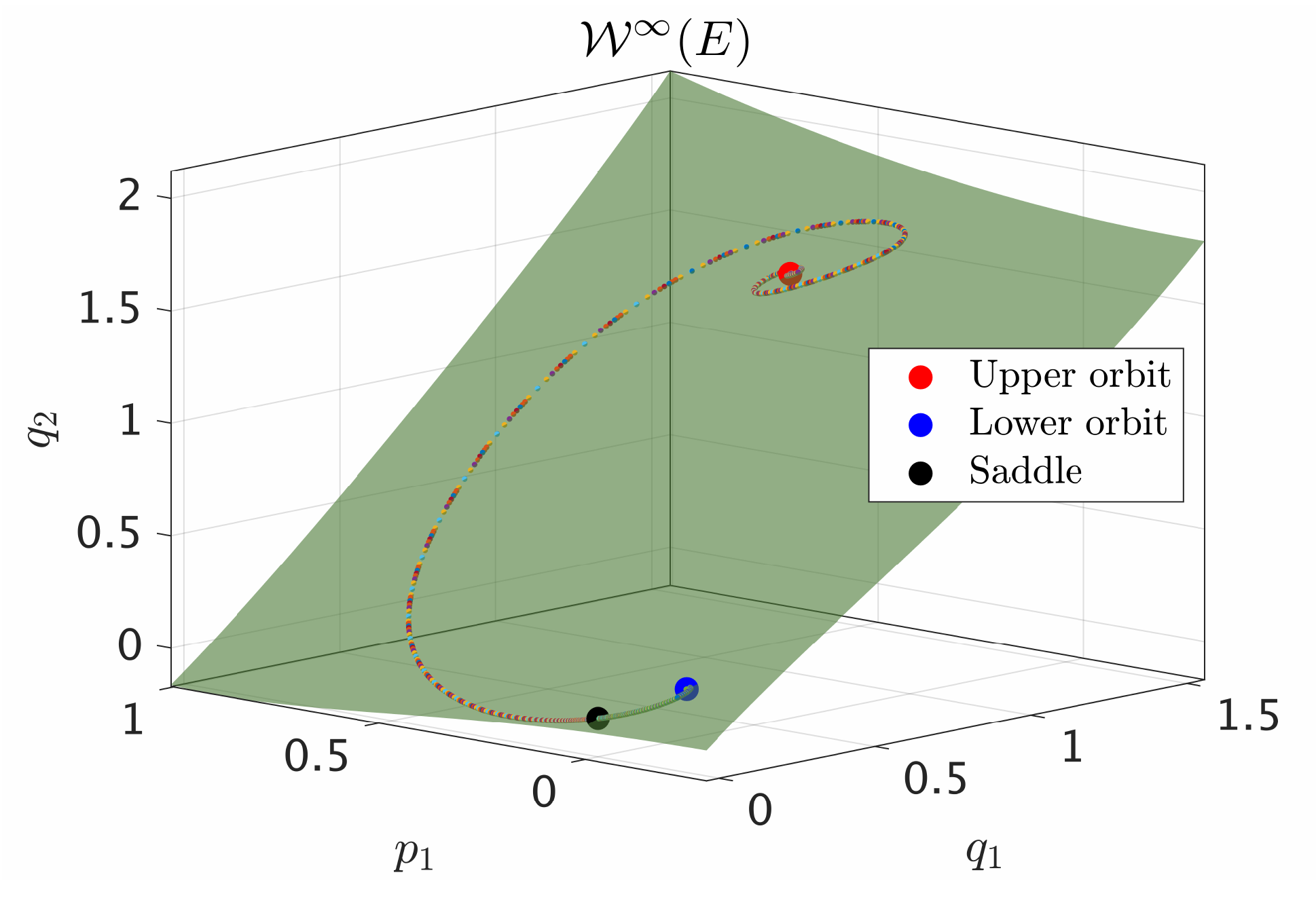}\caption{Fixed points and test trajectories on the primary, mixed-mode SSM,
$\mathcal{W^{\infty}}\left(E\right)$, of the Poincaré map of the
Shaw\textendash Pierre model (\ref{eq:shawpierreODE}). The test trajectories
are initialized near the saddle point and their iterates are denoted
by colored dots. The extracted SSM is the green surface.\label{fig:Shaw-Pierre period T map SSM}}
\end{figure}

To illustrate the accuracy of this reduced-order model on $\mathcal{W^{\infty}}\left(E\right)$,
we generate $20$ test trajectories that are initialized near the
saddle. Each initial condition has a distance of $2\cdot10^{-2}$
from the saddle point. Figure \ref{fig:Shaw-Pierre period T map SSM}
shows that these test trajectories indeed closely track the extracted
$\mathcal{W^{\infty}}\left(E\right)$ in a 3D projection from the
4D phase space of the Poincaré map. 

\begin{figure}[H]
\centering{}\includegraphics[width=1\textwidth]{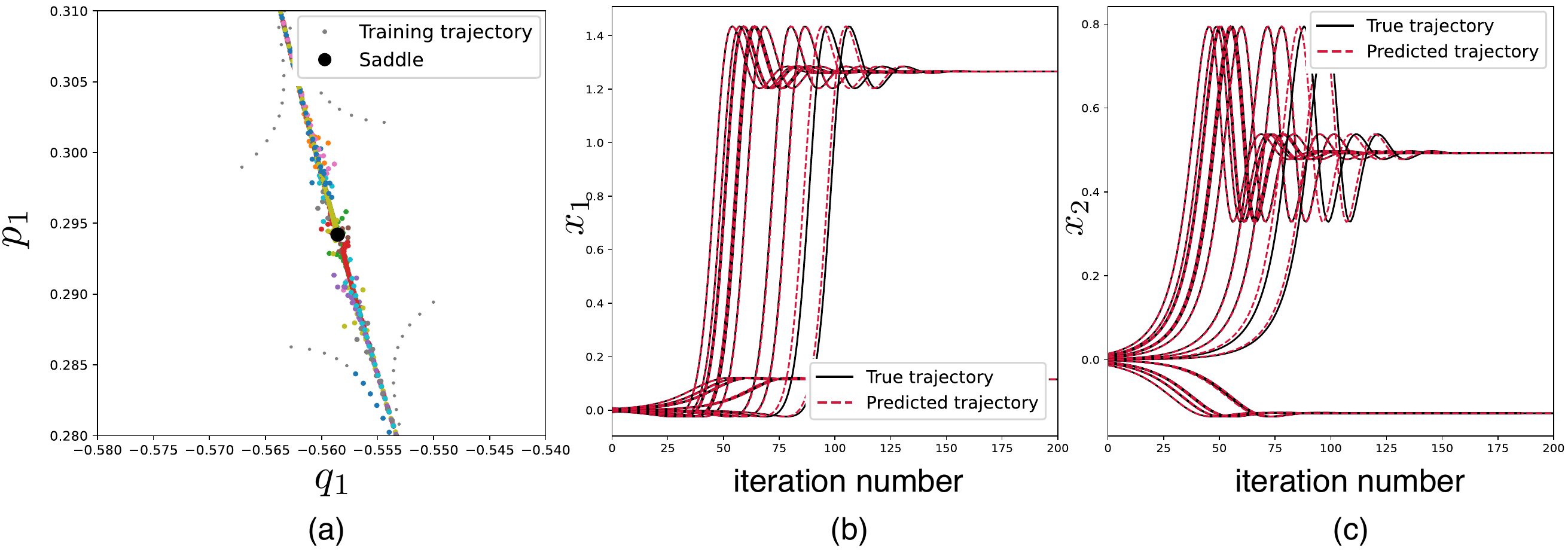}\caption{Test trajectories of the full Poincaré map and their corresponding
predictions from the 2D reduced model on $\mathcal{W^{\infty}}\left(E\right)$
for the Shaw\textendash Pierre model (\ref{eq:shawpierreODE}). (a)
Test trajectories shown as colored dots, initialized near the saddle-type
fixed point of the map. (b)-(c) Evolution of $q_{1}$ and $p_{1}$,
respectively, along with the predictions of the $5^{th}$-order $\mathcal{W^{\infty}}\left(E\right)$-reduced
model. The discrete iterates are connected into continuous lines for
visual clarity. \label{fig:shawpierretestPredictions}}
\end{figure}

We can also make predictions for the time evolution of these test
trajectories using solely the reduced discrete model. The resulting
trajectories obtained from the iterations of the full Poincaré map
and the predictions for them using the 2D reduced dynamics match closely,
as seen in Fig. \ref{fig:shawpierretestPredictions}.

\section{Conclusions}

We have extended the existing theory of (primary, or smoothest) spectral
submanifolds to the full family of continuous invariant manifolds
that emanate from a nonresonant fixed point of a generic finite-dimensional,
class $C^{\infty}$ continuous or discrete dynamical system. This
extended family now includes nonlinear continuations of spectral subspaces
of mixed stability type and/or of finite smoothness. For the latter
type of manifolds (fractional SSMs), we have derived a formal expansion
with specific fractional powers that can be inferred from the spectrum
of the linearized system. We have obtained similar results for discrete
dynamical systems defined by iterated mappings. We have illustrated
our results both for discrete and continuous dynamical systems on
the data-driven reduced-order modeling of non-linearizable behavior
in three physical systems: a canonical shear flow, a buckling beam
and a forced nonlinear rigid body system.

The main technical tool we rely on to prove these results is a classic
linearization result of \citet{sternberg58}, coupled with the explicit
construction of all continuous invariant manifolds through the origin
of a general, multi-dimensional linear system with a hyperbolic fixed
point. Neither this linearization result nor its refinements (see,
e.g. \citet{sell85}) are specific enough about the smoothness of
the linearizing transformation, except for the class of $C^{\infty}$
dynamical systems considered here, for which the linearization is
also of class $C^{\infty}$. We note that the linearizing transformation
only converges on a domain around the fixed point in which the original
system exhibits linearizable behavior, completely missing all interesting
dynamical phenomena that arise from the coexistence of several isolates
invariant sets. For this reason, we do not advocate linearization
for model reduction and only use linearization here to identify the
function basis over which an SSMs should be more globally approximated
from data. Specifically, we use truncated, finite fractional expansions
for the SSMs and its reduced dynamics in our three data-driven examples.
These truncated expansions continue to be well-defined outside the
domain of linearizable behavior as well and hence have the potential
to capture nontrivial dynamics, as we illustrate on examples. 

A conceptual limitation of the results in this paper is that it assumes
the underlying dynamical system to be finite dimensional. This is
not an issue for numerical data sets or solid or fluid mechanics,
as those are always generated from finite-dimensional discretization
of their originally infinite-dimensional governing equations. Nevertheless,
for experimental data obtained from a physical system, it would be
desirable to have a general, infinite-dimensional extension of Sternberg's
linearization results that we use. Such an extension is available
for the most frequent case in practice, wherein the spectrum of the
linearized system is either fully in the left-hand or the right-hand
side of the complex plane, i.e., satisfies the conditions of statement
(v) in Theorems 1 or 2 (see \citet{elbialy01}). This extension implies
that the finite-dimensional like-mode, primary and fractional SSMs
we have identified have infinite-dimensional limits under appropriate
spectral gap conditions. As for primary mixed-mode SSMs, their existence
and uniqueness in infinite dimensions has been known under appropriate
spectral gap conditions as long as they are also pseudo-unstable or
pseudo-stable manifolds, i.e., normally hyperbolic and purely repelling
or purely attracting, respectively (\citet{elbialy02}). Primary pseudo-stable
manifold SSMs, which also include primary like-mode SSMs, have also
been proven to exist by \citet{buza23} specifically for the Navier\textendash Stokes
equations. 

A more practical challenge for the fully data-driven application of
our results is that the fractional powers in the parametrization of
secondary SSMs depend on the real part of the spectrum of the linearization
outside the underlying spectral subspace $E$. While traces of frequencies
outside of $E$ can be extracted from the frequency analysis of the
available data, the real parts of the eigenvalues are more difficult
to identify. An option is to treat the leading-order fractional powers
as unknowns and determine them from a maximum likelihood regression.
This approach shows promise in some of our preliminary numerical tests,
but may also identify physically incorrect values for the fractional
powers in some other examples we have considered. Another option is
to estimate the real parts of the dominant eigenvalues outside $E$
from a linear data-driven method, such as DMD or EDMD. A third option
is to excite individual modes at small amplitudes and infer their
decay rates by fitting one-degree-of-freedom linear oscillators to
their forced response curves. The detailed implementation of these
strategies will be pursued in other publications.

The backbone curves and damping curves derived in Section 4 from a
fractional normal form should be helpful in explaining some unexpected
backbone curves and instantaneous damping curves observed in detailed
experimental data from gravity measurements and hydrogel oscillations
(\citet{cenedese23}, \citet{eriten22}). In those experiments, the
observed shapes of backbone curves and instantaneous damping curves
cannot be explained without including the fractional-powered terms
in formulas \eqref{eq:Omega(r)} and \eqref{eq:Omega(r)}, which provided
the original motivation of our present study. These experimental data
sets along with others will be analyzed elsewhere with the methods
developed here.

We close by noting that the fractional polynomial powers identified
in Theorems 1 and 2 generically arise in the reduced dynamics of any
non-unique invariant manifold emanating from hyperbolic fixed points.
This implies that even other data-driven modeling approaches ignoring
the existence of SSMs (e.g., EDMD and SINDy) would benefit from including
such fractional powers in the dictionary of functions that they seek
to fit to data.

\section{Supplementary Material}

In this section, we give technical details for Example 1 (Appendix
A) and Example 2 (Appendix B). We also prove Theorem 1 (Appendix C),
Propositions 1 and 2 (Appendix D), Theorem 2 (Appendix E), Propositions
3 and 4 (Appendix F), and Theorem 3 (Appendix G).\\
\\
\\
\\

\textbf{ACKNOWLEDGMENT}\\

We acknowledge helpful discussion with Mattia Cenedese and Melih Eriten.
We are also grateful to Gergely Buza for bringing Ref. \citep{buza23}
to our attention.\\

\textbf{AUTHOR DECLARATIONS} \\

\emph{Conflict of Interest }

The authors have no conflicts to disclose. \\

\emph{Author Contributions} 

George Haller: Conceptualization (lead), Methodology (lead), Formal
analysis (lead), Supervision (lead), Writing \textendash{} original
draft (lead). Bálint Kaszás: Conceptualization (equal), Formal analysis
(lead), Data Curation (equal), Writing \textendash{} original draft
(equal), Writing \textendash{} review \& editing (equal). Aihiu Liu:
Conceptualization (equal); Formal analysis (equal), Data Curation
(equal), Writing \textendash{} original draft (equal); Joar Axås:
Conceptualization (equal), Supervision (equal), Writing \textendash{}
review \& editing (equal). \\
\\

\textbf{DATA AVAILABILITY}

All code and data used in the paper will be made available upon publication
under https://github.com/haller-group.

\vfil\eject
\begin{verse}
{\huge{}Supplementary Material }{\huge\par}
\end{verse}

\section{Appendix A: Details for example (1)}

\subsection{Explicit solutions}

We seek invariant manifolds of system 
\begin{align}
\dot{x} & =x(y-b),\nonumber \\
\dot{y} & =cy(x-a),\label{eq:example-1}
\end{align}
with the parameters $a,b,c>0$ that can be written as 
\[
W(E_{1})=\{(x,h(x)):x\in\mathbb{R}\}.
\]

Substituting $y=h(x)$ into the differential equation we get the invariance
equation
\begin{equation}
x\left(h(x)-b\right)\frac{d}{dx}h(x)=ch(x)\left(x-a\right).\label{eq:separableODE}
\end{equation}

A solution to \eqref{eq:separableODE} is $h(x)=0$, which is the
analytic smooth SSM, as we will see later. Assuming that $h(x)\neq0$,
we may rearrange \eqref{eq:separableODE} to obtain
\[
\frac{c(x-a)}{x}=\frac{h'(x)[h(x)-b]}{h(x)},
\]
where differentiation with respect to $x$ is denoted by a prime.
We may integrate both sides of this equation, after which we have
\[
c(x-a\log x)=h(x)-b\log h(x)+C,
\]
for some $C\in\mathbb{R}$, a constant of integration. This is an
algebraic relationship, which already gives an implicit definition
for the manifold $y=h(x).$ However, we may rearrange this to search
for an explicit expression. With the substitution $\hat{h}=-\frac{h(x)}{b},$
we can obtain

\begin{equation}
-\frac{1}{b}e^{-\frac{cx}{b}+\frac{C}{b}}x^{\frac{ca}{b}}=e^{\hat{h}}\hat{h}(x).\label{eq:implicit equation for hh}
\end{equation}

Note that the solution $h(x)\equiv0$ is a singular solution of \eqref{eq:separableODE},
because it cannot be obtained from \eqref{eq:implicit equation for hh}
for any choice of the parameter $C$: for $x\neq0$, the left-hand
side is never zero, which means $h\neq0$. The implicit relation \eqref{eq:implicit equation for hh}
can be solved for $h(x)$ using the Lambert $W$ function, which is
also called the product logarithm function (see, e.g., \citet{bronshtein15}).
This is defined through the inverse relation
\[
z=we^{w},
\]
which can be solved for $w$ if $z\geq-\frac{1}{e}$. In this case,
the solution is $w=W_{m}(z)$ for $m=-1$ or $m=0$, denoting the
two branches of the $W$ function. This allows us to express the family
of invariant manifolds, parametrized by the constant of integration
$C$, as 
\begin{equation}
h(x)=-bW_{m}\left(-\frac{1}{b}e^{-\frac{cx}{b}+\frac{C}{b}}x^{\frac{ca}{b}}\right).\label{eq:exact planar model fractional graph family}
\end{equation}

To decide which branch of the $W$ function we need, we note that
at $x=0$, the $m=-1$ branch diverges. The $m=0$ branch on the other
hand is well defined at $x=0$ for all values of $C$. Since we are
interested in an invariant manifold that contains the fixed point
at the origin, we will use the principal branch $W_{0}$ of the $W$
function.

We can also determine the value of the constant $C$, if we require
$h(x)$ to also contain the fixed point at $(a,b)$. The relationship
$h(a)=b$ implies

\[
W_{0}\left(-\frac{1}{b}e^{-\frac{ca}{b}+\frac{C^{*}}{b}}a^{\frac{ca}{b}}\right)=-1,
\]
which can be solved for $C^{*}$ as 
\[
C^{*}=b\log\left(\frac{b}{e}a^{-\frac{ca}{b}}\right)+\frac{ca}{b},
\]
to arrive at the particular member of the family of the invariant
manifolds that contains both fixed points of system \eqref{eq:example-1}.
This finally gives the exact reduced order model 
\[
\dot{x}=-bx\left[W_{0}\left(-\frac{1}{b}e^{-\frac{cx}{b}+\frac{C^{*}}{b}}x^{\frac{ca}{b}}\right)-1\right]
\]
 on the attracting SSM that connects the two fixed points.

\subsection{POD, DMD and Koopman decompositions}

One can consider the proper orthogonal decomposition (POD) on this
problem to see if a projection to the dominant POD mode yields good
reduced-order dynamics. If POD yields a dominant POD mode along a
yet unknown line $y=Kx,$ then projection of the right-hand side of
(\ref{eq:example-1}) to this line gives the equation
\begin{align}
 & =\left.\frac{\left(x(y-b),cy(x-a)\right)\cdot\left(1,K\right)}{\sqrt{K^{2}+1}}\right|_{y=Kx}=\left.\frac{\left(x(Kx-b),cKx(x-a)\right)\cdot\left(1,K\right)}{\sqrt{K^{2}+1}}\right|_{y=Kx}\label{eq:pod rhs}\\
 & =\frac{K+cK^{2}}{\sqrt{K^{2}+1}}x\left[x-\frac{b+cK^{2}a}{K+cK^{2}}\right].\nonumber 
\end{align}
On the other hand, projecting the left-hand-side of \eqref{eq:example-1}
results in 
\begin{equation}
\left.\frac{\left(\dot{x},\dot{y}\right)\cdot\left(1,K\right)}{\sqrt{K^{2}+1}}\right|_{y=Kx}=\left.\frac{\dot{x}+K\dot{y}}{\sqrt{K^{2}+1}}\right|_{y=Kx}.\label{eq:pod lhs}
\end{equation}
Equating \eqref{eq:pod lhs} and \eqref{eq:pod rhs}, we obtain the
POD-based model 
\begin{equation}
\dot{x}=\frac{K+cK^{2}}{1+K^{2}}x\left[x-\frac{b+cK^{2}a}{K+cK^{2}}\right].\label{eq:pod reduced model}
\end{equation}

\begin{figure}[h]
\centering{}\includegraphics[width=1\textwidth]{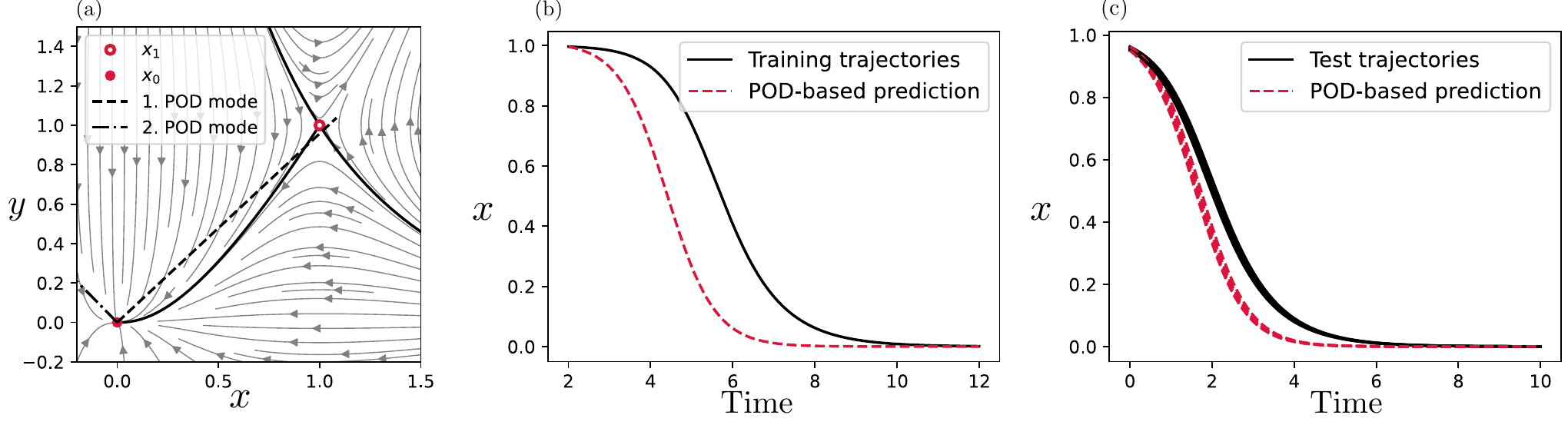}\caption{POD-based reduced order model example (\ref{eq:example-1}). Panel
(a) shows the phase space of system \eqref{eq:example-1} along with
the numerically computed POD modes (dashed lines). The slope of the
first POD mode is $K=0.95764$. In panels (b) and (c) the predicted
dynamics is compared with the true dynamics for training and testing
trajectories, respectively.\label{fig:POD for example 1}}
\end{figure}

This POD model will always have a fixed point at the origin and at
another location, which is not necessarily $x=b$. Figure \ref{fig:POD for example 1}
shows the results of POD-based reduced order modeling for example
(\ref{eq:example-1}) with $a=b=1$ and $c=2.5$. The model is not
accurate for the details of the transition orbit but has the correct
asymptotic behavior in forward time. Note, however, that deriving
this reduced-order model required the knowledge of the full set of
equations of the dynamical system, and hence this approach is inapplicable
in a data-driven setting.

\begin{figure}[H]
\centering{}\includegraphics[width=0.5\textwidth]{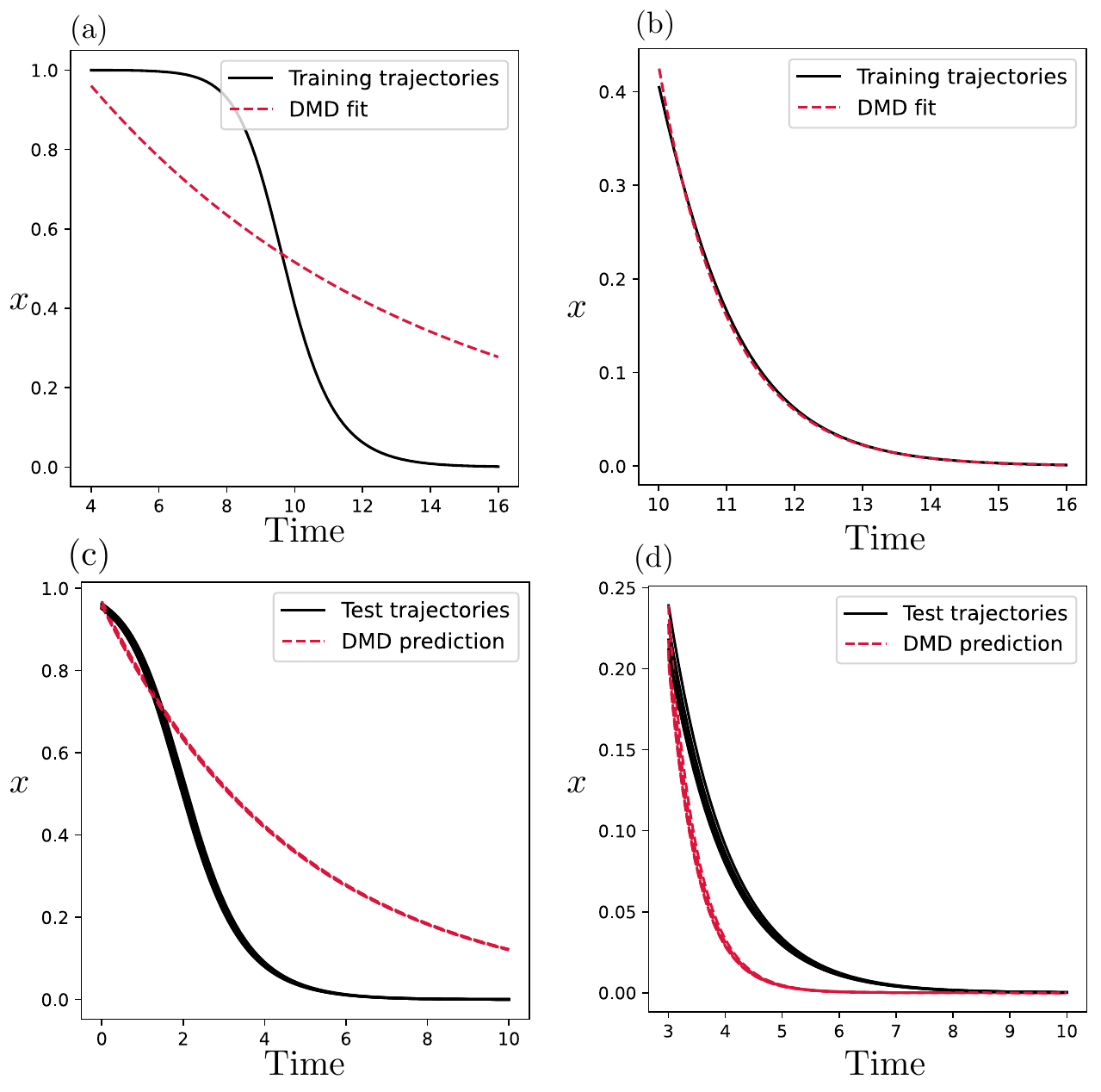}\caption{DMD-based reduced order models. (a) DMD fit to a full transition trajectory.
(b) DMD fit to a subset of the same transition trajectory closer to
the origin. (c) Predictions (red) from the DMD fit in (a) for test
trajectories of the full system. (d) Predictions (red) from the DMD
fit in (b) for truncated test trajectories of the full system.\label{fig:DMD fit for example 1}}
\end{figure}

Dynamic mode decomposition (DMD) is fully data driven and will fit
a linear dynamical system to trajectories. This system will clearly
have a stable node for small amplitude training data and a saddle
for larger-amplitude training data. Figure \ref{fig:DMD fit for example 1}
illustrates the prediction quality of two different DMD models: one
fitted to the full transition trajectory (Fig. \ref{fig:DMD fit for example 1}a)
and one fitted to a truncated subset of the transition trajectory
closer to the origin (Fig. \ref{fig:DMD fit for example 1}b). The
second fit can capture the asymptotic behavior (see the decay to the
origin in Fig. \ref{fig:DMD fit for example 1}d) but neither fit
can make reliable predictions outside a small neighborhood of the
origin, as expected. 

A numerically constructed Koopman decomposition will behave similarly,
given that it can only capture linearizable dynamics. Theoretically,
the two $C^{1}$-unique principal Koopman eigenfunctions (see \citet{kvalheim21})
associated with the domain of attraction of the origin are necessarily
the ones whose zero set is the $x$ axis. On this, the dynamics is
indeed globally linear and hence does not capture the nontrivial fixed
point and the heteroclinic connection of interest.

\subsection{Reduced-order models on fractional SSMs }

By Theorem 1 of the main text, there exist a family $\mathcal{W}\left(E\right)$
of continuous SSMs tangent to the spectral subspace $E=\left\{ \left(x,y\right)\in\mathbb{R}^{2}:\,y=0\right\} $
in example (\ref{eq:example-1}). By Proposition 1, for an arbitrary
integer $\mathcal{K}\geq2$, the $x\geq0$ branches of members of
this family can be written as 
\begin{align}
y & =C_{\left(0,1\right)}(x)x^{\frac{k_{4}ac}{b}}+\sum_{2\leq k_{1}+k_{4}\frac{ac}{b}\leq\mathcal{K}}C_{\mathbf{k}}(x)x^{k_{1}+k_{4}\frac{ac}{b}}+o\left(\left|x\right|^{\mathcal{K}}\right),\qquad\mathbf{k}=\left(k_{1},k_{4}\right),\label{eq:1D SSM in first example}
\end{align}
with the arbitrary, piecewise constant functions $C_{\mathbf{k}}(x)\in\mathbb{R}$
satisfying
\[
C_{\mathbf{k}}(x)=\left\{ \begin{array}{c}
C_{\mathbf{k}}^{+}\qquad x>0,\\
\\
C_{\mathbf{k}}^{-}\qquad x\leq0.
\end{array}\right.
\]
 Restriction of example (\ref{eq:example-1}) to this invariant family
of graphs gives the reduced-order model family
\begin{equation}
\dot{x}=x\left[-b+C_{\left(0,1\right)}\left|x\right|^{k_{4}\frac{ac}{b}}+\sum_{2\leq k_{1}+k_{4}\frac{ac}{b}\leq\mathcal{K}}C_{\mathbf{k}}(x)\left|x\right|^{k_{1}+k_{4}\frac{ac}{b}}+o\left(\left|x\right|^{\mathcal{K}}\right)\right].\label{eq:1D reduced model in first example}
\end{equation}
 One member of this family is the primary SSM $\mathcal{W}^{\infty}\left(E\right)=E$
satisfying $C_{\mathbf{k}}=0$ for all $\mathbf{k}$. Another member
is the heteroclinic connection $\mathcal{W}^{*}\left(E\right)$ that
is generally of class $C^{\mathrm{Int\left[ac/b\right]}}$ smooth.
In the data-driven setting, the parameters in the formulas (\ref{eq:1D SSM in first example})-(\ref{eq:1D reduced model in first example})
can be determined by regression to training trajectories that are
in a neighborhood of the heteroclinic orbit. 

\begin{figure}[H]
\centering{}\includegraphics[width=1\textwidth]{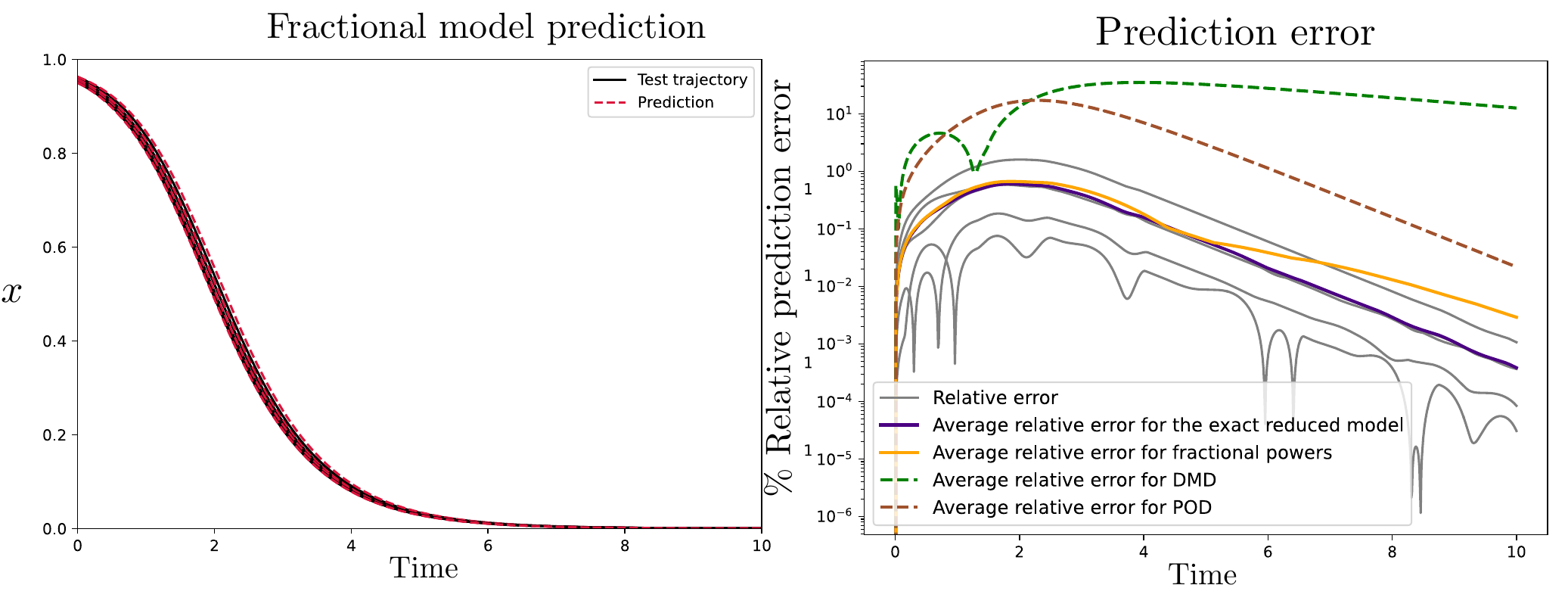}\caption{Fractional SSM-based reduced-order model compared with other model
reduction methods for example (1). In the left panel we show the time
evolution of 5 test trajectories along with the predicted time evolution
based on the 5th order model. In the right panel, we show the prediction
errors from data-driven fractional SSM reduction (yellow), DMD (green)
and POD reduction (brown). For comparison, we also show the prediction
for these nearby trajectories from the exact reduced mode obtained
by restricting the full dynamics of example (1) to the heteroclinic
connection obtained from formula (\ref{eq:exact planar model fractional graph family})
. \label{fig:Fractional-SSM-for-simple-exampe-1}}
\end{figure}

For comparison with the POD and DMD-based reduced-order models, we
show the predictions obtained by a data-driven reduced-order model
for example (\ref{eq:example-1}). In particular, using a single trajectory
initialized close to the heteroclinic orbit, we compute the fractional
SSM (\ref{eq:1D SSM in first example}) up to order $\mathcal{K}=5$
via regression. Restricting to this fractional SSM then yields the
reduced-order model (\ref{eq:1D reduced model in first example})
with specific coefficients up to order $\mathcal{K}=5$ . Figure \ref{fig:Fractional-SSM-for-simple-exampe-1}
shows the close agreement between the predicted and true time evolutions
of 5 test trajectories initialized close to the saddle. We also show
the close agreement with the exact, analytically obtained reduced-order
model, as well as the substantially larger errors produced by DMD
and POD.

\section{Appendix B: Details for example (2)}

\subsection{No 2D pseudo-unstable manifold }

For the motivating example 
\begin{align}
\dot{x}_{1} & =x_{2},\nonumber \\
\dot{x}_{2} & =-x_{1}-cx_{2}+x_{1}^{3},\nonumber \\
\dot{x}_{3} & =-kx_{3}+kax_{1}^{2}+2ax_{1}x_{2},\label{eq:3D motivating example-1}
\end{align}
of the main text, we first recall from \citet{delallave95} the precise
statement of Irwin's theorem for the existence of a pseudo-unstable
manifold for a $C^{r}$ diffeomorphism $f\colon\mathcal{B}\to\mathcal{B}$
on a Banach space $\mathcal{B}$ with $f(0)=0$. Assume that $\mathcal{B}=S\oplus U$
is an invariant direct-sum decomposition of $\mathcal{B}$ into the
closed subspaces $S$ and $U$ such that for some number $a<1$, we
have\linebreak{}
\begin{equation}
\left\Vert Df^{-1}(0)\vert_{U}\right\Vert <a,\qquad\left\Vert Df(0)\vert_{S}\right\Vert <a^{-r}.\label{eq:Irwin conditions-1-1}
\end{equation}
Then, for small enough $\epsilon>0$, the set
\[
W(S)=\left\{ \left\Vert x\right\Vert <\epsilon\colon\,\sup_{n\geq0}\left|f^{-n}(x)\right|a^{-n}<\infty\right\} 
\]
is a $C^{r}$manifold that is tangent to $S$ at $x=0$. This manifold
is, therefore, composed of orbits whose growth rates are bounded from
above by $a^{n}$. 

In example (\ref{eq:3D motivating example-1}), let us select $U$
to be the spectral subspace corresponding to the eigenvalues $\left\{ -1,+1\right\} $
and $S$ the spectral subspace corresponding to the eigenvalue $-k$,
with $k=\sqrt{2}/2$.  Then, noting that 
\[
Df(0)=\exp\left[\begin{array}{ccc}
0 & 1 & 0\\
-1 & 0 & 0\\
0 & 0 & -\sqrt{2}/2
\end{array}\right]
\]
holds for the time-one flow map $f$ of system (\ref{eq:3D motivating example-1}),
the requirements in (\ref{eq:Irwin conditions-1-1}) become
\[
\left\Vert Df^{-1}(0)\vert_{U}\right\Vert =e<a,\qquad\left\Vert Df(0)\vert_{S}\right\Vert =\frac{1}{e^{\sqrt{2}/2}}<\frac{1}{a^{r}}.
\]
These conditions are satisfied if we find $a,r\in\mathbb{R}^{+}$
such that 
\[
e<a<e^{\frac{\sqrt{2}/2}{r}}.
\]
Such an $a<1$, however, will not exists. In other words, the subspace
$U$ is not a pseudo-unstable subspace.

\subsection{2D mixed-mode SSM}

Theorem 1 guarantees the existence of infinitely many continuous,
mixed-mode SSMs tangent to the spectral subspace $E=\left\{ x\in\mathbb{R}^{3}:\,x_{3}=0\right\} $.
Precisely one of them, $\mathcal{W}^{\infty}\left(E\right),$ is of
class $C^{\infty}$. Indeed, we identified this unique, smooth invariant
manifold as $\mathcal{M}=\left\{ x\in\mathbb{R}^{3}\colon x_{3}=ax_{1}^{2}\right\} $
in our description of the example. This manifold can also be approximated
purely from data using the \emph{SSMLearn} package, as we noted in
Remark 1 after Theorem 1. We do not present this numerical identification
here, as we identify a similar $\mathcal{W}^{\infty}\left(E\right)$
from data via \emph{SSMLearn} for our more challenging buckling beam
example. 

\section{Appendix C: Proof of Theorem 1\label{sec:Proof of main theorem for continuous sytstems}}

\subsection{Secondary SSMs in the linearized system }

Recalling the original nonlinear system
\begin{equation}
\dot{x}=\mathcal{A}x+f(x),\qquad x\in\mathbb{R}^{n},\quad\mathcal{A}\in\mathbb{R}^{n\times n},\quad f=\mathcal{O}\left(\left|x\right|^{2}\right)\in C^{\infty},\qquad1\leq n<\infty,\label{eq:original system-3}
\end{equation}
we first describe all invariant graphs over the spectral subspace
$E$ in the linear part 
\begin{equation}
/\dot{y}=\mathcal{A}y,\qquad y\in\mathbb{R}^{n},\label{eq:linearized system}
\end{equation}
of (\ref{eq:original system-3}). We use here the notation
\[
y=\left(u,\left(a_{1},b_{1}\right),\ldots,\left(a_{q},b_{q}\right),v,\left(c_{1},d_{1}\right),\ldots,\left(c_{s},d_{s}\right)\right)\in\mathbb{R}^{p}\times\underbrace{\mathbb{R}^{2}\times\ldots\times\mathbb{R}^{2}}_{q}\times\mathbb{R}^{r}\times\underbrace{\mathbb{R}^{2}\times\ldots\times\mathbb{R}^{2}}_{s}
\]
for the phase space variable of the linear equation (\ref{eq:linearized system})
to distinguish its solutions from those of its nonlinear counterpart
(\ref{eq:original system-3}), for which we used the variable $x$
of the dame dimension. We also recall from the main text that $\mathcal{A}$
can be assumed in the form
\begin{equation}
\mathcal{A}=\left(\begin{array}{cccccccc}
A\\
 & B_{1}\\
 &  & \ddots &  &  &  & 0\\
 &  &  & B_{q}\\
 &  &  &  & C\\
 &  &  &  &  & D_{1}\\
 &  & 0 &  &  &  & \ddots\\
 &  &  &  &  &  &  & D_{s}
\end{array}\right),\label{eq:linsystem1-1-2}
\end{equation}
with
\begin{align}
A & =\mathrm{diag}\left[\lambda_{1},\ldots,\lambda_{p}\right]\in\mathbb{R}^{p\times p},\quad B_{k}=\left(\begin{array}{cc}
\alpha_{k} & -\omega_{k}\\
\omega_{k} & \alpha_{k}
\end{array}\right)\in\mathbb{R}^{2\times2},\quad k=1,\ldots,q,\nonumber \\
C & =\mathrm{diag}\left[\kappa_{1},\ldots,\kappa_{r}\right]\in\mathbb{R}^{r\times r},\quad D_{m}=\left(\begin{array}{cc}
\beta_{m} & -\nu_{m}\\
\nu_{m} & \beta_{m}
\end{array}\right)\in\mathbb{R}^{2\times2},\quad m=1,\ldots,s.\label{eq:linsystem1-2-2}
\end{align}
We also recall that the hyperbolicity condition
\begin{equation}
0\notin\mathrm{Re}\left[\mathrm{spect}\left(\mathcal{A}\right)\right],\label{eq:hyperbolicity-1}
\end{equation}
and the nonresonance conditions
\begin{equation}
\lambda_{j}\neq\sum_{k=1}^{n}m_{k}\lambda_{k},\quad m_{k}\in\mathbb{N},\quad\sum_{k=1}^{n}m_{k}\geq2,\quad j=1,\ldots,n\label{eq:nonresonance conditions-2}
\end{equation}
are assumed to hold.

To find all continuous invariant graphs in the linear system (\ref{eq:linearized system}),
we first introduce multi-dimensional polar coordinates $(r,\phi)$
and $(\rho,\theta)$ instead of the $(a,b)$ and $(c,d)$ coordinate
pairs by letting\textrm{
\begin{align}
a_{k} & =r_{k}\cos\phi_{k},\quad\,\,\,\,\,b_{k}=r_{k}\sin\phi_{k},\nonumber \\
c_{m} & =\rho_{m}\cos\theta_{m},\quad d_{m}=\rho_{m}\sin\theta_{m}.\label{eq:polar coordinates}
\end{align}
}This transforms the $\left(a,b,c,d\right)$ component of the linear
system (\ref{eq:linearized system}) to the polar form
\begin{align}
\dot{r}_{k} & =\alpha_{k}r_{k},\quad\dot{\phi}_{k}=\omega_{k},\quad\,\,\,\,\,k=1,\ldots,q,\label{eq:linsystem2+}\\
\dot{\rho}_{m} & =\beta_{m}\rho_{m},\quad\dot{\theta}_{m}=\nu_{m},\quad m=1,\ldots,s,\label{eq:linsystem4+}
\end{align}

Using the $(u,r,\phi)$ coordinates within $E$ and the lack of explicit
dependence of the right-hand side of eqs. (\ref{eq:linsystem2+})-(\ref{eq:linsystem4+})
on $\phi$, we can write the $v_{\ell}(t)$ component of a trajectory
contained in an invariant, continuous graph over $E$ as

\begin{equation}
v_{\ell}(t)=F_{\ell}\left(u(t),r(t)\right),\quad\ell=1,\ldots,r,\label{eq:inv. condition 1}
\end{equation}
for some continuous function $F\colon\mathbb{R}^{p}\times\mathbb{R}^{q}\to\mathbb{R}^{r}$.
Differentiating the invariance condition (\ref{eq:inv. condition 1})
with respect to time and using eqs. (\ref{eq:linsystem1-1-2}) and
(\ref{eq:linsystem1-2-2}) gives the linear system of PDEs 
\[
\kappa_{\ell}F_{\ell}=\sum_{j=1}^{p}\partial_{u_{j}}F_{\ell}\lambda_{j}u_{j}+\sum_{k=1}^{q}\partial_{r_{k}}F_{\ell}\alpha_{k}r_{k},\quad\ell=1,\ldots,r.
\]
The general solution of this PDE system for $u_{j}\in\mathbb{R}$
and $r_{k}\geq0$ is 
\begin{align}
F_{\ell} & =\sum_{j=1}^{p}\left[K_{\ell j}^{+}H(u_{j})+K_{\ell j}^{-}H(-u_{j})\right]\left|u_{j}\right|^{\frac{\kappa_{\ell}}{\lambda_{j}}}+\sum_{k=1}^{q}L_{\ell k}r_{k}^{\frac{\kappa_{\ell}}{\alpha_{k}}}\nonumber \\
 & =\sum_{j=1}^{p}\left[K_{\ell j}^{+}H(u_{j})+K_{\ell j}^{-}H(-u_{j})\right]\left|u_{j}\right|^{\frac{\kappa_{\ell}}{\lambda_{j}}}+\sum_{k=1}^{q}L_{\ell k}\left(a_{k}^{2}+b_{k}^{2}\right)^{\frac{\kappa_{\ell}}{2\alpha_{k}}},\label{eq:fuv}
\end{align}
 where 
\begin{equation}
H(\zeta)=\left\{ \begin{array}{c}
1\qquad\zeta>0,\\
\\
0\qquad\zeta\geq0,
\end{array}\right.\label{eq:heavyside}
\end{equation}
is the Heaviside step function and $K_{\ell j}^{\pm}$ and $L_{\ell k}$
are arbitrary real constants subject to the restriction
\begin{align}
K_{\ell j}^{\pm}=0,\quad\mathrm{if}\quad\frac{\kappa_{\ell}}{\lambda_{j}}<0,\nonumber \\
L_{\ell k}=0,\quad\mathrm{if}\quad\frac{\kappa_{\ell}}{\alpha_{k}}<0.\label{eq:contraints on K and L}
\end{align}
The need for using the piecewise constant function (\ref{eq:heavyside})
and the restrictions on the parameters in (\ref{eq:contraints on K and L})
are explained in the remarks after Theorem 1 in the main text. 

Next, using again the $(u,r,\phi)$ coordinates within $E$ and the
$S^{1}$ symmetry of the equations (\ref{eq:linsystem2+}) and (\ref{eq:linsystem4+})
with respect to $\phi$, we can write the $\left(\rho_{m},\theta_{m}\right)$
component of a trajectory contained in an invariant graph over $E$
as

\[
\rho_{m}(t)=G_{m}^{\rho}\left(u(t),r(t)\right),\quad\theta_{m}(t)=G_{m}^{\theta}\left(u(t),r(t)\right),\quad m=1,\ldots,s,
\]
for some continuous functions $G^{\rho}\colon\mathbb{R}^{p}\times\mathbb{R}^{q}\to\mathbb{R}^{s}$
and $G^{\theta}\colon\mathbb{R}^{p}\times\mathbb{R}^{q}\to\mathbb{\mathbb{T}}^{s}$.
Differentiating this invariance condition with respect to time and
using eqs. (\ref{eq:linearized system}), (\ref{eq:linsystem2+})
and (\ref{eq:linsystem4+}) gives now the set of $2m$ linear PDEs
\begin{align*}
\beta_{m}G_{m}^{\rho} & =\sum_{j=1}^{p}\partial_{u_{j}}G_{m}^{\rho}\lambda_{j}u_{j}+\sum_{k=1}^{q}\partial_{r_{k}}G_{m}^{\rho}\alpha_{k}r_{k},\\
\nu_{m} & =\sum_{j=1}^{p}\partial_{u_{j}}G_{m}^{\theta}\lambda_{j}u_{j}+\sum_{k=1}^{q}\partial_{r_{k}}G_{m}^{\theta}\alpha_{k}r_{k}.
\end{align*}
The general solutions of these PDEs are
\begin{align}
G_{m}^{\rho} & =\sum_{j=1}^{p}\left[M_{mj}^{+}H(u_{j})+M_{mj}^{-}H(-u_{j})\right]\left|u_{j}\right|^{\frac{\beta_{m}}{\lambda_{j}}}+\sum_{k=1}^{q}P_{mk}r_{k}^{\frac{\beta_{m}}{\alpha}}\nonumber \\
 & =\sum_{j=1}^{p}\left[M_{mj}^{+}H(u_{j})+M_{mj}^{-}H(-u_{j})\right]\left|u_{j}\right|^{\frac{\beta_{m}}{\lambda_{j}}}+\sum_{k=1}^{q}P_{mk}\left(a_{k}^{2}+b_{k}^{2}\right)^{\frac{\beta_{m}}{2\alpha_{k}}},\label{eq:grho}\\
G_{m}^{\theta} & =N_{m}+\sum_{j=1}^{p}\frac{\nu_{m}}{\lambda_{j}}\log\left|u_{j}\right|+\sum_{k=1}^{q}\frac{\nu_{m}}{\alpha_{k}}\log r_{k}\nonumber \\
 & =N_{m}+\sum_{j=1}^{p}\frac{\nu_{m}}{\lambda_{j}}\log\left|u_{j}\right|+\sum_{k=1}^{q}\frac{\nu_{m}}{2\alpha_{k}}\log\left(a_{k}^{2}+b_{k}^{2}\right).\label{eq:gtheta}
\end{align}
 Here $M_{mj}^{\pm}$, $N_{m}$ and $P_{mk}$ are arbitrary real constants
subject to the constraints
\begin{align}
M_{mj}^{+} & =M_{mj}^{-}=0,\quad\mathrm{if}\quad\frac{\beta_{m}}{\lambda_{j}}<0,\nonumber \\
P_{mk} & =0,\quad\mathrm{if}\quad\frac{\beta_{m}}{\alpha_{k}}<0,\label{eq:contraints on M and P}
\end{align}
that are needed for the same reason as those in eq. (\ref{eq:contraints on K and L}).

Using the formulas (\ref{eq:polar coordinates}), (\ref{eq:fuv}),
(\ref{eq:grho}) and (\ref{eq:gtheta}), we obtain the following general
form for the full family of continuous invariant graphs over $E$
that contain the origin:
\begin{align}
v_{\ell} & =\mathcal{V}_{\ell}(u,a,b;K,L)\nonumber \\
 & :=\sum_{j=1}^{p}\left[K_{\ell j}^{+}H(u_{j})+K_{\ell j}^{-}H(-u_{j})\right]\left|u_{j}\right|^{\frac{\kappa_{\ell}}{\lambda_{j}}}+\sum_{k=1}^{q}L_{\ell k}\left(a_{k}^{2}+b_{k}^{2}\right)^{\frac{\kappa_{\ell}}{2\alpha_{k}}},\quad\ell=1,\ldots,r,\nonumber \\
\left(\begin{array}{c}
c_{m}\\
d_{m}
\end{array}\right) & =\left(\begin{array}{c}
\mathcal{C}{}_{m}(u,a,b)\\
\mathcal{D}_{m}(u,a,b)
\end{array}\right)\nonumber \\
 & :=\left(\sum_{j=1}^{p}\left[M_{mj}^{+}H(u_{j})+M_{mj}^{-}H(-u_{j})\right]\left|u_{j}\right|^{\frac{\beta_{m}}{\lambda_{j}}}+\sum_{k=1}^{q}P_{mk}\left(a_{k}^{2}+b_{k}^{2}\right)^{\frac{\beta_{m}}{2\alpha_{k}}}\right)\label{eq:general_fractional_SSM1-3}\\
 & \,\,\,\,\,\,\times\left(\begin{array}{c}
\cos\left(N_{m}+\sum_{j=1}^{p}\frac{\nu_{m}}{\lambda_{j}}\log\left|u_{j}\right|+\sum_{k=1}^{q}\frac{\nu_{m}}{2\alpha_{k}}\log\left(a_{k}^{2}+b_{k}^{2}\right)\right)\\
\sin\left(N_{m}+\sum_{j=1}^{p}\frac{\nu_{m}}{\lambda_{j}}\log\left|u_{j}\right|+\sum_{k=1}^{q}\frac{\nu_{m}}{2\alpha_{k}}\log\left(a_{k}^{2}+b_{k}^{2}\right)\right)
\end{array}\right),\quad m=1,\ldots,s.\nonumber 
\end{align}
For SSMs that are symmetric in the $u_{j}$ variable, the formulas
simplify due to the relations 
\[
K_{\ell j}^{+}H(u_{j})+K_{\ell j}^{-}H(-u_{j})=K_{\ell j},\quad M_{mj}^{+}H(u_{j})+M_{mj}^{-}H(-u_{j})=M_{mj}.
\]

For fixed values of the parameter vectors $K,L,M,N$ and $P$, each
of the scalar equations 
\begin{align*}
v_{\ell} & =\mathcal{V}_{\ell}(u,a,b),\quad\ell=1,\ldots,r,\\
c_{m} & =\mathcal{C}_{m}(u,a,b),\quad m=1,\ldots,s,\\
d_{m} & =\mathcal{D}_{m}(u,a,b),\quad m=1,\ldots,s,
\end{align*}
define a codimension-one invariant surface in the phase space of the
linearized system (\ref{eq:linearized system}). The transverse intersections
of these $2s+r$ hypersurfaces are the individual SSMs, each of which
has codimension $2s+r$ and dimension $p+2q$, which is the same as
the dimension of the underlying spectral subspace $E$. Note that
under the nonresonance conditions, \eqref{eq:nonresonance conditions-2},
the only analytic invariant graphs of the linear system over $E$
are those satisfying the SSMs described by analytic graphs  with $K,L,M,N,P=0$.
All other invariant graphs are only finitely many times differentiable
at the origin.

\subsection{Full family of SSMs in the nonlinear system (\ref{eq:original system-3})}

The remaining question is if the full family of invariant graphs (\ref{eq:general_fractional_SSM1-3})
of the linearized system (\ref{eq:linearized system}) gives rise
to a similar family of invariant manifolds in the original nonlinear
system (\ref{eq:original system-3}). Under the assumption of hyperbolicity
(eq. (\ref{eq:hyperbolicity-1})) and nonresonance (eq. (\ref{eq:nonresonance conditions-2}))
for $\mathcal{A}$, \citet{sternberg58} proves the existence of a
locally defined $C^{\infty}$ linearizing diffeomorphism for system
(\ref{eq:original system-3}) with a formal (i.e., not necessarily
convergent) Taylor expansion
\begin{equation}
x=y+h(y)=y+\sum_{\left|\mathbf{k}\right|\geq2}H_{\mathbf{k}}y^{\mathbf{k}},\qquad\mathbf{k}\in\mathbb{N}^{n},\quad\left|\mathbf{k}\right|=\sum_{j=1}^{n}k_{j},\quad y^{\mathbf{k}}:=y_{1}^{k_{1}}y_{2}^{k_{2}}\cdots y_{n}^{k_{n}},\quad h_{\mathbf{k}}\in\mathbb{R}^{n},\label{eq:linearizing transformation}
\end{equation}
that transforms the full system (\ref{eq:original system-3}) into
its linear part (\ref{eq:linearized system}). 

Recalling that $y=(u,a,b,v,c,d)^{\mathrm{T}}$, we obtain from (\ref{eq:linearizing transformation})
\begin{align}
x & =y+\sum_{\left|\mathbf{k}\right|\geq2}H_{\mathbf{k}}u^{\mathbf{k}_{1}}a^{\mathbf{k}_{2}}b^{\mathbf{k}_{3}}v^{\mathbf{k}_{4}}c^{\mathbf{k}_{5}}d^{\mathbf{k}_{6}},\nonumber \\
 & \,\,\,\,\,\,\mathbf{k}_{1}\in\mathbb{N}^{p},\quad\mathbf{k}_{2},\mathbf{k}_{3}\in\mathbb{N}^{q},\quad\mathbf{k}_{4}\in\mathbb{N}^{r},\quad\mathbf{k}_{5},\mathbf{k}_{6}\in\mathbb{N}^{s},\quad\mathbf{k}=\left(\mathbf{k}_{1},\mathbf{k}_{2},\mathbf{k}_{3},\mathbf{k}_{4},\mathbf{k}_{5},\mathbf{k}_{6}\right).\label{eq:linearizing transformation-2}
\end{align}
Then restricting this mapping to the invariant graph family (\ref{eq:general_fractional_SSM1-3}),
we obtain the full family of SSMs, $\mathrm{SSM}_{KLMNP}(E)$, of
the nonlinear system (\ref{eq:original system-3}) in the parametrized
form
\begin{align}
x & =\mathrm{\mathrm{SSM}}(u,a,b)\nonumber \\
 & :=\left(u,\left(a,b\right),\mathcal{V},(\mathcal{C},\mathcal{D})\right)^{\mathrm{T}}+\sum_{\left|\mathbf{k}\right|\geq2}H_{\mathbf{k}}u^{\mathbf{k}_{1}}a^{\mathbf{k}_{2}}b^{\mathbf{k}_{3}}\mathcal{V}^{\mathbf{k}_{4}}(u,a,b)\mathcal{C}^{\mathbf{k}_{5}}(u,a,b)\mathcal{D}^{\mathbf{k}_{6}}(u,a,b),\label{eq:SSM in original coordinates real}
\end{align}
with the coordinate components of $\mathcal{V}$, $\mathcal{C}$ and
$\mathcal{D}$ defined in (\ref{eq:general_fractional_SSM1-3}). We
obtain locally all $C^{\infty}$ SSMs of the nonlinear system from
formula (\ref{eq:SSM in original coordinates real}) by setting $K,L,M,N,P=0$.
Note that this also includes mixed-mode SSMs and hence already provides
a larger family of invariant manifolds than those defined as SSMs
by \citet{haller16}, who used the invariant manifold results of \citet{cabre03}
(as opposed to the linearization results of \citet{sternberg58})
to construct invariant manifolds tangent to spectral subspaces. The
results of \citet{cabre03} require weaker nonresonance conditions
than those in \eqref{eq:nonresonance conditions-2} but do not yield
fractional or mixed-mode SSMs.

\subsection{Complex form of the SSM parameterization\label{subsec:Complex-form-of-SSM-CDS}}

We have derived the general SSM family (\ref{eq:general_fractional_SSM1-3})
in real coordinates to avoid difficulties with handling complex logarithms
and fractional powers in the derivation. The final form of the parameterization,
however, can readily be converted to a more compact complexified form.
Specifically, using the polar coordinates from eq. (\ref{eq:polar coordinates}\textrm{),
we introduce the complex vectors $z\in\mathbb{C}^{q}$ and $w\in\mathbb{C}^{s}$
by letting 
\begin{align*}
z_{k} & =r_{k}e^{i\phi_{k}},\quad w_{m}=\rho_{m}e^{i\theta_{m}},\qquad k=1,\ldots,q,\quad m=1,\ldots,s.
\end{align*}
}We also restrict the general formula (\ref{eq:general_fractional_SSM1-3})
to SSMs that are symmetric in all the real $u_{j}$ variables (see
Fig. 3 of the main text) to obtain 
\begin{align}
v_{\ell} & =\mathcal{V}_{\ell}(u,z):=\sum_{j=1}^{p}K_{\ell j}\left|u_{j}\right|^{\frac{\kappa_{\ell}}{\lambda_{j}}}+\sum_{k=1}^{q}L_{\ell k}\left|z_{k}\right|^{\frac{\kappa_{\ell}}{\alpha_{k}}},\nonumber \\
w_{m} & =\mathcal{\mathcal{E}}_{m}(u,z):=\left[\sum_{j=1}^{p}O_{mj}\left|u_{j}\right|^{\frac{\beta_{m}}{\lambda_{j}}}+\sum_{k=1}^{q}Q_{mk}\left|z_{k}\right|^{\frac{\beta_{m}}{\alpha_{k}}}\right]e^{i\left(\sum_{j=1}^{p}\frac{\nu_{m}}{\lambda_{j}}\log\left|u_{j}\right|+\sum_{k=1}^{q}\frac{\nu_{m}}{\alpha_{k}}\log\left|z_{k}\right|\right)},\label{eq:complex SSM 1-3}
\end{align}
for all $\ell=1,\ldots,r$, $m=1,\ldots,s$, and with the definition
$O_{mj}=M_{mj}e^{iN_{m}}$ and $Q_{mk}=P_{mk}e^{iN_{m}}$. In these
coordinates, the near-identity, linearizing transformation (\ref{eq:linearizing transformation-2})
becomes
\begin{align}
x & =y+\sum_{\left|\mathbf{k}\right|\geq2}H_{\mathbf{k}}u^{\mathbf{k}_{1}}z^{\mathbf{k}_{2}}\bar{z}^{\mathbf{k}_{3}}v^{\mathbf{k}_{4}}w^{\mathbf{k}_{5}}\bar{w}^{\mathbf{k}_{6}},\label{eq: complex near-identity transformation}\\
 & \,\,\,\,\,\,\mathbf{k}_{1}\in\mathbb{N}^{p},\quad\mathbf{k}_{2},\mathbf{k}_{3}\in\mathbb{N}^{q},\quad\mathbf{k}_{4}\in\mathbb{N}^{r},\quad\mathbf{k}_{5},\mathbf{k}_{6}\in\mathbb{N}^{s},\quad\mathbf{k}=\left(\mathbf{k}_{1},\mathbf{k}_{2},\mathbf{k}_{3},\mathbf{k}_{4},\mathbf{k}_{5},\mathbf{k}_{6}\right).\nonumber 
\end{align}
Substitution of the formulas (\ref{eq:complex SSM 1-3}) into (\ref{eq: complex near-identity transformation})
and denoting the $(v,w)$ coordinate components of the coefficients
$H_{\mathbf{k}}$ by $h_{\mathbf{k}}$ proves statement (i) of the
theorem. 

\subsection{Reduced dynamics on the SSM}

As the SSM family (\ref{eq:general_fractional_SSM1-3}) is a graph
over the $(u,z)$ components of the $x$ variable, the reduced dynamics
on any member of this family can be expressed locally near the origin
in terms of the $(u,z)$ variables. Specifically, substituting formula
(\ref{eq:linearizing transformation-2}) for the SSM family into the
right-hand side of the original nonlinear system (\ref{eq:original system-3})
and keeping only the $(u,z$) components of the system proves statement
(ii) of theorem. 

\subsection{Smoothness and analyticity}

The spectral subspace $E$ of the linearized system (\ref{eq:linearized system})
is also locally the only $C^{\infty}$ invariant graph over $E$,
as all other invariant graphs over $E$ have finite smoothness by
formula (\ref{eq:complex SSM 1-3}). Mapping this manifold under the
inverse of the $C^{\infty}$ linearizing transformation (\ref{eq:linearizing transformation})
gives a unique $C^{\infty}$ invariant graph over $E$ in the original
nonlinear system (\ref{eq:original system-3}). Under condition \emph{
\begin{equation}
\frac{\kappa_{\ell}}{\lambda_{j}}<0,\qquad\frac{\beta_{m}}{\lambda_{j}}<0,\qquad\frac{\kappa_{\ell}}{\alpha_{k}}<0,\qquad\frac{\beta_{m}}{\alpha_{k}}<0,\label{eq:condition for SSM uniqueness in any smoothness class-1}
\end{equation}
}there are no neighboring continuous invariant graphs over the spectral
subspace $E$, as one sees from formula (\ref{eq:complex SSM 1-3}),
and hence there is a unique (primary), class $C^{\infty}$ SSM, $\mathcal{W}\left(E\right)$,
in this case in system (\ref{eq:original system-3}). If at least
one of the inequalities in (\eqref{eq:condition for SSM uniqueness in any smoothness class-1})
are violated, then there are additional invariant graphs over $E$
but all have finite smoothness. Therefore, even in that case, the
$C^{\infty}$ SSM is still unique in its own smoothness class. 

Typical fractional SSMs will have fractional expansions involving
all terms in in eq. (\ref{eq:complex SSM 1-3}) that violate one of
the inequalities in statement (iii) of the Theorem 1. This implies
a generic smoothness class of $C^{\eta}$ for fractional SSMs, with
the positive integer $\eta$ defined in statement (iii) of the Theorem
1. Nongeneric fractional SSMs (that do not include all admissible
terms in their expansions in (\ref{eq:complex SSM 1-3}) can have
smoothness classes higher than $C^{\eta}$, as noted in one of the
remarks after Theorem 1. This concludes the proof of statement (iii)
of the theorem.

\citet{sell85} proves that if the original system (\ref{eq:original system-3})
is $C^{R}$ smooth with respect to a parameter vector $\mu$, then
the transformation is also $C^{R}$ smooth with respect to $\mu$.
This proves statement (iv) of Theorem 1.

\citet{sternberg57} proves the existence of $C^{\rho}$  linearizing
transformations for nonresonant fixed points that are asymptotically
stable in forward or backward time. This implies statement (v) of
Theorem 1, as long as the leading part of the expansions for the SSMs
is below $\rho$, as stated.

Under the assumptions of \citet{sternberg57} and under the added
assumption of analyticity for $f$, Poincaré's analytic linearization
theorem (see, e.g., \citep{arnold83}) implies that the formal Taylor
expansion in eq. (\ref{eq:linearizing transformation}) converges
in a neighborhood of the origin. This proves the statement (vi) about
the convergence of the expansion in formula for the SSMs. 

\section{Appendix D: Proofs of Propositions 1 and 2 \label{sec:Proofs of propositions}}

\subsection{Proof of Proposition 1}

Under the assumptions of Proposition 1, we obtain from the general
formula in statement (i) of Theorem 1 the specific expression 
\begin{align}
v & =\mathcal{V}+\sum_{2\leq\mathrm{k_{1}+\sum_{\ell=1}^{r}k_{4\ell}\frac{\kappa_{\ell}}{\lambda_{1}}}\leq\mathcal{K}}h_{\mathbf{k}}u_{1}^{k_{1}}\mathcal{V}^{\mathbf{k}_{4}}+o\left(\left|u_{1}\right|^{\mathcal{K}}\right),\label{eq:SSM in original coordinates real-1-5-2}\\
 & \,\,\,\,\,\,k_{1}\in\mathbb{N},\quad\mathbf{k}_{4}\in\mathbb{N}^{r},\quad\mathbf{k}=\left(k_{1},\mathbf{k}_{4}\right),\nonumber 
\end{align}
with 
\begin{align}
\mathcal{V}_{\ell}(u,z) & =K_{\ell1}(u_{1})\left|u_{1}\right|^{\frac{\kappa_{\ell}}{\lambda_{1}}}\quad\ell=1,\ldots,r,\label{eq:complex SSM 1-1}
\end{align}
where
\[
K_{\ell j}(u_{j})\equiv0,\quad\mathrm{if}\quad\frac{\kappa_{\ell}}{\lambda_{j}}<0.
\]

Substitution of (\ref{eq:complex SSM 1-1}) into (\ref{eq:SSM in original coordinates real-1-5-2})
and a subsequent re-indexing gives 
\begin{align}
v & =\left(\begin{array}{c}
K_{11}(u_{1})\left|u_{1}\right|^{\frac{\kappa_{1}}{\lambda_{1}}}\\
\vdots\\
K_{r1}(u_{1})\left|u_{1}\right|^{\frac{\kappa_{r}}{\lambda_{1}}}
\end{array}\right)+\sum_{2\leq\mathrm{k_{1}+\sum_{\ell=1}^{r}k_{4\ell}\frac{\kappa_{\ell}}{\lambda_{1}}}\leq\mathcal{K}}h_{\mathbf{k}}u_{1}^{k_{1}}\prod_{\ell=1}^{r}K_{\ell1}^{k_{4\ell}}(u_{1})\left|u_{1}\right|^{\sum_{\ell=1}^{r}k_{4\ell}\frac{\kappa_{\ell}}{\lambda_{1}}}+o\left(\left|u_{1}\right|^{\mathcal{K}}\right)\nonumber \\
 & =\left(\begin{array}{c}
K_{11}(u_{1})\left|u_{1}\right|^{\frac{\kappa_{1}}{\lambda_{1}}}\\
\vdots\\
K_{r1}(u_{1})\left|u_{1}\right|^{\frac{\kappa_{r}}{\lambda_{1}}}
\end{array}\right)+\sum_{2\leq\mathrm{k_{1}+\sum_{\ell=1}^{r}k_{4\ell}\frac{\kappa_{\ell}}{\lambda_{1}}}\leq\mathcal{K}}h_{\mathbf{k}}\prod_{\ell=1}^{r}K_{\ell1}^{k_{4\ell}}(u_{1})u_{1}^{k_{1}}\left|u_{1}\right|^{\sum_{\ell=1}^{r}k_{4\ell}\frac{\kappa_{\ell}}{\lambda_{1}}}+o\left(\left|u_{1}\right|^{\mathcal{K}}\right)\nonumber \\
 & =\sum_{1\leq\mathrm{k_{1}+\sum_{\ell=1}^{r}k_{4\ell}\frac{\kappa_{\ell}}{\lambda_{1}}}\leq\mathcal{K}}C_{\mathbf{k}}(u_{1})u_{1}^{k_{1}}\left|u_{1}\right|^{\sum_{\ell=1}^{r}k_{4\ell}\frac{\kappa_{\ell}}{\lambda_{1}}}+o\left(\left|u_{1}\right|^{\mathcal{K}}\right),\label{eq:2D w_SSM-1-1}
\end{align}
with 
\[
C_{\mathbf{k}\ell}(u_{1})=\left\{ \begin{array}{c}
C_{\mathbf{k}\ell}^{+}\qquad u_{1}>0,\\
\\
C_{\mathbf{k}\ell}^{-}\qquad u_{1}\leq0,
\end{array}\right.\qquad C_{\mathbf{k}\ell}^{\pm}=0:\quad\frac{\kappa_{\ell}}{\lambda_{1}}<0\,\,\,\,\,\,\,\mathrm{or}\,\,\,\,\,\left|\mathbf{k}\right|=1,\,k_{1}=1.
\]
as claimed.

\subsection{Proof of Proposition 2}

Under the assumptions of Proposition 2, we obtain from the general
formula in statement (i) of Theorem 1 the specific expression 
\begin{align}
w & =\mathcal{\mathcal{E}}+\sum_{1\leq\mathrm{k_{2}+k_{3}+\sum_{m=1}^{s}\left(k_{5m}+k_{6m}\right)\frac{\beta_{m}}{\alpha_{1}}}\leq\mathcal{K}}h_{\mathbf{k}}z_{1}^{k_{2}}\bar{z}_{1}^{k_{3}}\mathcal{\mathcal{E}}^{\mathbf{k}_{5}}\bar{\mathcal{\mathcal{E}}}^{\mathbf{k}_{6}}+o\left(\left|z_{1}\right|^{\mathcal{K}}\right),\label{eq:SSM in original coordinates real-1-4}\\
 & \,\,\,\,\,\,k_{2},k_{3}\in\mathbb{N},\quad\mathbf{k}_{5},\mathbf{k}_{6}\in\mathbb{N}^{s},\quad\mathbf{k}=\left(k_{2},k_{3},\mathbf{k}_{5},\mathbf{k}_{6}\right),\nonumber 
\end{align}
 with
\begin{align}
\mathcal{\mathcal{E}}_{m}(z) & =Q_{m1}\left|z_{1}\right|^{\frac{\beta_{m}}{\alpha_{1}}}e^{i\frac{\nu_{m}}{\alpha_{1}}\log\left|z_{1}\right|},\quad m=1,\ldots,s,\label{eq:Epsilon for q=00003D1}
\end{align}
where 
\begin{align}
 & Q_{m1}=0,\quad\mathrm{if}\quad\frac{\beta_{m}}{\alpha_{1}}<0.\label{eq:general_fractional_SSM2-1}
\end{align}
To put the scope of formula (\ref{eq:SSM in original coordinates real-1-4})
in perspective for this case, all manifold-based model-reduction approaches
in the literature have been focusing exclusively on SSMs satisfying
$Q_{mk}=0$ for all $m=1,\ldots,s$. Substitution of (\ref{eq:Epsilon for q=00003D1})
into (\ref{eq:SSM in original coordinates real-1-4}) and a subsequent
re-indexing gives 
\begin{align}
w & =\left(\begin{array}{c}
Q_{11}\left|z\right|^{\frac{\beta_{1}}{\alpha_{1}}}e^{i\frac{\nu_{1}}{\alpha_{1}}\log\left|z\right|}\\
\vdots\\
Q_{s1}\left|z\right|^{\frac{\beta_{s}}{\alpha_{1}}}e^{i\frac{\nu_{s}}{\alpha_{1}}\log\left|z\right|}
\end{array}\right)\nonumber \\
 & +\sum_{2\leq\mathrm{k_{2}+k_{3}+\sum_{m=1}^{s}\left(k_{5m}+k_{6m}\right)\frac{\beta_{m}}{\alpha_{1}}}\leq\mathcal{K}}D_{\mathbf{k}}z_{1}^{k_{2}}\bar{z}_{1}^{k_{3}}\prod_{m=1}^{s}\left|z_{1}\right|^{\left(k_{5m}+k_{6m}\right)\frac{\beta_{m}}{\alpha_{1}}}e^{i\left(k_{5m}-k_{6m}\right)\frac{\nu_{m}}{\alpha_{1}}\log\left|z\right|}+o\left(\left|z_{1}\right|^{\mathcal{K}}\right)\nonumber \\
 & =\sum_{2\leq\mathrm{k_{2}+k_{3}+\sum_{m=1}^{s}\left(k_{5m}+k_{6m}\right)\frac{\beta_{m}}{\alpha_{1}}}\leq\mathcal{K}}D_{\mathbf{k}}z_{1}^{k_{2}}\bar{z}_{1}^{k_{3}}\prod_{m=1}^{s}\left|z_{1}\right|^{\left(k_{5m}+k_{6m}\right)\frac{\beta_{m}}{\alpha_{1}}}e^{i\left(k_{5m}-k_{6m}\right)\frac{\nu_{m}}{\alpha_{1}}\log\left|z\right|}+o\left(\left|z_{1}\right|^{\mathcal{K}}\right)\label{eq:2D w_SSM-1}\\
 & =\sum_{2\leq\mathrm{k_{2}+k_{3}+\sum_{m=1}^{s}\left(k_{5m}+k_{6m}\right)\frac{\beta_{m}}{\alpha_{1}}}\leq\mathcal{K}}D_{\mathbf{k}}z_{1}^{k_{2}}\bar{z}_{1}^{k_{3}}\left|z_{1}\right|^{\sum_{m=1}^{s}\left(k_{5m}+k_{6m}\right)\frac{\beta_{m}}{\alpha_{1}}}e^{i\sum_{m=1}^{s}\left(k_{5m}-k_{6m}\right)\frac{\nu_{m}}{\alpha_{1}}\log\left|z_{1}\right|}+o\left(\left|z_{1}\right|^{\mathcal{K}}\right),\nonumber 
\end{align}
with 
\begin{align*}
D_{\mathbf{k}m} & =0:\quad\frac{\beta_{m}}{\alpha_{1}}<0\,\,\,\,\,\,\mathrm{or}\,\,\,\,\,\,\,\left|\mathbf{k}\right|=1,\,\,k_{2}+k_{3}=1.
\end{align*}
as claimed.

\section{Appendix E: Proof of Theorem 2\label{sec:Proof of main theorem for discrete systems}}

\subsection{Generalized SSMs in the linearized mapping}

To construct the full SSM family arising from the spectral subspace
$E$ in the nonlinear mapping 
\begin{equation}
x(\imath+1)=\mathcal{A}x(\imath)+f\left(x(\imath)\right),\qquad x\in\mathbb{R}^{n},\quad\mathcal{A}\in\mathbb{R}^{n\times n},\quad f=\mathcal{O}\left(\left|x\right|^{2}\right)\in C^{\infty},\qquad1\leq n<\infty,\label{eq:original system-1-1}
\end{equation}
we first describe all invariant graphs over the spectral subspace
$E$ in the linear part 
\begin{equation}
x(\imath+1)=\mathcal{A}x(\imath),\qquad x\in\mathbb{R}^{n},\label{eq:linearized system-2}
\end{equation}
of the original nonlinear system, where 
\begin{equation}
\mathcal{A}=\left(\begin{array}{cccccccc}
A\\
 & B_{1}\\
 &  & \ddots &  &  &  & 0\\
 &  &  & B_{q}\\
 &  &  &  & C\\
 &  &  &  &  & D_{1}\\
 &  & 0 &  &  &  & \ddots\\
 &  &  &  &  &  &  & D_{s}
\end{array}\right),\label{eq:linsystem1-1-1-1}
\end{equation}
where
\begin{align}
A & =\mathrm{diag}\left[\lambda_{1},\ldots,\lambda_{p}\right]\in\mathbb{R}^{p\times p},\quad B_{k}=\left(\begin{array}{cc}
\alpha_{k} & -\omega_{k}\\
\omega_{k} & \alpha_{k}
\end{array}\right)\in\mathbb{R}^{2\times2},\quad k=1,\ldots,q,\nonumber \\
C & =\mathrm{diag}\left[\kappa_{1},\ldots,\kappa_{r}\right]\in\mathbb{R}^{r\times r},\quad D_{m}=\left(\begin{array}{cc}
\beta_{m} & -\nu_{m}\\
\nu_{m} & \beta_{m}
\end{array}\right)\in\mathbb{R}^{2\times2},\quad m=1,\ldots,s.\label{eq:linsystem1-2-1-1}
\end{align}
We also recall that the hyperbolicity condition 
\begin{equation}
\mathrm{spect}\left(\mathcal{A}\right)\cap\left\{ \lambda\in\mathbb{C}:\,\,\left|\lambda\right|=1\right\} =\emptyset,\label{eq:hyperbolicity-2-1}
\end{equation}
the nonresonance conditions 
\begin{equation}
\lambda_{j}\neq\prod_{k=1}^{n}\lambda_{k}^{m_{k}},\quad m_{k}\in\mathbb{N},\quad\sum_{k=1}^{n}m_{k}\geq2,\quad j=1,\ldots,n,\label{eq:nonresonance conditions-1-1}
\end{equation}
and the simplifying assumption 
\begin{equation}
\kappa_{\ell}>0,\quad\ell=1,\ldots,r,\label{eq:positivre real part-1}
\end{equation}
 are assumed to be satisfied.

We first introduce multi-dimensional polar coordinates $(r,\phi)$
and $(\rho,\theta)$ instead of the $(a,b)$ and $(c,d)$ coordinate
pairs by letting\textrm{
\begin{align}
a_{k} & =r_{k}\cos\phi_{k},\quad\,\,\,\,\,b_{k}=r_{k}\sin\phi_{k},\nonumber \\
c_{m} & =\rho_{m}\cos\theta_{m},\quad d_{m}=\rho_{m}\sin\theta_{m}.\label{eq:polar coordinates-2}
\end{align}
}This transform the $\left(a,b,c,d\right)$ component of the linear
system (\ref{eq:linearized system-2}) to the polar form
\begin{align}
r_{k} & (\imath+1)=\sqrt{\alpha_{k}^{2}+\omega_{k}^{2}}r_{k}(\imath),\quad\phi_{k}(\imath+1)=\phi_{k}(\imath)+\arctan\frac{\omega_{k}}{\alpha_{k}},\quad\,\,\,\,\,k=1,\ldots,q,\label{eq:linsystem2+-1}\\
\rho_{m} & (\imath+1)=\sqrt{\beta_{m}^{2}+\nu_{m}^{2}}\rho_{m}(\imath),\quad\theta_{m}(\imath+1)=\theta_{m}(\imath)+\arctan\frac{\nu_{m}}{\beta_{m}},\quad m=1,\ldots,s.\label{eq:linsystem4+-1}
\end{align}

As in the continuous-time case, we can write the $v_{\ell}(\imath+1)$
component of a trajectory contained in an invariant graph over $E$
as

\begin{equation}
v_{\ell}(\imath+1)=F_{\ell}\left(u(\imath+1),r(\imath+1)\right),\quad\ell=1,\ldots,r.\label{eq:discrete invariance eq. 1}
\end{equation}
Using eqs. (\ref{eq:linsystem1-1-1-1}), (\ref{eq:linsystem1-2-1-1}),
(\ref{eq:linsystem2+-1}) and (\ref{eq:linsystem4+-1}), we find that
eq. (\ref{eq:discrete invariance eq. 1}) is equivalent to
\begin{equation}
\kappa_{\ell}F_{\ell}\left(u(\imath),r(\imath)\right)=F_{\ell}\left(\lambda_{1}u_{1}(\imath),\ldots,\lambda_{p}u_{p}(\imath),\sqrt{\alpha_{1}^{2}+\omega_{1}^{2}}r_{1}(\imath),\ldots,\sqrt{\alpha_{q}^{2}+\omega_{q}^{2}}r_{q}(\imath)\right).\label{eq:discrete_SSM_invariance_eq1}
\end{equation}
We find the general solution of the functional equation (\ref{eq:discrete_SSM_invariance_eq1})
for $r_{k}\geq0$ to be
\begin{align}
F_{\ell} & =\sum_{j=1}^{p}\left[K_{\ell j}^{+}H(u_{j})+K_{\ell j}^{-}H(-u_{j})\right]\left|u_{j}\right|^{\log_{\left|\lambda_{j}\right|}\kappa_{\ell}}+\sum_{k=1}^{q}L_{\ell k}r_{k}^{\log_{\sqrt{\alpha_{k}^{2}+\omega_{k}^{2}}}\kappa_{\ell}}\nonumber \\
 & =\sum_{j=1}^{p}\left[K_{\ell j}^{+}H(u_{j})+K_{\ell j}^{-}H(-u_{j})\right]\left|u_{j}\right|^{\frac{\log\kappa_{\ell}}{\log\left|\lambda_{j}\right|}}+\sum_{k=1}^{q}L_{\ell k}\left(a_{k}^{2}+b_{k}^{2}\right)^{\frac{\log\kappa_{\ell}}{\log\left(\alpha_{k}^{2}+\omega_{k}^{2}\right)}},\label{eq:fuv-1}
\end{align}
with the Heaviside function $H(x)$ defined in (\ref{eq:heavyside})
and with the real constants $K_{\ell j}^{\pm}$ and $L_{\ell k}$
subject to the restrictions
\begin{align}
K_{\ell j}^{\pm}=0,\quad\mathrm{if}\quad\frac{\log\kappa_{\ell}}{\log\left|\lambda_{j}\right|}<0,\nonumber \\
L_{\ell k}=0,\quad\mathrm{if}\quad\frac{\log\kappa_{\ell}}{\log\left(\alpha_{k}^{2}+\omega_{k}^{2}\right)}<0.\label{eq:contraints on K and L-1}
\end{align}
We recall that $\kappa_{\ell}$ is positive by assumption (\ref{eq:positivre real part-1}). 

Next, we seek the $\left(\rho_{m},\theta_{m}\right)$ component of
a trajectory of the linearized map (\ref{eq:linearized system-2})
contained in an invariant graph over $E$ as

\[
\rho_{m}(\imath+1)=G_{m}^{\rho}\left(u(\imath+1),r(\imath+1)\right),\quad\theta_{m}(\imath+1)=G_{m}^{\theta}\left(u(\imath+1),r(\imath+1)\right),\quad m=1,\ldots,s,
\]
 or, equivalently, 
\begin{align}
\sqrt{\beta_{m}^{2}+\nu_{m}^{2}}G_{m}^{\rho}\left(u(\imath),r(\imath)\right) & =G_{m}^{\rho}\left(\lambda_{1}u_{1}(\imath),\ldots,\lambda_{p}u_{p}(\imath),\sqrt{\alpha_{1}^{2}+\omega_{1}^{2}}r_{1}(\imath),\ldots,\sqrt{\alpha_{q}^{2}+\omega_{q}^{2}}r_{q}(\imath)\right),\nonumber \\
G_{m}^{\rho}\left(u(\imath),r(\imath)\right)+\arctan\frac{\nu_{m}}{\beta_{m}} & =G_{m}^{\rho}\left(\lambda_{1}u_{1}(\imath),\ldots,\lambda_{p}u_{p}(\imath),\sqrt{\alpha_{1}^{2}+\omega_{1}^{2}}r_{1}(\imath),\ldots,\sqrt{\alpha_{q}^{2}+\omega_{q}^{2}}r_{q}(\imath)\right).\label{eq:discrete_SSM_invariance_eq2}
\end{align}
The general solutions of these functional equations are
\begin{align}
G_{m}^{\rho} & =\sum_{j=1}^{p}\left[M_{mj}^{+}H(u_{j})+M_{mj}^{-}H(-u_{j})\right]\left|u_{j}\right|^{\log_{\left|\lambda_{j}\right|}\sqrt{\beta_{m}^{2}+\nu_{m}^{2}}}+\sum_{k=1}^{q}P_{mk}r_{k}^{\log_{\sqrt{\alpha_{k}^{2}+\omega_{k}^{2}}}\sqrt{\beta_{m}^{2}+\nu_{m}^{2}}}\nonumber \\
 & =\sum_{j=1}^{p}\left[M_{mj}^{+}H(u_{j})+M_{mj}^{-}H(-u_{j})\right]\left|u_{j}\right|^{\frac{\log\sqrt{\beta_{m}^{2}+\nu_{m}^{2}}}{\log\left|\lambda_{j}\right|}}+\sum_{k=1}^{q}P_{mk}\left(a_{k}^{2}+b_{k}^{2}\right)^{\frac{\log\sqrt{\beta_{m}^{2}+\nu_{m}^{2}}}{\log\left(\alpha_{k}^{2}+\omega_{k}^{2}\right)}},\label{eq:grho-1}\\
G_{m}^{\theta} & =N_{m}+\frac{1}{p+q}\arctan\frac{\nu_{m}}{\beta_{m}}\left[\sum_{j=1}^{p}\frac{\log\left|u_{j}\right|}{\log\left|\lambda_{j}\right|}+\sum_{k=1}^{q}\frac{\log r_{k}}{\log\sqrt{\alpha_{k}^{2}+\omega_{k}^{2}}}\right]\nonumber \\
 & =N_{m}+\frac{1}{p+q}\arctan\frac{\nu_{m}}{\beta_{m}}\left[\sum_{j=1}^{p}\frac{\log\left|u_{j}\right|}{\log\left|\lambda_{j}\right|}+\sum_{k=1}^{q}\frac{\log\left(a_{k}^{2}+b_{k}^{2}\right)}{\log\left(\alpha_{k}^{2}+\omega_{k}^{2}\right)}\right].\label{eq:gtheta-1}
\end{align}
 Here $M_{mj}^{\pm}$ and $P_{mk}$ are arbitrary real constants subject
to the constraints
\begin{align}
M_{mj}^{+} & =M_{mj}^{-}=0,\quad\mathrm{if}\quad\frac{\log\sqrt{\beta_{m}^{2}+\nu_{m}^{2}}}{\log\left|\lambda_{j}\right|}<0,\label{eq:contraints on M and P-1-1}
\end{align}
that are needed for the same reason as (\ref{eq:contraints on K and L}).

Using the formulas (\ref{eq:polar coordinates-2}), (\ref{eq:fuv-1}),
(\ref{eq:grho-1}) and (\ref{eq:gtheta-1}), we obtain the following
general form for the full family of invariant graphs over $E$ that
contain the origin:
\begin{align}
v_{\ell} & =\mathcal{V}_{\ell}(u,a,b;K,L)\nonumber \\
 & :=\sum_{j=1}^{p}\left[K_{\ell j}^{+}H(u_{j})+K_{\ell j}^{-}H(-u_{j})\right]\left|u_{j}\right|^{\frac{\log\kappa_{\ell}}{\log\left|\lambda_{j}\right|}}+\sum_{k=1}^{q}L_{\ell k}\left(a_{k}^{2}+b_{k}^{2}\right)^{\frac{\log\kappa_{\ell}}{\log\left(\alpha_{k}^{2}+\omega_{k}^{2}\right)}},\quad\ell=1,\ldots,r,\nonumber \\
\left(\begin{array}{c}
c_{m}\\
d_{m}
\end{array}\right) & =\left(\begin{array}{c}
\mathcal{C}{}_{m}(u,a,b;M,N,P)\\
\mathcal{D}_{m}(u,a,b;M,N,P)
\end{array}\right)\nonumber \\
 & :=\left[\sum_{j=1}^{p}\left[M_{mj}^{+}H(u_{j})+M_{mj}^{-}H(-u_{j})\right]\left|u_{j}\right|^{\frac{\log\sqrt{\beta_{m}^{2}+\nu_{m}^{2}}}{\log\left|\lambda_{j}\right|}}+\sum_{k=1}^{q}P_{mk}\left(a_{k}^{2}+b_{k}^{2}\right)^{\frac{\log\sqrt{\beta_{m}^{2}+\nu_{m}^{2}}}{\log\left(\alpha_{k}^{2}+\omega_{k}^{2}\right)}}\right]\label{eq:general_fractional_SSM1-3-1}\\
 & \,\,\,\,\,\,\times\left(\begin{array}{c}
\cos\left(N_{m}+\frac{1}{p+q}\arctan\frac{\nu_{m}}{\beta_{m}}\left[\sum_{j=1}^{p}\frac{\log\left|u_{j}\right|}{\log\left|\lambda_{j}\right|}+\sum_{k=1}^{q}\frac{\log\left(a_{k}^{2}+b_{k}^{2}\right)}{\log\left(\alpha_{k}^{2}+\omega_{k}^{2}\right)}\right]\right)\\
\sin\left(N_{m}+\frac{1}{p+q}\arctan\frac{\nu_{m}}{\beta_{m}}\left[\sum_{j=1}^{p}\frac{\log\left|u_{j}\right|}{\log\left|\lambda_{j}\right|}+\sum_{k=1}^{q}\frac{\log\left(a_{k}^{2}+b_{k}^{2}\right)}{\log\left(\alpha_{k}^{2}+\omega_{k}^{2}\right)}\right]\right)
\end{array}\right),\quad m=1,\ldots,s,\nonumber 
\end{align}
where 
\begin{align*}
K_{\ell j}^{+} & =K_{\ell j}^{-}=0,\quad\mathrm{if}\quad\frac{\log\kappa_{\ell}}{\log\left|\lambda_{j}\right|}<0,\\
L_{\ell k} & =0,\quad\mathrm{if}\quad\frac{\log\kappa_{\ell}}{\log\left(\alpha_{k}^{2}+\omega_{k}^{2}\right)}<0,\\
M_{mj}^{+}=M_{mj}^{-} & =0,\quad\mathrm{if}\quad\frac{\log\sqrt{\beta_{m}^{2}+\nu_{m}^{2}}}{\log\left|\lambda_{j}\right|}<0.
\end{align*}
For SSMs symmetric in the $u_{j}$ variable, the formulas again simplify
in the same fashion as in the continuous time case.

\subsection{Full family of SSMs in the nonlinear mapping (\ref{eq:general_fractional_SSM1-3-1})}

As in the continuous time case, the $C^{\infty}$ linearization results
of \citet{sternberg58} imply that the full family of invariant graphs
(\ref{eq:general_fractional_SSM1-3-1}) of the linearized mapping
(\ref{eq:linearized system-2}) are smoothly mapped into a similar
family of invariant manifolds of the original nonlinear mapping (\ref{eq:original system-1-1}).
This means that a formal (i.e., not necessarily convergent) Taylor
expansion,
\begin{equation}
x=y+h(y)=y+\sum_{\left|\mathbf{k}\right|\geq2}H_{\mathbf{k}}y^{\mathbf{k}},\qquad\mathbf{k}\in\mathbb{N}^{n},\quad y^{\mathbf{k}}:=y_{1}^{k_{1}}y_{2}^{k_{2}}\cdots y_{n}^{k_{n}},\quad h_{\mathbf{k}}\in\mathbb{R}^{n},\label{eq:linearizing transformation-1}
\end{equation}
transforms the full mapping (\ref{eq:general_fractional_SSM1-3-1})
into its linear part (\ref{eq:linearized system-2}). 

As in the continuous time case, restricting the specific form 
\begin{align*}
x & =y+\sum_{\left|\mathbf{k}\right|\geq2}H_{\mathbf{k}}u^{\mathbf{k}_{1}}a^{\mathbf{k}_{2}}b^{\mathbf{k}_{3}}v^{\mathbf{k}_{4}}c^{\mathbf{k}_{5}}d^{\mathbf{k}_{6}},\\
 & \,\,\,\,\,\,\mathbf{k}_{1}\in\mathbb{N}^{p},\quad\mathbf{k}_{2},\mathbf{k}_{3}\in\mathbb{N}^{q},\quad\mathbf{k}_{4}\in\mathbb{N}^{r},\quad\mathbf{k}_{5},\mathbf{k}_{6}\in\mathbb{N}^{s},\quad\mathbf{k}=\left(\mathbf{k}_{1},\mathbf{k}_{2},\mathbf{k}_{3},\mathbf{k}_{4},\mathbf{k}_{5},\mathbf{k}_{6}\right),
\end{align*}
of the linearizing transformation to the invariant graph family (\ref{eq:general_fractional_SSM1-3-1}),
we obtain the full family of SSMs, $\mathrm{SSM}_{KLMNP}(E)$, of
the nonlinear mapping (\ref{eq:original system-1-1}) in the parametrized
form
\begin{align}
x & =\mathrm{\mathrm{SSM}_{KLMNP}}(u,a,b)\nonumber \\
 & :=\left(u,a,b,\mathcal{V},\mathcal{C},\mathcal{D}\right)^{\mathrm{T}}+\sum_{\left|\mathbf{k}\right|\geq2}H_{\mathbf{k}}u^{\mathbf{k}_{1}}a^{\mathbf{k}_{2}}b^{\mathbf{k}_{3}}\mathcal{V}^{\mathbf{k}_{4}}\mathcal{C}^{\mathbf{k}_{5}}\mathcal{D}^{\mathbf{k}_{6}},\label{eq:SSM in original coordinates real-2}
\end{align}
with the coordinate components of $\mathcal{V}$, $\mathcal{C}$ and
$\mathcal{D}$ defined in (\ref{eq:general_fractional_SSM1-3-1}).
We obtain locally all $C^{\infty}$ SSMs of the nonlinear system from
formula (\ref{eq:SSM in original coordinates real-2}) by setting
$K,L,M,N,P=0$. 

\subsection{Complex form of the SSM parameterization}

The derivation of the complexified form of the SSM family representation
again requires the\textrm{ the complex vectors $z\in\mathbb{C}^{q}$
and $w\in\mathbb{C}^{s}$ defined as
\begin{align*}
z_{k} & =r_{k}e^{i\phi_{k}},\quad w_{m}=\rho_{m}e^{i\theta_{m}},\qquad k=1,\ldots,q,\quad m=1,\ldots,s.
\end{align*}
}We then obtain
\begin{align}
v_{\ell} & =\mathcal{V}_{\ell}(u,z;K,L):=\sum_{j=1}^{p}K_{\ell j}\left(u_{j}\right)\left|u_{j}\right|^{\frac{\log\kappa_{\ell}}{\log\left|\lambda_{j}\right|}}+\sum_{k=1}^{q}L_{\ell k}\left(a_{k}^{2}+b_{k}^{2}\right)^{\frac{\log\kappa_{\ell}}{\log\left(\alpha_{k}^{2}+\omega_{k}^{2}\right)}},\quad\ell=1,\ldots,r\nonumber \\
w_{m} & =\mathcal{\mathcal{E}}_{m}(u,z;O,Q):=\left[\sum_{j=1}^{p}O_{mj}\left(u_{j}\right)\left|u_{j}\right|^{\frac{\log\left(\beta_{m}^{2}+\nu_{m}^{2}\right)}{2\log\left|\lambda_{j}\right|}}+\sum_{k=1}^{q}Q_{mk}\left|z_{k}\right|^{\frac{\log\left(\beta_{m}^{2}+\nu_{m}^{2}\right)}{\log\left(\alpha_{k}^{2}+\omega_{k}^{2}\right)}}\right]\times\label{eq:complex SSM 1-2}\\
 & \,\,\,\,\,\,\,\,\,\,\,\,\,\,\,\,\,\,\,\,\,\,\,\,\,\,\,\,\,\,\,\,\,\,\,\,\,\,\,\,\,\,\,\,\,\,\,\,\times e^{i\left(\frac{1}{p+q}\arctan\frac{\nu_{m}}{\beta_{m}}\left[\sum_{j=1}^{p}\frac{\log\left|u_{j}\right|}{\log\left|\lambda_{j}\right|}+\sum_{k=1}^{q}\frac{2\log\left|z_{k}\right|}{\log\left(\alpha_{k}^{2}+\omega_{k}^{2}\right)}\right]\right)},\nonumber 
\end{align}
for all $\ell=1,\ldots,r$, $m=1,\ldots,s$, and with the definition
$O_{mj}=M_{mj}e^{iN_{m}}$ and $Q_{mk}=P_{mk}e^{iN_{m}}$. In these
coordinates, the general expansion for the SSM family in (\ref{eq:SSM in original coordinates real-2})
becomes
\begin{align}
x & =\mathrm{\mathrm{SSM}_{KLMQ}}(u,z)\nonumber \\
 & :=\left(u,z,\mathcal{V},\mathcal{\mathcal{E}},\bar{\mathcal{\mathcal{E}}}\right)^{\mathrm{T}}+\sum_{\left|\mathbf{k}\right|\geq2}H_{\mathbf{k}}u^{\mathbf{k}_{1}}z^{\mathbf{k}_{2}}\mathcal{V}^{\mathbf{k}_{3}}\left(u,z\right)\mathcal{E}^{\mathbf{k}_{4}}\left(u,z\right),\label{eq:SSM in original coordinates real-1-2}\\
 & \,\,\,\,\,\,\mathbf{k}_{1}\in\mathbb{N}^{p},\quad\mathbf{k}_{2}\in\mathbb{N}^{q},\quad\mathbf{k}_{3}\in\mathbb{N}^{r},\quad\mathbf{k}_{4},\mathbf{k}_{6}\in\mathbb{N}^{s},\quad\mathbf{k}=\left(\mathbf{k}_{1},\mathbf{k}_{2},\mathbf{k}_{3},\mathbf{k}_{4}\right),
\end{align}
which proves statement (i) of the theorem if we denote the $(v,w)$
coordinate components of the coefficients $H_{\mathbf{k}}$ by $h_{\mathbf{k}}$. 

\subsection{Reduced dynamics, smoothness and analyticity}

Statements (ii)-(v) follow from the same reasoning we used in proving
the corresponding statements of Theorem 1. This is because with the
appropriate re-formulation, all linearization results we have used
for continuous dynamical systems are equally valid for discrete dynamical
systems defined by iterated maps.

\section{Appendix F: Proofs of Propositions 3 and 4}

\subsection{Proof of Proposition 3}

Under the assumptions of Proposition 3, we obtain from the general
formula in statement (i) of Theorem 2 the specific expression 
\begin{align}
v & =\mathcal{V}+\sum_{2\leq\mathrm{k_{1}+\sum_{\ell=1}^{r}k_{4\ell}\frac{\log\kappa_{\ell}}{\log\left|\lambda_{1}\right|}}\leq\mathcal{K}}h_{\mathbf{k}}u_{1}^{k_{1}}\mathcal{V}^{\mathbf{k}_{4}}+o\left(\left|u_{1}\right|^{\mathcal{K}}\right),\label{eq:SSM in original coordinates real-1-5-2-1}\\
 & \,\,\,\,\,\,k_{1}\in\mathbb{N},\quad\mathbf{k}_{4}\in\mathbb{N}^{r},\quad\mathbf{k}=\left(k_{1},\mathbf{k}_{4}\right),\nonumber 
\end{align}
with 
\begin{align}
\mathcal{V}_{\ell}(u,z) & =K_{\ell1}(u_{1})\left|u_{1}\right|^{\frac{\log\kappa_{\ell}}{\log\left|\lambda_{1}\right|}}\quad\ell=1,\ldots,r,\label{eq:complex SSM 1-1-2}
\end{align}
where
\[
K_{\ell j}(u_{j})\equiv0,\quad\mathrm{if}\quad\frac{\log\kappa_{\ell}}{\log\left|\lambda_{1}\right|}<0.
\]

Substitution of (\ref{eq:complex SSM 1-1-2}) into (\ref{eq:SSM in original coordinates real-1-5-2-1})
and a subsequent re-indexing gives 
\begin{align}
v & =\left(\begin{array}{c}
K_{11}(u_{1})\left|u_{1}\right|^{\frac{\log\kappa_{1}}{\log\left|\lambda_{1}\right|}}\\
\vdots\\
K_{r1}(u_{1})\left|u_{1}\right|^{\frac{\log\kappa_{r}}{\log\left|\lambda_{1}\right|}}
\end{array}\right)+\sum_{2\leq\mathrm{k_{1}+\sum_{\ell=1}^{r}k_{4\ell}\frac{\kappa_{\ell}}{\lambda_{1}}}\leq\mathcal{K}}h_{\mathbf{k}}u_{1}^{k_{1}}\prod_{\ell=1}^{r}K_{\ell1}^{k_{4\ell}}(u_{1})\left|u_{1}\right|^{\sum_{\ell=1}^{r}k_{4\ell}\frac{\log\kappa_{\ell}}{\log\left|\lambda_{1}\right|}}+o\left(\left|u_{1}\right|^{\mathcal{K}}\right)\nonumber \\
 & =\sum_{1\leq\mathrm{k_{1}+\sum_{\ell=1}^{r}k_{4\ell}\frac{\log\kappa_{\ell}}{\log\left|\lambda_{1}\right|}}\leq\mathcal{K}}C_{\mathbf{k}}(u_{1})u_{1}^{k_{1}}\left|u_{1}\right|^{\sum_{\ell=1}^{r}k_{4\ell}\frac{\log\kappa_{\ell}}{\log\left|\lambda_{1}\right|}}+o\left(\left|u_{1}\right|^{\mathcal{K}}\right),\label{eq:2D w_SSM-1-1-1}
\end{align}
with 
\[
C_{\mathbf{k}\ell}(u_{1})=\left\{ \begin{array}{c}
C_{\mathbf{k}\ell}^{+}\qquad u_{1}>0,\\
\\
C_{\mathbf{k}\ell}^{-}\qquad u_{1}\leq0,
\end{array}\right.\qquad C_{\mathbf{k}\ell}^{\pm}=0:\quad\frac{\log\kappa_{\ell}}{\log\left|\lambda_{1}\right|}<0\,\,\,\,\,\,\,\mathrm{or}\,\,\,\,\,\left|\mathbf{k}\right|=1,\,k_{1}=1.
\]
as claimed.

\subsection{Proof of Proposition 4}

As in our discussion for continuous dynamical systems, we consider
here in more detail SSMs arising from 2D spectral subspaces in underdamped
oscillatory systems ($p=0$, $q=1$, $r=0$, $\frac{\beta_{m}}{2\alpha_{1}}>0$).
We then obtain from the general formula (\ref{eq:general_fractional_SSM1-3-1})
the specific expression 
\begin{align}
\left(\begin{array}{c}
c_{m}(a,b)\\
d_{m}(a,b)
\end{array}\right) & =\left(\begin{array}{c}
\mathcal{C}{}_{m}(a_{1},b_{1};N,P)\\
\mathcal{D}_{m}(a_{1},b_{1};N,P)
\end{array}\right)\nonumber \\
 & =P_{m1}\left(a_{1}^{2}+b_{1}^{2}\right)^{\frac{\log\left(\beta_{m}^{2}+\nu_{m}^{2}\right)}{2\log\left(\alpha_{1}^{2}+\omega_{1}^{2}\right)}}\left(\begin{array}{c}
\cos\left(N_{m}+\arctan\frac{\nu_{m}}{\beta_{m}}\frac{\log\left(a_{1}^{2}+b_{1}^{2}\right)}{\log\left(\alpha_{1}^{2}+\omega_{1}^{2}\right)}\right)\\
\sin\left(N_{m}+\arctan\frac{\nu_{m}}{\beta_{m}}\frac{\log\left(a_{1}^{2}+b_{1}^{2}\right)}{\log\left(\alpha_{1}^{2}+\omega_{1}^{2}\right)}\right)
\end{array}\right),\quad m=1,\ldots,s,\label{eq:case of 2D SSMs-1}
\end{align}
in which we can individually select the $P_{m1}$ and $N_{m}$ coefficients
to be zero or nonzero. 

For the complex form of the SSM family, eq. (\ref{eq:complex SSM 1-2})
gives, with $z\in\mathbb{C}$, the formula
\begin{align}
w_{m} & =\mathcal{\mathcal{E}}_{m}(z;Q):=Q_{m1}\left|z\right|^{\frac{\log\left(\beta_{m}^{2}+\nu_{m}^{2}\right)}{\log\left(\alpha_{1}^{2}+\omega_{1}^{2}\right)}}e^{i\arctan\frac{\nu_{m}}{\beta_{m}}\frac{\log\left|z\right|}{\log\sqrt{\alpha_{1}^{2}+\omega_{1}^{2}}}},\label{eq:complex SSM 1-1-1}
\end{align}
for all $m=1,\ldots,s$. Then the general expression (\ref{eq:SSM in original coordinates real-1-2})
for the full SSM family becomes
\begin{align}
x & =\mathrm{\mathrm{SSM}_{Q}}(z,\bar{z})\nonumber \\
 & :=\left(z,\bar{z},\mathcal{\mathcal{E}},\bar{\mathcal{\mathcal{E}}}\right)^{\mathrm{T}}+\sum_{\sum_{i=1}^{4}\left|k_{i}\right|\geq2}H_{\mathbf{k}}z^{k_{1}}\bar{z}^{k_{2}}\mathcal{E}^{\mathbf{k}_{3}}\left(z\right)\mathcal{\bar{E}}^{\mathbf{k}_{4}}\left(z\right),\label{eq:SSM in original coordinates real-1-1-1}\\
 & \,\,\,\,\,\,k_{1},k_{2}\in\mathbb{N},\quad\mathbf{k}_{3},\mathbf{k}_{4}\in\mathbb{N}^{s},\quad\mathbf{k}=\left(k_{1},k_{2},\mathbf{k}_{3},\mathbf{k}_{4}\right).\nonumber 
\end{align}
Substitution of (\ref{eq:complex SSM 1-1-1}) into eq. (\ref{eq:SSM in original coordinates real-1-1-1})
and a subsequent reindexing then leads to the expression 
\begin{equation}
\left(\begin{array}{c}
z\\
\bar{z}\\
w\\
\bar{w}
\end{array}\right)=\left(\begin{array}{c}
z\\
\bar{z}\\
\mathcal{\mathcal{E}}(z;Q)\\
\bar{\mathcal{\mathcal{E}}}(z;Q)
\end{array}\right)+\sum_{\left|\mathbf{k}\right|\geq2}\left(\begin{array}{c}
A_{\mathbf{k}}\\
B_{\mathbf{k}}\\
\mathbf{C}_{\mathbf{k}}\\
\mathbf{D}_{\mathbf{k}}
\end{array}\right)z^{k_{1}}\bar{z}^{k_{2}}\prod_{m=1}^{s}\left|z\right|^{\left(k_{3m}+k_{4m}\right)\frac{\log\left(\beta_{m}^{2}+\nu_{m}^{2}\right)}{\log\left(\alpha_{1}^{2}+\omega_{1}^{2}\right)}}e^{i\left(k_{3m}-k_{4m}\right)\arctan\frac{\nu_{m}}{\beta_{m}}\frac{\log\left|z\right|}{\log\sqrt{\alpha_{1}^{2}+\omega_{1}^{2}}}},\label{eq:2D SSM in complex coordinates-1}
\end{equation}
with 
\[
\bar{A}_{k_{1}k_{2}k_{3}k_{4}}=B_{k_{2}k_{1}k_{4}k_{3}}\in\mathbb{C},\quad\bar{\mathbf{C}}_{k_{1}k_{2}k_{3}k_{4}}=\mathbf{D}_{k_{2}k_{1}k_{4}k_{3}}\in\mathbb{C}^{s},
\]
\[
k_{1},k_{2}\in\mathbb{N},\quad\mathbf{k}_{3},\mathbf{k}_{4}\in\mathbb{N}^{s},\quad\mathbf{k}=\left(k_{1},k_{2},\mathbf{k}_{3},\mathbf{k}_{4}\right).
\]

Formula (\ref{eq:2D SSM in complex coordinates-1}) gives the embedding
of the manifolds into the space of the original coordinates as
\begin{align}
w_{\mathrm{SSM}}(z,\bar{z}) & =\left(\begin{array}{c}
Q_{11}\left|z\right|^{\frac{\log\left(\beta_{1}^{2}+\nu_{1}^{2}\right)}{\log\left(\alpha_{1}^{2}+\omega_{1}^{2}\right)}}e^{i\arctan\frac{\nu_{1}}{\beta_{1}}\frac{\log\left|z\right|}{\log\sqrt{\alpha_{1}^{2}+\omega_{1}^{2}}}}\\
\vdots\\
Q_{s1}\left|z\right|^{\frac{\log\left(\beta_{s}^{2}+\nu_{s}^{2}\right)}{\log\left(\alpha_{1}^{2}+\omega_{1}^{2}\right)}}e^{i\arctan\frac{\nu_{s}}{\beta_{s}}\frac{\log\left|z\right|}{\log\sqrt{\alpha_{1}^{2}+\omega_{1}^{2}}}}
\end{array}\right)\nonumber \\
 & +\sum_{\left|\mathbf{k}\right|\geq2}\mathbf{C}_{\mathbf{k}}z^{k_{1}}\bar{z}^{k_{2}}\prod_{m=1}^{s}\left|z\right|^{\left(k_{3m}+k_{4m}\right)\frac{\log\left(\beta_{m}^{2}+\nu_{m}^{2}\right)}{\log\left(\alpha_{1}^{2}+\omega_{1}^{2}\right)}}e^{i\left(k_{3m}-k_{4m}\right)\arctan\frac{\nu_{m}}{\beta_{m}}\frac{\log\left|z\right|}{\log\sqrt{\alpha_{1}^{2}+\omega_{1}^{2}}}}\nonumber \\
 & =\sum_{\left|\mathbf{k}\right|\geq1}\mathbf{C}_{\mathbf{k}}z^{k_{1}}\bar{z}^{k_{2}}\left|z\right|^{\sum_{m=1}^{s}\left(k_{3m}+k_{4m}\right)\frac{\log\left(\beta_{m}^{2}+\nu_{m}^{2}\right)}{\log\left(\alpha_{1}^{2}+\omega_{1}^{2}\right)}}e^{i\sum_{m=1}^{s}\left(k_{3m}-k_{4m}\right)\arctan\frac{\nu_{m}}{\beta_{m}}\frac{\log\left|z\right|}{\log\sqrt{\alpha_{1}^{2}+\omega_{1}^{2}}}},\label{eq:2D complex par}
\end{align}
where
\begin{align*}
\mathbf{C}_{\mathbf{k}} & =\mathbf{0},\quad\left|\mathbf{k}\right|=1,\quad k_{1}+k_{2}=1,\\
C_{\mathbf{k}}^{m} & =\mathbf{0}:\quad\left|\mathbf{k}\right|=1,\quad k_{4m}\neq0,\quad m=1,\ldots s.
\end{align*}

\section{Appendix G: Proof of Theorem 3 \label{sec:Normal forms on fractional SSMs}}

Under the complex parametrization (\ref{eq:2D complex par}) of 2D
secondary SSMs of underdamped nonlinear oscillatory systems, Theorem
1 yields the reduced dynamics on the SSM in the general form
\begin{equation}
\dot{z}=\gamma z+\sum_{\left|\mathbf{k}\right|\geq2}P_{\mathbf{k}}(z,\bar{z}),\quad\gamma=\alpha_{1}+i\beta_{1},\quad P_{\mathbf{k}}(z,\bar{z})=E_{\mathbf{k}}z^{k_{1}}\bar{z}^{k_{2}}\left|z\right|^{\sum_{m=1}^{s}\left(k_{3m}+k_{4m}\right)\frac{\beta_{m}}{\alpha_{1}}}e^{i\sum_{m=1}^{s}\left(k_{3m}-k_{4m}\right)\frac{\nu_{m}}{\alpha_{1}}\log\left|z\right|},\label{eq:reduced dynamics on 2D SSM}
\end{equation}
where we have used the variable $z=z_{1}$ for notational simplicity.
We want to perform a normal form transformation to remove as many
terms from this class of generalized monomials as possible. Specifically,
we seek a near-identity transformation defined on the SSM in the form
\begin{equation}
z=\xi+\sum_{\left|\mathbf{k}\right|\geq2}R_{\mathbf{k}}(\xi,\bar{\xi}),\quad R_{\mathbf{k}}(\xi,\bar{\xi})=F_{\mathbf{k}}\xi^{k_{1}}\bar{\xi}^{k_{2}}\left|\xi\right|^{\sum_{m=1}^{s}\left(k_{3m}+k_{4m}\right)\frac{\beta_{m}}{\alpha_{1}}}e^{i\sum_{m=1}^{s}\left(k_{3m}-k_{4m}\right)\frac{\nu_{m}}{\alpha_{1}}\log\left|\xi\right|}\label{eq:fractional monomial space}
\end{equation}
to simplify the dynamics on the SSM to
\[
\dot{\xi}=\gamma\xi+\sum_{\left|\mathbf{k}\right|\geq2}S_{\mathbf{k}}(\xi,\bar{\xi}).
\]
Substituting this transformation into eq. (\ref{eq:reduced dynamics on 2D SSM}),
we obtain
\begin{align*}
\dot{\xi} & =-\partial_{\xi}\sum_{\left|\mathbf{k}\right|\geq2}R_{\mathbf{k}}(\xi,\bar{\xi})\dot{\xi}-\partial_{\bar{\xi}}\sum_{\left|\mathbf{k}\right|\geq2}R_{\mathbf{k}}(\xi,\bar{\xi})\dot{\bar{\xi}}+\gamma\left(\xi+\sum_{\left|\mathbf{k}\right|\geq2}R_{\mathbf{k}}(\xi,\bar{\xi})\right)+\sum_{\left|\mathbf{k}\right|\geq2}E_{\mathbf{k}}P_{\mathbf{k}}(z,\bar{z})\\
 & =\gamma\xi+\gamma\sum_{\left|\mathbf{k}\right|\geq2}R_{\mathbf{k}}(\xi,\bar{\xi})-\partial_{\xi}\sum_{\left|\mathbf{k}\right|\geq2}R_{\mathbf{k}}(\xi,\bar{\xi})\gamma\xi-\partial_{\bar{\xi}}\sum_{\left|\mathbf{k}\right|\geq2}R_{\mathbf{k}}(\xi,\bar{\xi})\bar{\gamma}\bar{\xi}+\sum_{\left|\mathbf{k}\right|\geq2}E_{\mathbf{k}}P_{\mathbf{k}}(z,\bar{z})\\
 & -\partial_{\xi}\sum_{\left|\mathbf{k}\right|\geq2}R_{\mathbf{k}}(\xi,\bar{\xi})\sum_{\left|\mathbf{k}\right|\geq2}S_{\mathbf{k}}(\xi,\bar{\xi})-\partial_{\bar{\xi}}\sum_{\left|\mathbf{k}\right|\geq2}R_{\mathbf{k}}(\xi,\bar{\xi})\sum_{\left|\mathbf{k}\right|\geq2}\bar{S}_{\mathbf{k}}(\xi,\bar{\xi})
\end{align*}
The order-\textbf{$\mathbf{k}$ }fractional monomial on the right-hand
side of this equation drops out if $\mathcal{N}$ can be selected
such a way that 
\begin{equation}
\gamma R_{\mathbf{k}}(\xi,\bar{\xi})-\left[\partial_{\xi}R_{\mathbf{k}}(\xi,\bar{\xi})\gamma\xi+\partial_{\bar{\xi}}R_{\mathbf{k}}(\xi,\bar{\xi})\bar{\gamma}\bar{\xi}\right]=T_{\mathbf{k}}(\xi,\bar{\xi}),\label{eq:homological equation}
\end{equation}
where $T_{\mathbf{k}}(\xi,\bar{\xi})$ only depends on fractional
monomials of order less than $\mathcal{O}\left(\left|\mathbf{k}\right|\right)$.
As a consequence, a member of the family of (\ref{eq:homological equation})
of homological equations can be solved for $R_{\mathbf{k}}$ if the
linear operator 
\begin{equation}
\mathcal{N}\colon R_{\mathbf{k}}\mapsto\gamma R_{\mathbf{k}}-\left[\partial_{\xi}R_{\mathbf{k}}\gamma\xi+\partial_{\bar{\xi}}R_{\mathbf{k}}\bar{\gamma}\bar{\xi}\right]\label{eq:homological operator}
\end{equation}
is invertible on the space $\mathcal{F}^{\mathbf{k}}$ of the fractional
monomial family defined in (\ref{eq:fractional monomial space}). 

First, note that the operator $\mathcal{N}$maps the space $\mathcal{F}^{\mathbf{k}}$
into itself. Second, note that
\begin{align*}
\partial_{\xi}R_{\mathbf{k}}\gamma\xi & =F_{\mathbf{k}}\gamma\xi e^{i\sum_{m=1}^{s}\left(k_{3m}-k_{4m}\right)\frac{\nu_{m}}{\alpha_{1}}\log\left|\xi\right|}\times\\
 & \times\left[k_{1}\xi^{k_{1}-1}\bar{\xi}^{k_{2}}\left|\xi\right|^{\sum_{m=1}^{s}\left(k_{3m}+k_{4m}\right)\frac{\beta_{m}}{\alpha_{1}}}+\sum_{m=1}^{s}\frac{\beta_{m}}{\alpha_{1}}\left(k_{3m}+k_{4m}\right)\xi^{k_{1}}\bar{\xi}^{k_{2}}\left|\xi\right|^{\sum_{m=1}^{s}\left(k_{3m}+k_{4m}\right)\frac{\beta_{m}}{\alpha_{1}}-1}\frac{\bar{\xi}}{2\left|\xi\right|}\right.\\
 & \,\,\,\,\,\,\,\,\,\,\,\left.+\xi^{k_{1}}\bar{\xi}^{k_{2}}\left|\xi\right|^{\sum_{m=1}^{s}\left(k_{3m}+k_{4m}\right)\frac{\beta_{m}}{\alpha_{1}}}i\sum_{m=1}^{s}\left(k_{3m}-k_{4m}\right)\frac{\nu_{m}}{\alpha_{1}}\frac{1}{\left|\xi\right|}\frac{\bar{\xi}}{2\left|\xi\right|}\right]\\
 & =\gamma\left[k_{1}+\sum_{m=1}^{s}\frac{\beta_{m}}{2\alpha_{1}}\left(k_{3m}+k_{4m}\right)+i\sum_{m=1}^{s}\left(k_{3m}-k_{4m}\right)\frac{\nu_{m}}{2\alpha_{1}}\right]R_{\mathbf{k}},\\
\partial_{\bar{\xi}}R_{\mathbf{k}}\bar{\gamma}\bar{\xi} & =\bar{\gamma}\left[k_{2}+\sum_{m=1}^{s}\frac{\beta_{m}}{2\alpha_{1}}\left(k_{3m}+k_{4m}\right)+i\sum_{m=1}^{s}\left(k_{3m}-k_{4m}\right)\frac{\nu_{m}}{2\alpha_{1}}\right]R_{\mathbf{k}}.
\end{align*}
This implies that all monomials $R_{\mathbf{k}}(\xi,\bar{\xi})$ are
eigenfunctions of the operator $\mathcal{N}$ defined in (\ref{eq:homological operator}).
Specifically, using the extended normal form introduced in the data-driven
SSM construction by \citep{cenedese22a} (which approximates the eigenvalue
pair along the SSM as $\gamma=\pm i\omega$ to maximize the domain
of definition of the normal from transformation), we obtain
\begin{align*}
\mathcal{N}R_{\mathbf{k}} & =\left(i\omega-\left[k_{1}i\omega-k_{2}i\omega\right]\right)R_{\mathbf{k}}=i\omega\left[k_{2}+1-k_{1}\right]R_{\mathbf{k}}.
\end{align*}
When the eigenvalue $i\omega\left[k_{2}+1-k_{1}\right]$ of $\mathcal{N}$
vanishes, the term $T_{\mathbf{k}}(\xi,\bar{\xi})$ in the homological
equation (\ref{eq:homological equation}) is called resonant and will
have to stay in the normal form after the transformation (\ref{eq:fractional monomial space})
is carried out. The general form of all such unremovable, resonant
terms in the normal form is, therefore, 
\begin{equation}
T_{\mathbf{k}}^{res}(\xi,\bar{\xi})=\xi^{k_{1}}\bar{\xi}^{k_{1}-1}\left|\xi\right|^{\sum_{m=1}^{s}\left(k_{3m}+k_{4m}\right)\frac{\beta_{m}}{\alpha_{1}}}e^{i\sum_{m=1}^{s}\left(k_{3m}-k_{4m}\right)\frac{\nu_{m}}{\alpha_{1}}\log\left|\xi\right|},\label{eq:resonant fractional terms-1}
\end{equation}
which includes, as a subclass, the classic normal form terms for $k_{3m}=k_{4m}=0$.
Note that the lowest-order resonant fractional terms in the normal
form arise for $k_{1}=k_{3m}=1$, $k_{4m}=0$ and $k_{1}=k_{4m}=1$,
$k_{3m}=0$ yielding the fractional normal terms $\xi\left|\xi\right|^{\frac{\beta_{m}}{2\alpha_{1}}}e^{i\frac{\nu_{m}}{\alpha_{1}}\log|\xi|}$
and $\xi\left|\xi\right|^{\frac{\beta_{m}}{2\alpha_{1}}}e^{-i\frac{\nu_{m}}{\alpha_{1}}\log|\xi|}$
in the complex normal form ODE for $\xi$ and $\bar{\xi}\left|\xi\right|^{\frac{\beta_{m}}{2\alpha_{1}}}e^{-i\frac{\nu_{m}}{\alpha_{1}}\log|\xi|}$
and $\bar{\xi}\left|\xi\right|^{\frac{\beta_{m}}{2\alpha_{1}}}e^{i\frac{\nu_{m}}{\alpha_{1}}\log|\xi|}$
in the complex normal form ODE for $\bar{\xi}$ on the SSM. Up to
a given polynomial order $\mathcal{K}$ of normalization, only fractional
terms with multi-indices $\mathbf{k}$ satisfying
\[
2k_{1}-1+\sum_{m=1}^{s}\left(k_{3m}+k_{4m}\right)\frac{\beta_{m}}{\alpha_{1}}\leq\mathcal{K}
\]
will appear in the normal form. The smoothness of the normal form
transformation follows directly from the form of the fractional exponents
in eq. (\ref{eq:fractional monomial space}), which concludes the
proof of statement (i) of the theorem.

To prove statement (ii) of the theorem, let us consider SSMs constructed
over a 2D spectral subspace of a 4D dynamical system ($n=4$) with
two pairs of complex conjugate eigenvalues at the origin ($p=r=0$
and $q=s=1$). Using the notation of Section \ref{subsec:Complex-form-of-SSM-CDS},
the complexified dynamical system can be written as 

\begin{equation}
\frac{d}{dt}\left(\begin{array}{c}
z\\
\bar{z}\\
w\\
\bar{w}
\end{array}\right)=\left(\begin{array}{cccc}
\alpha+i\omega & 0 & 0 & 0\\
0 & \alpha-i\omega & 0 & 0\\
0 & 0 & \beta+i\nu & 0\\
0 & 0 & 0 & \beta-i\nu
\end{array}\right)\left(\begin{array}{c}
z\\
z\\
w\\
\bar{w}
\end{array}\right)+\mathbf{f}(z,\bar{z},w,\bar{w}),\label{eq:normalform_nonlineareq}
\end{equation}
where the nonlinearity $\mathbf{f}(z,\bar{z},w,\bar{w})$ is at least
quadratic in its arguments. (We have dropped the subscripts of $z_{1}$
and $w_{1}$ to simplify our notation.) By Proposition 2, the full
SSM family over the 2D subspace 
\[
E=\left\{ (z,\bar{z},w,\bar{w})\in\mathbb{C}^{4}\colon w=\bar{w}=0\right\} 
\]
can be described as
\begin{equation}
w\left(z,\bar{z}\right)=\sum_{\left|\mathbf{k}\right|=1}^{\infty}C_{\mathbf{k}}z^{k_{1}}\bar{z}^{k_{2}}|z|^{(k_{3}+k_{4})\frac{\beta}{\alpha}}e^{i\frac{\nu}{\alpha}\log|z|(k_{3}-k_{4})}.\label{eq:complex SSM}
\end{equation}
We will compute the normal form of the SSM-reduced dynamics up to
cubic orders ($\mathcal{K}=3$), which includes all terms in (\ref{eq:complex SSM})
satisfying
\[
k_{1}+k_{2}+(k_{3}+k_{4})\frac{\beta}{\alpha}\leq3.
\]
Note that the exponential term in (\ref{eq:complex SSM}) is of unit
norm and thus does not contribute to the magnitude of $w.$ 

Under the assumption \emph{
\begin{equation}
\frac{1}{2}<\frac{\beta_{1}}{\alpha_{1}}<2\label{eq:beta over alpha assumption-1-1}
\end{equation}
}of Theorem 3, the relevant terms in the expansion (\ref{eq:complex SSM})
of $w$ are
\begin{align*}
w\left(z,\bar{z}\right) & =C_{1100}z\bar{z}+C_{2000}z^{2}+C_{0200}\bar{z}^{2}+C_{2100}z^{2}\bar{z}+C_{1200}z\bar{z}^{2}+C_{3000}z^{3}+C_{0300}\bar{z}^{3}\\
 & +C_{0010}|z|^{\frac{\beta}{\alpha}}e^{i\frac{\nu}{\alpha}\log|z|}+C_{0001}|z|^{\frac{\beta}{\alpha}}e^{-i\frac{\nu}{\alpha}\log|z|}+C_{1010}z|z|^{\frac{\beta}{\alpha}}e^{i\frac{\nu}{\alpha}\log|z|}+C_{1001}z|z|^{\frac{\beta}{\alpha}}e^{-i\frac{\nu}{\alpha}\log|z|}\\
 & +C_{0110}\bar{z}|z|^{\frac{\beta}{\alpha}}e^{i\frac{\nu}{\alpha}\log|z|}+C_{0101}\bar{z}|z|e^{-i\frac{\nu}{\alpha}\log|z|}+C_{0011}|z|^{2\frac{\beta}{\alpha}}+C_{0020}|z|e^{2i\frac{\nu}{\alpha}\log|z|}\\
 & +C_{0002}|z|^{2\frac{\beta}{\alpha}}e^{-2i\frac{\nu}{\alpha}\log|z|}+o\left(\left|z\right|^{3}\right).
\end{align*}
We obtain the reduced dynamics on the SSM by restricting \eqref{eq:normalform_nonlineareq}
to (\ref{eq:complex SSM}), which gives 
\begin{align*}
\dot{z} & =(\alpha+i\omega)z+az^{2}\bar{z}+bz|z|^{\frac{\beta}{\alpha}}e^{i\theta(|z|)}+cz|z|^{\frac{\beta}{\alpha}}e^{-i\theta(|z|)}+dz|z|^{\frac{2\beta}{\alpha}}+ez|z|^{\frac{2\beta}{\alpha}}e^{i\theta(|z|)}+fz|z|^{\frac{2\beta}{\alpha}}e^{-i\theta(|z|)},\\
\dot{\bar{z}} & =(\alpha-i\omega)\bar{z}+\bar{a}\bar{z}^{2}z+\bar{b}\bar{z}|z|^{\frac{\beta}{\alpha}}e^{i\theta(|z|)}+\bar{c}\bar{z}|z|^{\frac{\beta}{\alpha}}e^{-i\theta(|z|)}+\bar{d}\bar{z}|z|^{\frac{2\beta}{\alpha}}+\bar{e}\bar{z}|z|^{\frac{2\beta}{\alpha}}e^{2i\theta(|z|)}+\bar{f}\bar{z}|z|^{\frac{2\beta}{\alpha}}e^{-2i\theta(|z|)}.
\end{align*}
Introducing polar coordinates via $z=re^{i\phi}$, we obtain we obtain
the real normal form
\begin{align*}
\dot{r} & =\alpha r+Ar^{3}+r^{\frac{\beta}{\alpha}+1}\left[A_{1}\cos\theta+A_{2}\sin\theta\right]+r^{\frac{2\beta}{\alpha}+1}\left[A_{3}+A_{4}\cos2\theta+A_{5}\sin2\theta\right],\\
\dot{\phi} & =\omega+Br^{2}+r^{\frac{\beta}{\alpha}}\left[B_{1}\cos\theta+B_{2}\sin\theta\right]+r^{\frac{2\beta}{\alpha}}\left[B_{3}+B_{4}\cos2\theta+B_{5}\sin2\theta\right],
\end{align*}
or, equivalently, by the definition of $\theta(|z|)$, 
\begin{align}
\dot{r} & =\alpha r+Ar^{3}+P_{1}r^{\frac{\beta}{\alpha}+1}\sin2\left[Q+\frac{\nu}{\alpha}\log r\right]+P_{2}r^{\frac{2\beta}{\alpha}+1}+P_{3}r^{\frac{2\beta}{\alpha}+1}\sin2\left[Q+\frac{\nu}{\alpha}\log r\right],\nonumber \\
\dot{\phi} & =\omega+Br^{2}+R_{1}r^{\frac{\beta}{\alpha}}\sin2\left[Q+\frac{\nu}{\alpha}\log r\right]+R_{2}r^{\frac{2\beta}{\alpha}}+R_{3}r^{\frac{2\beta}{\alpha}}\sin2\left[Q+\frac{\nu}{\alpha}\log r\right],\label{eq:final polar normal form}
\end{align}
for appropriately redefined constants that depend on both the original
full system and the specific SSM chosen. This completes the proof
of statement (ii) of the theorem.

\bibliographystyle{plainnat}

\end{document}